\documentclass[10pt]{article}

\usepackage[unicode]{hyperref}
\usepackage{amsmath}
\usepackage{graphicx}								% to make \includegraphics work
\usepackage{subfig}
\usepackage{amsfonts} 							% to make \mathbb work
\usepackage{dsfont}
\usepackage{bigints} 								% to use big integrals 
\usepackage{algorithm}              % algorithm environment
\usepackage{algorithmic}
\usepackage{booktabs} 							% correct spacing in tables
\usepackage{multirow}								% multirows in tables
\usepackage{setspace}
\usepackage{mathrsfs}
\usepackage{amssymb}
\usepackage{tikz}
\usepackage{color}
\usepackage{xcolor}
\usepackage{mathtools}
\usepackage{xfrac}

\usepackage{framed}
\usepackage{afterpage}
\usepackage{capt-of}
\usepackage{stackengine}

\usepackage{geometry}
\usepackage{adjustbox}
\usepackage{sidecap}

\usepackage{tikz}	
\usetikzlibrary{arrows,chains,matrix,positioning,scopes}

\makeatletter
\tikzset{join/.code=\tikzset{after node path={%
\ifx\tikzchainprevious\pgfutil@empty\else(\tikzchainprevious)% 
edge[every join]#1(\tikzchaincurrent)\fi}}}

\makeatother
\tikzset{>=stealth',every on chain/.append style={join}, every join/.style={->}}
\tikzstyle{labeled}=[execute at begin node=$\scriptstyle,   execute at end node=$]

\newcommand*\rfrac[2]{{}^{#1}\!/_{#2}}
\newcommand{\minus}{\scalebox{0.5}[1.0]{\( - \)}}

\newcommand {\Int}   {\int\limits}
\newcommand {\Sum}   {\sum\limits}

\newcommand {\IntO}  {\Int_\Omega}

\newcommand {\Rd}    {{\mathds{R}}^d}

\newcommand {\Rtwo}  {{\mathds{R}}^2}
\newcommand {\Rthree}  {{\mathds{R}}^3}

\newcommand {\Ieff}  {I_{\rm eff}}

\def \NormA#1  {{\mid\!\mid\!\mid #1 \mid\!\mid\!\mid}^2_A }   % ||...||
\def \NormAinverse#1  { {\mid\!\mid\!\mid #1 \mid\!\mid\!\mid}^2_{A^{-1}} }   % ||...||_*
\def \Normenergy#1  {\mid\!\mid\!\mid #1 \mid\!\mid\!\mid}   % |||...|||

\def \NormQT#1 {{\mid\!\mid\!\mid #1 \mid\!\mid\!\mid}^2_{Q_T}}   % ||...||
\def \NormQO#1 {{\mid\!\mid\!\mid #1 \mid\!\mid\!\mid}^2_{Q^0}}   % ||...||
\def \NormQk#1 {{\mid\!\mid\!\mid #1 \mid\!\mid\!\mid}^2_{Q^k}}   % ||...||
\def \Normf#1  {\Big \lceil #1 \Big\rceil_\Omega}
\def \dvrg     {\mathrm{div}}

\def \dx       {\mathrm{\:d}x}

\def \ds       {\mathrm{\:d}s}

\def \Normt#1  {\mid\!\mid\!\mid #1 \mid\!\mid\!\mid}   % |||...|||

\def\L#1{L^{#1}}
\def\H#1{H^{#1}}

\def\tildeH#1{\widetilde{H}^{#1}}

\def\H#1{H^{#1}}

\newcommand {\CPTrefpithree}  {C^{\, \mathrm{P}}_{\widehat{\T}, {\rfrac{\pi}{3}}}}
\newcommand {\CPTrefpitwo}  {C^{\, \mathrm{P}}_{\widehat{\T}, {\rfrac{\pi}{2}}}}
\newcommand {\CPTrefpifour}  {C^{\, \mathrm{P}}_{\widehat{\T}, {\rfrac{\pi}{4}}}}

\newcommand {\cpithree}       {\chi^{\mathrm{P}}_{{\rfrac{\pi}{3}}}}
\newcommand {\cpitwo}       {\chi^{\mathrm{P}}_{{\rfrac{\pi}{2}}}}
\newcommand {\cpifour}       {\chi^{\mathrm{P}}_{{\rfrac{\pi}{4}}}}

\newcommand {\mupitwo}      {\mu_{\rfrac{\pi}{2}}}
\newcommand {\mupithree}      {\mu_{\rfrac{\pi}{3}}}
\newcommand {\mupifour}      {\mu_{\rfrac{\pi}{4}}}

\newcommand {\CPoincare}  {C^{\mathrm P}_{\Omega}}
\newcommand {\CPGamma}    {C^{\mathrm{P}}_{\Gamma}}
\newcommand {\CtrGamma}   {C^{\mathrm{Tr}}_{\Gamma}}
\newcommand {\CP}      		{\overline{C}^{\, \mathrm{P}}_{\Gamma}}
\newcommand {\CG}      		{\overline{C}^{\, \mathrm{Tr}}_{\Gamma}}

\newcommand {\CPtetr}     {\widetilde{C}^{\, \mathrm{P}}_{\Gamma}}
\newcommand {\CGtetr}     {\widetilde{C}^{\, \mathrm{Tr}}_{\Gamma}}

\newcommand {\CPTrefleg}  {C^{\, \mathrm{P}}_{\widehat{\Gamma}, \rfrac{\pi}{2}}}
\newcommand {\CPTrefhyp}  {C^{\, \mathrm{P}}_{\widehat{\Gamma}, \rfrac{\pi}{4}}}
\newcommand {\CGTrefleg}  {C^{\, \mathrm{Tr}}_{\widehat{\Gamma}, \rfrac{\pi}{2}}}
\newcommand {\CGTrefhyp}  {C^{\, \mathrm{Tr}}_{\widehat{\Gamma}, \rfrac{\pi}{4}}}

\newcommand {\CPThatalphahat} {C^{\, \mathrm{P}}_{\widehat{\Gamma}, \rfrac{\pi}{2}, \hat{\alpha}}}
\newcommand {\CGThatalphahat} {C^{\, \mathrm{Tr}}_{\widehat{\Gamma}, \rfrac{\pi}{2}, \hat{\alpha}}}

\newcommand {\CPTapproxhatalphahat} {C^{\, \mathrm{P}, M}_{\widehat{\Gamma}, \rfrac{\pi}{2}, \hat{\alpha}}}
\newcommand {\CGTapproxhatalphahat} {C^{\, \mathrm{Tr}, M}_{\widehat{\Gamma}, \rfrac{\pi}{2}, \hat{\alpha}}}

\newcommand {\CPTleg}     {\overline{C}^{\, \mathrm{P}}_{\rfrac{\pi}{2}}}
\newcommand {\CPThyp}     {\overline{C}^{\, \mathrm{P}}_{\rfrac{\pi}{4}}}
\newcommand {\CGTleg}     {\overline{C}^{\, \mathrm{Tr}}_{\rfrac{\pi}{2}}}
\newcommand {\CGThyp}     {\overline{C}^{\, \mathrm{Tr}}_{\rfrac{\pi}{4}}}
\newcommand {\CPLS}       {\overline{C}^{\rm P, \triangle}_{\T}}
\newcommand {\CPlower}    {\underline{C}^{\rm P}_{\T}}
\newcommand {\CPupper}    {\overline{C}^{\rm P, \oplus}_{\T}}
\newcommand {\CPMR}       {\overline{C}^{\rm P}_{\T}}

\newcommand {\cpleg} {\gamma^{\mathrm{P}}_{\rfrac{\pi}{2}}}
\newcommand {\cphyp} {\gamma^{\mathrm{P}}_{\rfrac{\pi}{4}}}
\newcommand {\cgleg} {\gamma^{\mathrm{Tr}}_{\rfrac{\pi}{2}}}
\newcommand {\cghyp} {\gamma^{\mathrm{Tr}}_{\rfrac{\pi}{4}}}

\newcommand {\cpalphahat} {\gamma^{\mathrm{P}}_{\rfrac{\pi}{2}, \hat{\alpha}}}
\newcommand {\cgalphahat} {\gamma^{\mathrm{Tr}}_{\rfrac{\pi}{2}, \hat{\alpha}}}

\newcommand {\Tref} {\widehat{{T}}}
\newcommand {\Gref} {\widehat{\Gamma}}
\newcommand {\T} {{T}}

\newcommand {\approxCPT} {\underline{C}^{M, \mathrm{P}}_{\Gamma}}

\newcommand {\approxCGT} {\underline{C}^{M, \mathrm{Tr}}_{\Gamma}}

\newcommand {\approxC} 	 {\underline{C}^{M, \mathrm{P}}_{\T}}

\newcommand {\CPT} {C^{\mathrm{P}}_{\Gamma}}
\newcommand {\CGT} {C^{\mathrm{Tr}}_{\Gamma}}

\newcommand {\muleg} {\mu_{\rfrac{\pi}{2}}}
\newcommand {\muhyp} {\mu_{\rfrac{\pi}{4}}}

\newcommand {\diam} {\mathrm{diam}}

\def \mean#1{\left\{\!\left|#1\right|\!\right\}}
\def \Mean#1#2 {{ \Big \{ #1 \Big\} }_{#2}}
\def \smallMean#1#2 {{ \big \{  #1 \big\} }_{#2}}

\def\ProofBegin{\noindent{\bf Proof:} \:}
\def\ProofEnd{{\hfill $\square$}}

{\bf}{\it}
\newtheorem{lemma}{Lemma}{\bf}{\it}
{\bf}{\it}
\DeclareFixedFont{\ttb}{T1}{txtt}{bx}{n}{12} % for bold
\DeclareFixedFont{\ttm}{T1}{txtt}{m}{n}{12}  % for normal

\usepackage{color}
\definecolor{deepblue}{rgb}{0,0,0.5}
\definecolor{deepred}{rgb}{0.6,0,0}
\definecolor{deepgreen}{rgb}{0,0.5,0}
\usepackage{listings} 

\newcommand\pythonstyle{
\lstset{
		language=Python,
		basicstyle=\footnotesize\ttfamily,
		otherkeywords={self},             % Add keywords here
		keywordstyle=\bfseries\color{deepblue}\bfseries,
		commentstyle=\itshape\color{purple!40!black},
		stringstyle=\color{deepgreen},
		frame=single,                         % Any extra options here
		showstringspaces=false,            % 
		belowcaptionskip=.75\baselineskip
}}
\lstnewenvironment{python}[1][]
{
\pythonstyle
\lstset{#1}
}
{}

\newcommand\pythoninline[1]{{\pythonstyle\lstinline!#1!}}

\newcommand\myatop[2]{\genfrac{}{}{0pt}{}{#1}{#2}}

\title{Sharp bounds of constants in Poincar\'{e}-type inequalities
for simplicial domains}
\author{S. Matculevich and S. Repin\\\\
Department of Mathematical Information Technology, 
University of Jyv\"askyl\"a\\
FIN-40100 Jyv\"askyl\"a, FINLAND\\
e-mails: svetlana.v.matculevich@jyu.fi, sergey.repin@jyu.fi\\\\
St. Petersburg Dept. of V.A. Steklov Institute of Mathematics of RAS\\
St. Petersburg, Russia}
\begin{document}
	
\maketitle 

%----------------------------------------------------------------------------------%
	
\begin{abstract}	

The paper is concerned with sharp estimates of constants in classical Poincar\'{e}
inequalities and Poincar\'{e}-type
inequalities for functions having zero mean value in a simplicial domain or on a part
of its boundary. These estimates are important for quantitative analysis of problems
generated by differential equations, where numerical approximations are typically
constructed with the help of simplicial meshes.
We suggest easily computable relations that provide sharp bounds of the respective
constants and compare these results with analytical estimates (if they are known). 
In the last
section, we present an example that shows possible applications of the results
and derive a computable majorant of the
difference between the exact solution of a boundary value problem and an arbitrary finite
dimensional approximation computed on a simplicial mesh, which uses above mentioned 
constants. 

\end{abstract}

%=====================================================================================%
%=====================================================================================%

\section{Introduction}

%--------------------------------------------------------------------------------------%
Let $T$ be a bounded domain in $\Rd$ ($d \geq 2$) with Lipschitz
boundary $\partial T$. It is well known that the Poincar\'e inequality
(\cite{Poincare1890,Poincare1894})
\begin{equation}
    \| w \|_ {T} \leq {C^{\mathrm P}_{T}} \, \| \nabla w \|_{T}
    \label{eq:classical-poincare-constant}
\end{equation}
holds for any
\begin{equation*}
    w\in\tildeH{1}(T)
    := \Big\{ w \in H^1(T)\,  \big | \, \mean{ w }_{T} = 0 \; \Big\},
\end{equation*}
where $\| w \|_ {T}$ denotes the norm in $\L{2}(T)$,
$\mean{w}_{T}: = \tfrac{1}{|T|} \int_{T} w \, \mathrm{d} x$ is the mean
value of $w$ over $T$, and $|T|$ is the Lebesgue measure of $T$. The constant
${C^{\mathrm P}_{T}}$ depends only on $T$ and $d$.
%--------------------------------------------------------------------------------------%

Poincar\'e-type inequalities also hold for 
\begin{equation*}
	w \in \tildeH{1}(T, \Gamma) :=
	\Big\{ w \in H^1(T)\, \big | \, \mean{w}_{\Gamma} = 0 \;\Big\},
	%\label{eq:space-with-boundary-mean}
\end{equation*}
where $\Gamma$ is a measurable part of $\partial T$ such that
$\mathrm{meas}_{d - 1} \Gamma > 0$ (in particular, $\Gamma$ may coincide with the whole
boundary).
%--------------------------------------------------------------------------------------%
For any $w \in \widetilde{H}^1(T, \Gamma)$, we have two inequalities similar to
\eqref{eq:classical-poincare-constant}. The first one
\begin{equation}
    \|w\|_{T} \, \leq \CPGamma \|\nabla w\|_{T}
    \label{eq:Comega}
\end{equation}
is another form of the Poincar\'e inequality \eqref{eq:classical-poincare-constant},
which is stated for a different set of functions and contains a different constant, 
i.e. ${C^{\mathrm P}_{T}} \leq \CPGamma$. 
The constant $\CPGamma$ is associated with the minimal positive eigenvalue of the 
problem %for any $u \in \H{1} (T, \Gamma)$
\begin{equation}
-\Delta u = \lambda u \;\; {\rm in} \;\; T; \quad 
\partial_n u = \lambda \mean{u}_{T} \;\; {\rm on} \;\; \Gamma; \quad 
\partial_n u = 0 \;\; {\rm on} \;\; \partial T \backslash \Gamma; \quad
\forall u \in \tildeH{1} (T, \Gamma).
\label{eq:eigenvalue-problem-cp}
\end{equation}
We note that inequalities of this type arose in finite element analysis many years 
ago (see, e.g., \cite{BabuskaAziz1976}), where \eqref{eq:Comega} was considered for
simplexes in $\Rtwo$.
The
second inequality
\begin{equation}
    \|w\|_{\Gamma} \, \leq \CtrGamma    \|\nabla w\|_{T}
    \label{eq:Cgamma}
\end{equation}
estimates the trace of $w \in \widetilde{H}^1(T, \Gamma)$ on $\Gamma$.
It is associated with the minimal nonzero eigenvalue of the problem 
\begin{equation}
-\Delta u = 0 \;\; {\rm in} \;\; T; \quad 
\partial_n u = \lambda u \;\; {\rm on} \;\; \Gamma; \quad 
\partial_n u = 0 \;\; {\rm on} \;\; \partial T \backslash \Gamma; \quad
\forall u \in \tildeH{1} (T, \Gamma).
\label{eq:eigenvalue-problem-ctr}
\end{equation}
%
%--------------------------------------------------------------------------------------%
The problem \eqref{eq:eigenvalue-problem-ctr} is a special case of the Steklov problem 
\cite{Stekloff1902}, where the spectral parameter appears in the boundary condition.  
Sometimes \eqref{eq:eigenvalue-problem-ctr} is associated with the so-called 
{\em sloshing problem}, which describes oscillations of a fluid in a container. 
Eigenvalues and eigenfunctions of the sloshing problem have been studied in 
\cite{FoxKuttler1983,BanuelosKulczyckiPolterovichSiudeja2010,KozlovKuznetsov2004,KozlovKuznetsovMotygin2004,KuznetsovKulczyckiKwasnickiNazarovPoborchiPolterovichSiudeja2014, 
ArxivGirouardPolterovich}
and some other papers cited therein.

%--------------------------------------------------------------------------------------%
Exact values of $\CPGamma$, $\CtrGamma$, and ${C^{\mathrm P}_{T}}$ are important from 
both analytical and computational points of view.
Poincar\'{e}-type inequalities are often used in analysis of nonconforming
approximations (e.g., discontinuous Galerkin or mortar methods), domain decomposition
methods (see, e.g., \cite{Klawonnatall2008,Dohrmann2008} and \cite{ToselliWidlund2005}), 
%analysis of problems described in terms of vector valued functions (see, e.g., 
%\cite{Fuchs2011,Pauly2015}), 
a posteriori estimates \cite{RepinBoundaryMeanTrace2015}, and other applications related 
to quantitative analysis of partial differential equations. 
Analysis of interpolation constants and their estimates for 
piecewise constant and linear interpolations over triangular finite elements can be 
found in \cite{XuefengOishi2013} and literature cited therein. Finally, we note that 
\cite{CarstensenGedicke2014} introduces a method of computing lower bounds 
for the eigenvalues of the Laplace operator based on nonconforming (Crouzeix-Raviart) 
approximations. This method yields guaranteed upper bounds of the constant in the 
Friedrichs' inequality.
%Therefore, exact values of respective constants (or sharp and guaranteed bounds of them) are interesting from both analytical and computational points of view.

%--------------------------------------------------------------------------------------%
It is known (see \cite{PayneWeinberger1960}) that for convex domains 
$${C^{\mathrm P}_{T}} \leq \tfrac{\diam (T)}{\pi}.$$
For triangles this estimate was improved
in \cite{LaugesenSiudeja2010} to
$${C^{\mathrm P}_{T}} \leq \tfrac{\diam (T)}{j_{1, 1}},$$ 
where $j_{1, 1} \approx 3.8317$ is the smallest positive
root of the Bessel function $J_1$. Moreover, for isosceles triangles from  
\cite{Bandle1980, LaugesenSiudeja2010} it follows that 
\begin{equation}
    {C^{\mathrm P}_{T}} \leq \CPLS := \diam (T) \cdot \,
    \begin{cases}
    \tfrac{1}{j_{1, 1}} &  \alpha \in (0, \tfrac{\pi}{3}],\\
    \min \Big\{ \tfrac{1}{j_{1, 1}}, \tfrac{1}{j_{0, 1}}
                            \big(2 (\pi - \alpha) \tan(\tfrac{\alpha}{2})\big)^{-\rfrac{1}{2}} \Big\} &
                                                    \alpha \in (\tfrac{\pi}{3}, \tfrac{\pi}{2}], \\
    \tfrac{1}{j_{0, 1}} \big(2 (\pi - \alpha) \tan(\tfrac{\alpha}{2})\big)^{-\rfrac{1}{2}} &
    \alpha \in (\tfrac{\pi}{2}, \pi).\\
    \end{cases}
    \label{eq:improved-estimates}
\end{equation}
Here, $j_{0, 1} \approx 2.4048$ is the smallest positive
root of the Bessel function $J_0$.
%
%--------------------------------------------------------------------------------------%
A lower bound of ${C^{\mathrm P}_{T}}$ for convex domains in $\Rtwo$ was derived in 
\cite{Cheng1975} and it reads
\begin{equation}
    {C^{\mathrm P}_{T}} \,\geq\, \tfrac{\diam \,(T)}{2 \,j_{0, 1}}.
    \label{eq:cheng}
\end{equation}
%
%This estimate compliments the Payne--Weinberger upper bound,
%and according to \cite{BanuelosBurdzy1999}, is known to be the best lower 
%bound for convex in $\Rtwo$ domains with diameter scaling among all known so far.
Analogously, work \cite{LaugesenSiudeja2009} provides lower bound
\begin{equation}
    {C^{\mathrm P}_{T}} \,\geq\, \tfrac{P}{4 \,\pi}, 
    \label{eq:ls}
\end{equation}
which improves \eqref{eq:cheng} for some cases. Here, $P$ is perimeter of $T$.

%--------------------------------------------------------------------------------------%
In \cite{NazarovRepin2014}, exact values of $\CPGamma$ and $\CtrGamma$ are found for
parallelepipeds, rectangles, and right triangles. Subsequently, we exploit the following 
two results:
%--------------------------------------------------------------------------------------%
\begin{itemize}
\item[1.]
If $T$ is based on vertexes $A = (0, 0)$, $B = (h, 0)$, $C = (0, h)$ and %\linebreak
$\Gamma := \big \{ x_1 \in [0, h], \, x_2 = 0 \big \}$ (i.e., $\Gamma$ coincides with
one of the legs of the isosceles right triangle), then
\begin{equation}
    \CPGamma = \tfrac{h}{\zeta_0} \quad \mathrm{and} \quad
    \CtrGamma = \left(\tfrac{h}{\hat{\zeta}_0 \,
           \tanh({\hat{\zeta}}_0)}\right)^{\rfrac{1}{2}},
    \label{eq:exact-cp-cg-t-leg}
\end{equation}
where $\zeta_0$ and $\hat{\zeta}_0$ are unique roots of the equations
\begin{equation}
z \cot(z) + 1 = 0 \quad \mbox{and} \quad  \tan(z) + \tanh(z) = 0,
\label{eq:roots}
\end{equation}
respectively, in the interval $(0, \pi)$ .
%
%--------------------------------------------------------------------------------------%
\item[2.]
If $T$ is based on vertexes $A = (0, 0)$, $B = (h, 0)$, 
$C = \big(\tfrac{h}{2}, \tfrac{h}{2} \big)$, and $\Gamma$ coincides with the hypotenuse of the isosceles right triangle, then
\begin{equation*}
    \CPGamma = \tfrac{h}{2 \zeta_0} \quad \mathrm{and} \quad
    \CtrGamma = \big(\tfrac{h}{2}\big)^{\rfrac{1}{2}}.
    %\label{eq:exact-cp-cg-t-hyp}
\end{equation*}
%
%For the case $h = 1$, we obtain the reference triangle $\Tref_{\rfrac{\pi}{4}}$.
\end{itemize}
It is worth emphasizing that values of $\CtrGamma$ for right isosceles triangles
follow from the exact solutions of the Steklov problem related to the square. 
This specific case 
was discussed in the work \cite{ArxivGirouardPolterovich}. 

Exact value of constants in the classical Poincar\'{e} inequality are also known for
certain triangles: 
\begin{itemize}
\item[1.] For the equilateral triangle
$\Tref_{\rfrac{\pi}{3}}$ based on vertexes 
$\hat{A} =(0, 0)$, $\hat{B} = (1, 0)$,
$\hat{C} = \big(\tfrac{1}{2}, \tfrac{\sqrt{3}}{2}\big)$, where \linebreak
%
%\begin{equation}
$\hat{\Gamma} := \big\{ x_1 \in [0, 1] ; \;\; x_2 = 0 \big\}$,
%\label{eq:arbitrary-gamma-hat}
%\end{equation}
%
the constant
$$C^{\mathrm{P}}_{\Tref,\, \rfrac{\pi}{3}} = \tfrac{3}{4 \pi}$$
is derived in \cite{Pinsky1980}.
\item[2.] For the right isosceles triangles 
$\Tref_{\rfrac{\pi}{4}}$ based on vertexes 
$\hat{A} =(0, 0)$, $\hat{B} = (1, 0)$, 
$\hat{C} = \big(\tfrac{1}{2}, \tfrac{1}{2}\big)$  
and $\Tref_{\rfrac{\pi}{2}}$ based on 
$\hat{A} = (0, 0)$, $\hat{B} = (1, 0)$, 
$\hat{C} = (0, 1)$,
we have
$$\CPTrefpifour = \tfrac{1}{\sqrt{2}\pi} \quad \mbox{and} \quad \CPTrefpitwo = \tfrac{1}{\pi},$$
respectively. 
Proofs can be found in \cite{HoshikawaUrakawa2010} and 
\cite{NakaoYamamoto2001I}. 
\end{itemize}
Explicit formulas of the same constants for certain 
three-dimensional domains are presented in papers \cite{Berard1980} and 
\cite{HoshikawaUrakawa2010}.

%--------------------------------------------------------------------------------------%
%--------------------------------------------------------------------------------------%
The above mentioned results form a basis for deriving sharp bounds of the constants
$\CPGamma$, $\CtrGamma$, and ${C^{\mathrm P}_{T}}$ for arbitrary non-degenerate 
triangles and
tetrahedrons, which are typical objects in various discretization methods. In Section
\ref{sc:arbitrary-triangle}, we deduce guaranteed and easily computable bounds of
$\CPGamma$, $\CtrGamma$, and ${C^{\mathrm P}_{T}}$ for triangular domains. 
The efficiency of these
bounds is tested in Section \ref{eq:numerial-tests-2d}, where $\CPGamma$ and $\CtrGamma$ 
are compared with
lower bounds computed numerically by solving a generalized eigenvalue problem generated
by Rayleigh quotients discretized 
%by minimization of the respective Rayleigh 
%quotients 
over sufficiently representative sets of trial functions. 
In the same section, we make a
similar comparison of numerical lower bounds related to the constant ${C^{\mathrm P}_{T}}$ with obtained upper bounds and existing estimates known from
\cite{LaugesenSiudeja2009,LaugesenSiudeja2010} and \cite{Cheng1975}. 
Lower bounds of the constants presented in Section \ref{eq:numerial-tests-2d} have been 
computed by two independent codes: the first code is based on the MATLAB Symbolic Math 
Toolbox \cite{Matlab}, and the second one uses The FEniCS Project 
\cite{LoggMardalWells2012}.
Section \ref{eq:numerial-tests-3d} is devoted to tetrahedrons. We combine numerical and
theoretical estimates in order to derive two-sided bounds of the constants. 
Finally,
in Section \ref{sec:example} we present an example that shows one possible application
of the estimates considered in previous sections. Here, the constants are used in order
to deduce a guaranteed and fully computable upper bound of the distance between the
exact solution of an elliptic boundary value problem and an arbitrary function
(approximation) in the respective energy space.

%======================================================================================%
%======================================================================================%

\section{Majorants of $\CPGamma$ and $\CtrGamma$ for triangular domains}
\label{sc:arbitrary-triangle}
Let $T$ be based on vertexes $A = (0, 0)$, $B = (h ,0)$, and 
$C = \big(h \rho \cos\alpha, \, h \rho \, \sin\alpha \big)$ 
%\begin{equation}
%T = {\rm conv}
%\big \{ (0, 0), (h ,0), \big(h \rho \cos\alpha, \, h \rho \, \sin\alpha \big) \big\}
%\label{eq:arbitrary-triangle}
%\end{equation}
%
and
\begin{equation}
\Gamma := \big\{ x_1 \in [0, h] ; \;\; x_2 = 0 \big\},
\label{eq:arbitrary-gamma}
\end{equation}
where $\rho>0$, $h>0$, and $\alpha \in (0,\pi)$ are geometrical parameters
that fully define a triangle $T$ (see Fig. \ref{eq:2d-simplex}). 
Easily computable bounds of $C^{\mathrm{P}}_{\Gamma}$ and $C^{\mathrm{Tr}}_{\Gamma}$
are presented in Lemma \ref{th:lemma-poincare-type-constants}
 below, which uses mappings of reference triangles to 
$T$ and well-known integral transformations
(see, e.g., \cite{Ciarlet1978}). 
%
%--------------------------------------------------------------------------------------%
%
%
%--------------------------------------------------------------------------------------%
% tikz picture
%--------------------------------------------------------------------------------------%
\begin{figure}[h]
\centering
%\adjustbox{valign=b}{
%\begin{minipage}[b]{\linewidth}
\begin{tikzpicture}[scale=0.6]

% coordinate axis Ox
\draw[->] (0.0, 0.0) -- (4.0, 0.0) node[anchor=west] {$x_1$};
\draw[->] (0.0, 0.0) -- (0.0, 4.0) node[anchor=west] {$x_2$};

% another triangle
\draw (0,0) node[anchor=north east] {$A (0, 0)$} --
      (2.8,0) node[anchor=north west] {$B (h, 0)$} --
      (1.8,2.5)  node[anchor=south west]
            {$C \big(h \rho \cos\alpha, h \rho \sin\alpha \big)$}-- (0,0);

\filldraw [black] (0,0) circle (0.8pt)
                  (2.8,0) circle (0.8pt)
                  (1.8,2.5) circle (0.8pt);

\draw (0.4,0.0) .. controls (0.3,0.3) and (0.2, 0.3) .. (0.2,0.3);
\draw (0.6,0.1) node[anchor=south] {$\alpha$};

\draw (1.5,0.5) node[anchor=west] {$T$};
\draw [ultra thick] (0,0) -- (2.8,0);
\draw (2,-1.0) node[anchor=south] {${\Gamma}$};
\end{tikzpicture}

\caption{Simplex in $\Rtwo$.}
\label{eq:2d-simplex}
%\end{minipage}}
%-----------------------------------------------------------------------------------%
%\adjustbox{valign=b}{
%\begin{minipage}[b]{0.55\linewidth} 
%\input{3d-simplex}
%\caption{Simplex in $\Rthree$.}
%\label{eq:3d-simplex}
%\end{minipage}}
\end{figure}

%--------------------------------------------------------------------------------------%
%--------------------------------------------------------------------------------------%
\begin{lemma}
\label{th:lemma-poincare-type-constants}
For any $w \in \tildeH{1}( T, \Gamma)$, the upper bounds of constants in the inequalities
\begin{alignat}{2}
    \|w\|_{ T} \, & \leq\, \CPGamma \, h \,\|\nabla w\|_{ T}
    \quad \mathrm{and} \quad
    \|w\|_\Gamma \, & \leq\, \CtrGamma \, h^{\rfrac{1}{2}} \, \|\nabla w\|_{ T}
    \label{eq:poincare-type-inequalities}
\end{alignat}
are defined as
\begin{equation*}
    \CPGamma \leq \CP = \min \Big \{ \cpleg \, \CPTrefleg, \: \cphyp \, \CPTrefhyp \Big \}
    \quad \mathrm{and} \quad
    \CtrGamma \leq \CG = \min \Big \{ \cgleg \, \CGTrefleg, \: \cghyp \,  \CGTrefhyp \Big \},
    %\label{eq:constants-poincare-type-inequalities}
\end{equation*}
respectively.
Here,
\begin{equation*}
    \cpleg = \muleg^{\rfrac{1}{2}},\quad
    \cphyp = \muhyp^{\rfrac{1}{2}}, \quad		
    \cgleg = \big( \rho \, \sin\alpha \big)^{-\rfrac{1}{2}} \, \cpleg, \quad
    \cghyp = \big( 2 \rho \,\sin\alpha\big)^{-\rfrac{1}{2}}\, \cphyp,
\end{equation*}
%--------------------------------------------------------------------------------------%
where
\begin{alignat}{2}
    \muleg (\rho, \alpha)
    = \,& \tfrac{1}{2}
    \Big( 1 + \rho^2 +
    \big( 1 + \rho^4 + 2  \, \rho^2 \,\cos2\alpha \big)^{\rfrac{1}{2}} \Big),
    \label{eq:mu-leg}\\
    \muhyp (\rho, \alpha)
    = \,& 2 \rho^2 - 2 \rho \, \cos\alpha + 1 +
    \big( (2 \rho^2 + 1)
    (2 \rho^2 + 1 - 4 \rho \, \cos\alpha + 4 \rho^2 \, \cos2\alpha) \big )^{\rfrac{1}{2}},
    \label{eq:mu-hyp}
\end{alignat}
%
%--------------------------------------------------------------------------------------%
and $\CPTrefleg \approx 0.49291$, $\CGTrefleg \approx 0.65602$ and 
$\CPTrefhyp \approx 0.24646$, $\CGTrefhyp \approx 0.70711$, where $\hat{\Gamma}$ is 
defined as
\begin{equation}
\hat{\Gamma} := \big\{ x_1 \in [0, 1] ; \;\; x_2 = 0 \big\}.
\label{eq:arbitrary-gamma-hat}
\end{equation}

\end{lemma}
%
%--------------------------------------------------------------------------------------%
\ProofBegin
Consider the linear mapping $\mathcal{F}_{\rfrac{\pi}{2}} : 
\Tref_{\rfrac{\pi}{2}} \rightarrow T$ with
\begin{equation*}
    x = \mathcal{F}_{\rfrac{\pi}{2}} \, (\hat{x})  = B_{\rfrac{\pi}{2}} \, \hat{x},
    \quad \mbox{where} \quad
    B_{\rfrac{\pi}{2}} =
    \begin{pmatrix}
            \: h & \rho h \cos\alpha \\[0.3em]
                 0 & \rho h \sin\alpha \:
    \end{pmatrix}
    , \quad
    \mathrm{det} B_{\rfrac{\pi}{2}}  = \rho h^2 \, \sin\alpha.
    \label{eq:transformation-1}
\end{equation*}
%
%--------------------------------------------------------------------------------------%
For any $\hat{w} \in \tildeH{1}(\Tref_{\rfrac{\pi}{2}}, \Gref)$, we have the estimate
\begin{equation}
    \|\, \hat{w} \,\|_{ \Tref_{\rfrac{\pi}{2}} }
    \leq \CPTrefleg \, \|\, \nabla \hat{w} \,\|_{\Tref_{\rfrac{\pi}{2}}},
    \label{eq:poicare-inequality-1-hatv}
\end{equation}
where $\CPTrefleg$ is the constant associated with the basic simplex 
$\Tref_{\rfrac{\pi}{2}}$ based on $\hat{A} =(0, 0)$, $\hat{B} = (1, 0)$, and 
$\hat{C} = (0, 1)$,
%
%\begin{equation}
%\Tref_{\rfrac{\pi}{2}} := {\rm conv} \big\{ (0, 0), (1, 0), (0, 1) \big\}.
%\label{eq:t-pi2}
%\end{equation}
%
%--------------------------------------------------------------------------------------%
%
Note that
\begin{alignat}{2}
    \| \, \hat{w} \, \|^2_{\Tref_{\rfrac{\pi}{2}}}
    = \tfrac{1}{\rho h^2 \sin\alpha } \| \,w\, \|^2_{ T},
    \label{eq:v-norm-transformation}
\end{alignat}
and
\begin{equation}
    \| \,\nabla \hat w \, \|^2_{\Tref_{\rfrac{\pi}{2}}}
    \leq \tfrac{1}{\rho h^2 \sin\alpha} \Int_{T}
    A_{\rfrac{\pi}{2}}(h, \rho,\alpha) \nabla w \cdot \nabla w \dx,
    \label{eq:grad-hatv-lower-estimate-1}
\end{equation}
where
\begin{equation*}
    A_{\rfrac{\pi}{2}} (h, \rho,\alpha) = h^2 \,
    \begin{pmatrix}
    1 + \rho^2 \, \cos^2\alpha  \qquad & \rho^2 \sin\alpha \, \cos\alpha \; \\[0.3em]
    \rho^2 \sin\alpha \, \cos\alpha & \rho^2 \sin^2\alpha
    \end{pmatrix}.
    %\label{eq:a-operator}
\end{equation*}
%--------------------------------------------------------------------------------------%
It is not difficult to see that 
$\lambda_{\rm max} (A_{\rfrac{\pi}{2}}) = h^2 \muleg(\rho, \alpha)$, 
%
%\begin{equation*}
    %\lambda_{\rm max} (A_{\rfrac{\pi}{2}}) = h^2 \muleg(\rho, \alpha), \quad
    %\muleg(\rho, \alpha) = \tfrac{1}{2} \Big( 1 + \rho^2 +
    %\big( 1 + \rho^4 + 2 \,\cos2\alpha \, \rho^2 \big)^{\rfrac{1}{2}} \Big).
%\end{equation*}
%
where $\muleg(\rho, \alpha)$ is defined in \eqref{eq:mu-leg}.
From  \eqref{eq:poicare-inequality-1-hatv}, \eqref{eq:v-norm-transformation}, and
\eqref{eq:grad-hatv-lower-estimate-1}, it follows that
%--------------------------------------------------------------------------------------%
\begin{equation}
    \|\, w \,\|_{ T} \, \leq \,
    \cpleg  \, \CPTrefleg \, h \, \, \|\, \nabla w \,\|_{ T}, \quad
    \cpleg (\rho, \alpha) = \muleg^{\rfrac{1}{2}}(\rho, \alpha).
    \label{eq:poicare-inequality-1-v}
\end{equation}
Notice that $\hat{w} \in \tildeH{1} (\Tref, \Gref)$ yields
\begin{equation*}
\mean{w}_{\Gamma} := \Int_{\Gamma} w(x) \ds = 
h \Int_{\Gref} w(x(\hat{x})) \, {\rm d \hat{s}} = 
h \Int_{\Gref} \hat{w} \, \rm{d \hat{s}} = 0.
\end{equation*}
Therefore, above mapping keeps $w \in \tildeH{1} (T, \Gamma)$. 

In view of inequality (\ref{eq:Cgamma}),
for any $\hat{w} \in \tildeH{1}(\Tref_{\rfrac{\pi}{2}}, \Gref)$ we have
\begin{equation*}
    \|\, \hat{w} \,\|_{\Gref} 
		\leq \CGTrefleg \|\, \nabla \hat{w} \,\|_{\Tref_{\rfrac{\pi}{2}}},
    \label{eq:poicare-inequality-2-hatv}
\end{equation*}
where $\CGTrefleg$ is the constant associated with the reference simplex 
$\Tref_{\rfrac{\pi}{2}}$.
%
%
%\begin{equation}
%\Tref_{\rfrac{\pi}{4}} 
%:= {\rm conv} \big\{ (0, 0), (1, 0), \big(\tfrac12, \tfrac12 \big) \big\}.
%\label{eq:t-pi4}
%\end{equation}
%
Since
%
%--------------------------------------------------------------------------------------%
\begin{equation*}
    \| \,\hat{w}\, \|^2_{\Gref} = \tfrac{1}{h} \| \, w \, \|^2_{\Gamma},
		%\label{eq:v-norm-transformation-on-gamma}
\end{equation*}
we obtain
\begin{alignat}{2}
    \|\, w \,\|_{\Gamma} \,
    \leq \, \cgleg \, \CGTrefleg \, h^{\rfrac{1}{2}} \|\, \nabla w \,\|_{ T}, \quad
    \cgleg (\rho, \alpha) = \Big(\tfrac{\muleg(\rho, \alpha)}
    {\rho \sin\alpha}\Big)^{\rfrac{1}{2}}.
    \label{eq:poincare-2-v-cgleg}
\end{alignat}
%
%--------------------------------------------------------------------------------------%
Now, we consider the mapping $\mathcal{F}_{\rfrac{\pi}{4}} : 
\Tref_{\rfrac{\pi}{4}} \rightarrow T$, where $\Tref_{\rfrac{\pi}{4}}$ is based on 
$\hat{A} =(0, 0)$, $\hat{B} = (1, 0)$, and $\hat{C} = (\tfrac{1}{2}, \tfrac{1}{2})$, i.e., 
\begin{equation*}
    x = \mathcal{F}_{\rfrac{\pi}{4}} (\hat{x}) = B_{\rfrac{\pi}{4}} \, \hat{x},
    \quad {\rm where}\quad
    B_{\rfrac{\pi}{4}} =
    \begin{pmatrix}
    \: h & \; 2 \rho h \cos\alpha \, \minus \, h \\[0.3em]
     0   & \; 2 \rho h \sin\alpha \:
    \end{pmatrix},
    \quad 
    \mathrm{det} \, B_{\rfrac{\pi}{4}} = 2 \rho h^2 \, \sin\alpha,
\end{equation*}
which yields another pair of estimates for the functions in
$\tildeH{1}( T, \Gamma)$:
%--------------------------------------------------------------------------------------%
%
\begin{equation}
    \|\, w \,\|_{T} \, \leq \, \cphyp \, \CPTrefhyp \, h \, \|\, \nabla w \,\|_{T}, \quad
    \cphyp (\rho, \alpha) = \muhyp^{\rfrac{1}{2}}(\rho, \alpha),
    \label{eq:poicare-inequality-2-v}
\end{equation}
and
\begin{alignat}{2}
    \|\, w \,\|_{\Gamma} \, \leq \, \cghyp \, \CGTrefhyp \,
    h^{\rfrac{1}{2}} \|\, \nabla w \,\|_{ T}, \quad
    \cghyp (\rho, \alpha)
    = \Big(\tfrac{\muhyp(\rho, \alpha)}{2 \rho \sin\alpha}\Big)^{\rfrac{1}{2}},
    \label{eq:poincare-2-v-cghyp}
\end{alignat}
where $\muhyp(\rho, \alpha)$ is defined in (\ref{eq:mu-hyp}).
Now, (\ref{eq:poincare-type-inequalities}) follows from
(\ref{eq:poicare-inequality-1-v}), (\ref{eq:poincare-2-v-cgleg}),
(\ref{eq:poicare-inequality-2-v}), and (\ref{eq:poincare-2-v-cghyp}).
\ProofEnd

%\begin{remark}
%\rm
%The selection of $\Gamma$ and $\alpha$ in $\T$ depends on 
%the finite element implementation, i.e., we select $\Gamma$ out three edges of $T$ 
%such that it satisfies the condition $\mean{w}_{\Gamma} = 0$. 
%According to the results of numerical experiments, 
%to provide optimal values of $\CPT$, $\Gamma$ must coincide with  
%with longest side of $\T$, and between two adjacent angles we select minimum one
%to be $\alpha$. For the constant $\CGT$, minimum values are attained for $\alpha$ lying 
%in the interval $(\tfrac{\pi}{3}, \tfrac{2\pi}{3})$ 
%(see Section \ref{eq:numerial-tests-2d}).
%\end{remark}

%\begin{remark}
%\label{rm:mean-on-boundary}
%\rm
%\end{remark}

\vskip10pt
%%----------------------------------------------------------------------------------%
Analogously to Lemma \ref{th:lemma-poincare-type-constants}, one  
can obtain an upper bound of the constant in 
\eqref{eq:classical-poincare-constant}. For that we consider three reference 
triangles $\Tref_{\rfrac{\pi}{2}}$, $\Tref_{\rfrac{\pi}{4}}$ (defined earlier), 
and 
$\Tref_{\rfrac{\pi}{3}}$ based on vertexes 
$A = (0, 0)$, $B = (1, 0)$, $C = (\tfrac{1}{2}, \tfrac{\sqrt{3}}{2})$.
%$\T_{\rfrac{\pi}{2}} := {\rm conv} \big\{ (0, 0), (1, 0), (0, 1) \big\}$.
%The exact values of the constants for latter one are obtained in 
%\cite{HoshikawaUrakawa2010}, \cite{Pinsky1980}, 
%and \cite{KikuchiLiu2007}. 
%
\begin{lemma}
\label{th:lemma-poincare-constants}
For any $w \in \tildeH{1}(\T)$, the constant in 
\begin{alignat}{2}
    \|w\|_{\T} \, & \leq\, \CPoincare h \,\|\nabla w\|_{\T},
    \label{eq:poincare-inequality}
\end{alignat}
is estimated as
\begin{equation}
    \CPoincare \leq \CPMR
		= \min \Big \{ \cpifour \, \CPTrefpifour , \: 
		\cpithree \, \CPTrefpithree, \: 
		\cpitwo \, \CPTrefpitwo \Big\}.
		\label{eq:classical-poincare-constant-estimate}
\end{equation}
Here, 
%
%\begin{equation*}
$\cpifour = \mupifour^{\rfrac{1}{2}}$, \quad
$\cpithree = \mupithree^{\rfrac{1}{2}}$,\quad     
$\cpitwo = \mupitwo^{\rfrac{1}{2}}$, where 
%\end{equation*}
%----------------------------------------------------------------------------------%
$\mupitwo$ and $\mupifour$ being defined in \eqref{eq:mu-leg} and 
\eqref{eq:mu-hyp}, respectively, and 
\begin{alignat}{2}
		%\mupifour (\rho, \alpha)
    %= \,& 2 \rho^2 - 2 \rho \, \cos\alpha + 1 +
    %\big( (2 \rho^2 + 1)
    %(2 \rho^2 + 1 - 4 \rho \, \cos\alpha + 4 \rho^2 \, \cos2\alpha) \big )^{\rfrac{1}{2}},
		%\label{eq:mu-pi4} \\		
		\mupithree (\rho, \alpha)
			= \,& \tfrac{2}{3} (1 + \rho^2 - \rho \,\cos\alpha) +
			2 \big( \tfrac{1}{9} (1 + \rho^2 - \rho \,\cos\alpha)^2 - \tfrac{1}{3} \rho^2 \sin^2\alpha \big )^{\rfrac{1}{2}},
		\label{eq:mu-pi3} 
		%\\
    %\mupitwo (\rho, \alpha)
    %= \,& \tfrac{1}{2}
    %\Big( 1 + \rho^2 +
    %\big( 1 + \rho^4 + 2  \, \rho^2 \,\cos2\alpha \big)^{\rfrac{1}{2}} \Big),
    %\label{eq:mu-pi2}
\end{alignat}
and $\CPTrefpifour = \tfrac{1}{\sqrt{2}\pi}$, 
$\CPTrefpithree = \tfrac{3}{4 \pi}$, and
$\CPTrefpitwo = \tfrac{1}{\pi}$.
\end{lemma}
%%----------------------------------------------------------------------------------%
%
\ProofBegin
The mapping $\mathcal{F}_{\rfrac{\pi}{2}}: \Tref_{\rfrac{\pi}{2}} \rightarrow \T$ 
coincides with \eqref{eq:transformation-1} from Lemma 
\ref{th:lemma-poincare-type-constants}. 
It is easy to see that $w \in \tildeH{1} (T)$ provides that $\hat{w} \in \tildeH{1} (\Tref)$.
The estimate 
%%----------------------------------------------------------------------------------%
%For any $\hat{w} \in \tildeH{1}(\Tref_{\rfrac{\pi}{2}})$, we have estimate
%%
%\begin{equation}
    %\|\, \hat{w} \,\|_{ \Tref_{\rfrac{\pi}{2}} }
    %\leq \CPTrefpitwo \, \|\, \nabla \hat{w} \,\|_{\Tref_{\rfrac{\pi}{2}}},
    %\label{eq:classical-poincare-inequality-ref-pi2-hatv}
%\end{equation}
%%
%where $\CPTrefpitwo$ is a constant corresponding to the simplex
%$\Tref_{\rfrac{\pi}{2}} := {\rm conv} \big\{ (0, 0), (1, 0), (0, 1) \big\}$.
%%----------------------------------------------------------------------------------%
%%
%We note that
%%
%\begin{alignat}{2}
    %\| \, \hat{w} \, \|^2_{\Tref_{\rfrac{\pi}{2}}}
    %= \tfrac{1}{\rho h^2 \sin\alpha } \| \,w\, \|^2_{\T},
    %\label{eq:classical-poincare-v-norm-transformation-ref-pi2}
%\end{alignat}
%%
%and
%%
%\begin{equation}
    %\| \,\nabla \hat w \, \|^2_{\Tref_{\rfrac{\pi}{2}}}
    %\leq \tfrac{1}{\rho h^2 \sin\alpha} \Int_{\T}
    %A_{\rfrac{\pi}{2}}(h, \rho,\alpha) \nabla w \cdot \nabla w \dx,
    %\label{eq:classical-poincare-grad-hatv-ref-pi2-lower-estimate-1}
%\end{equation}
%%
%where $A_{\rfrac{\pi}{2}}$ and $\lambda_{\rm max} (A_{\rfrac{\pi}{2}})$ 
%are analyzed in \refeq{eq:a-operator} 
%(Lemma \ref{th:lemma-poincare-type-constants}). 
%Here, $\muleg$ is defined in \eqref{eq:mu-pi2}. Using 
%\eqref{eq:classical-poincare-inequality-ref-pi2-hatv}, 
%\eqref{eq:classical-poincare-v-norm-transformation-ref-pi2}, and 
%\eqref{eq:classical-poincare-grad-hatv-ref-pi2-lower-estimate-1}, we obtain
%%
%----------------------------------------------------------------------------------%
\begin{equation}
    \|\, w \,\|_{\T} \, \leq \,
    \cpitwo  \, \CPTrefpitwo \, h \, \, \|\, \nabla w \,\|_{\T}, \quad
    \cpitwo (\rho, \alpha) = \mupitwo^{\rfrac{1}{2}}(\rho, \alpha)
    \label{eq:classical-poincare-inequality-1-v}
\end{equation}
is obtained by following steps of the previous proof. 
%----------------------------------------------------------------------------------%
From analysis of mappings
\begin{equation*}
    x = \mathcal{F}_{\rfrac{\pi}{3}} (\hat{x}) = B_{\rfrac{\pi}{3}} \, \hat{x},
    \quad \mbox{ where}\quad
		% check in the old note the matrix constant
    B_{\rfrac{\pi}{3}} =
    \begin{pmatrix}
    \: h & \; \tfrac{h}{\sqrt{3}}  (2 \rho \cos\alpha - 1) \, \minus \, h \\[0.3em]
     0   & \; \tfrac{2 h}{\sqrt{3}} \rho \sin\alpha \:
    \end{pmatrix}
    ,\quad
    \mathrm{det} \, B_{\rfrac{\pi}{3}} = \tfrac{2 h^2}{\sqrt{3}} \, \sin\alpha > 0,
\end{equation*}
and 
\begin{equation*}
    x = \mathcal{F}_{\rfrac{\pi}{4}} (\hat{x}) = B_{\rfrac{\pi}{4}} \, \hat{x},
    \quad {\rm where}\quad
    B_{\rfrac{\pi}{4}} =
    \begin{pmatrix}
    \: h & \; 2 \rho h \cos\alpha \, \minus \, h \\[0.3em]
     0   & \; 2 \rho h \sin\alpha \:
    \end{pmatrix}
    ,\quad 
    \mathrm{det} \, B_{\rfrac{\pi}{4}} = 2 \rho h^2 \, \sin\alpha > 0,
\end{equation*}
we obtain alternative estimates 
%----------------------------------------------------------------------------------%
%
\begin{alignat}{2}
    \|\, w \,\|_{\T} \, & \leq \, 
		\cpithree \, \CPTrefpithree \, h \, \|\, \nabla w \,\|_{\T}, \quad
    \cpithree (\rho, \alpha) = \mupithree^{\rfrac{1}{2}}(\rho, \alpha),
    \label{eq:classical-poincare-inequality-2-v} \\
		\|\, w \,\|_{\T} \, & \leq \, 
		\cpifour \, \CPTrefpifour \, h \, \|\, \nabla w \,\|_{\T}, \quad
    \cpifour (\rho, \alpha) = \mupifour^{\rfrac{1}{2}}(\rho, \alpha),
    \label{eq:classical-poincare-inequality-3-v}		
\end{alignat}
where $\mupithree(\rho, \alpha)$ and $\mupifour(\rho, \alpha)$ are defined in 
\eqref{eq:mu-pi3} and \eqref{eq:mu-hyp}, respectively.
Therefore, (\ref{eq:classical-poincare-constant-estimate}) follows from combination of
(\ref{eq:classical-poincare-inequality-1-v}), 
(\ref{eq:classical-poincare-inequality-2-v}), and 
(\ref{eq:classical-poincare-inequality-3-v}). 
%Analogously, if
%$\hat{w} \in \tildeH{1}(\Tref)$,
%it follows that $w \in \tildeH{1}(\T)$.
%
\ProofEnd
%%----------------------------------------------------------------------------------%

%======================================================================================%
%======================================================================================%
\section{Minorants of $\CPGamma$ and $\CtrGamma$ for triangular domains}
\label{eq:numerial-tests-2d}

\subsection{Two-sided bounds of $\CPGamma$ and $\CtrGamma$}
%--------------------------------------------------------------------------------------%
Majorants of $\CPGamma$ and $\CtrGamma$ provided by Lemma 
\ref{th:lemma-poincare-type-constants} should be compared with the corresponding 
minorants, which can be found by means of
%solving generalized eigenvalue problem generated by 
the Rayleigh quotients
\begin{equation}
    \mathcal{R}^{\mathrm{P}}_{\Gamma} [w]
    = \tfrac{\|\nabla w\|_{\T}}{\|w - \mean{w}_{\Gamma} \|_{\T}}
    \quad \mathrm{and} \quad
    \mathcal{R}^{\mathrm{Tr}}_{\Gamma} [w]
    = \tfrac{\|\nabla w\|_{\T}}{\|w - \mean{w}_{\Gamma} \|_{\Gamma}}.
		\label{eq:quotients-2d}
\end{equation}
Lower bounds are obtained if the quotients are minimized on finite dimensional 
subspaces
%Here, $w$ is approximated by the basis of finite dimensional subspaces \linebreak
$V^N \subset \H{1}(T)$ formed by sufficiently
representative collections of suitable test functions. For this purpose, we use either 
power or
Fourier series and introduce the spaces
\begin{equation*}
    V^{N}_1 := \mathrm{span} \big\{\: x^{i} y^{j}\: \big\} \quad {\rm and } \quad
    V^{N}_2 := \mathrm{span} \big\{\: \cos (\pi i x) \cos (\pi j y)\: \big\}, \;\;
    %\label{eq:basis}
\end{equation*}
where $i, j = 0, \ldots, N, \;\; (i, j) \neq (0, 0)$
and 
$$\mathrm{dim} \, V^{N}_1 = \mathrm{dim} \, V^{N}_2 = M(N) :=
(N + 1)^2 - 1.$$
The corresponding constants are denoted by $\approxCPT$ and $\approxCGT$, where $M$ indicates 
on number of basis functions in auxiliary subspace used.
%, where 
%$M(N)$ indicates the amount of used basis functions in the finite dimensional subspace. 
Since $ V^{N}_1$ and $V^{N}_2$ are limit dense
in $\H{1}(\T)$, the respective minorants tend to the exact constants as $M(N)$ tends to infinity.

We note that  
\begin{equation}
		\inf\limits_{w \in \H{1}(T)} \mathcal{R}^{\mathrm{P}}_{\Gamma} [w] = 
    \inf\limits_{w \in \H{1}(T)} 
		\tfrac{\|\nabla w\|_{\T}}{\|w - \mean{w}_{\Gamma} \|_{\T}}
     =
    \inf\limits_{w \in \tildeH{1}(T, \Gamma)} 
		\tfrac{\|\nabla w\|_{\T}}{\|w\|_{\T}} = \tfrac{1}{\CPGamma}.
\end{equation}
%
%Indeed, where $\mathcal{R}^{\mathrm{p}}_{\Gamma} [w]$ follows from the definition of the constant 
%$\CPGamma$ 
%For any $w \in \H{1} (T)$, we have
%%
%\begin{equation}
	%\|w - \mean{w}_{\Gamma} \|_{T} \, \leq \CPGamma \|\nabla w\|_{T}.
	%%{\rm and} \quad 
	%%\|w - \mean{w}_{\Gamma} \|_{\Gamma} \, \leq \CtrGamma \|\nabla w\|_{T}.
	%%\label{eq:full-poincare-type}
%\end{equation}
%
%Embeddings \eqref{eq:full-poincare-type} are justified by the equivalence of quantity 
%\linebreak
%$\| w \|_l := \| \nabla w \|_{\T} + \big| \Int_{\Gamma} w \ds \big|$ to the norm of 
%$\H{1} (T)$, so that the existence of $\CPGamma$ and $\CtrGamma$ follows automatically.
%
Therefore, minimization of the first quotient in \eqref{eq:quotients-2d} on 
$V^{N}_1$ or $V^{N}_2$ yields a lower bound of $\CPGamma$.
% $w \in \tildeH{1}(T, \Gamma)$
%we can compare the first quotient in \eqref{eq:quotients-2d} with estimate obtained in 
%left estimates in \eqref{eq:poincare-type-inequalities} 
%for function in $\tildeH{1}(T, \Gamma)$. 
For the quotient 
$\mathcal{R}^{\mathrm{Tr}}_{\Gamma} [w]$, we apply similar arguments. 
 
%--------------------------------------------------------------------------------------%
Numerical results presented below are obtained with the help of two different codes 
based on the MATLAB Symbolic Math Toolbox \cite{Matlab} and The FEniCS Project 
\cite{LoggMardalWells2012}.
%--------------------------------------------------------------------------------------%
Table \ref{tab:const-convergence-from-basis} demonstrates the ratios between the 
exact 
constants and respective approximate values (for the selected $\rho$ and $\alpha$).
They are quite close to $1$ even for relatively small
$N$. Henceforth, we select $N = 6$ or $7$ in the tests discussed below. 

%We consider two arbitrary simplexes. The first one coincides with the reference triangle
%$\Tref_{\rm I}$ ($\alpha = \tfrac{\pi}{2}$, $h = 1$, and $\rho = 1$), and the second
%is identical to $\Tref_{\rm II}$ ($\alpha = \tfrac{\pi}{4}$, $h = 1$, and
%$\rho = \tfrac{\sqrt{2}}{{2}}$). For both cases, $\cpleg$, $\cgleg$ and $\cphyp$,
%$\cghyp$ in (\ref{eq:constants-poincare-type-inequalities}) must be equal to 1. This
%fact is confirmed by the Table \ref{tab:const-convergence-from-basis}, in which lower
%bounds of the weighting parameters converge to $1$, if $M$ increases.

\begin{table}[!ht]
\centering
\footnotesize
\begin{tabular}{cc|cc|cc}
\multicolumn{2}{c|}{$ $}
& \multicolumn{2}{c|}{ $\alpha = \tfrac{\pi}{2}$, $\rho = 1$ }
& \multicolumn{2}{c}{ $\alpha = \tfrac{\pi}{4}$, $\rho = \tfrac{\sqrt{2}}{2}$} \\
\midrule
$N$ & $M(N)$ 
& %$\rfrac{\approxCPT}{\CPTrefleg}$ 
${\approxCPT}/{\CPTrefleg}$ 
& ${\approxCGT}/{\CGTrefleg}$ 
& ${\approxCPT}/{\CPTrefhyp}$ 
& ${\approxCGT}/{\CGTrefhyp}$ \\
\midrule
1 & 3  & 0.8801 & 0.9561 & 0.8647 & 1.0000 \\
2 & 8  & 0.9945 & 0.9898 & 0.9925 & 1.0000 \\
3 & 15 & 0.9999 & 0.9998 & 0.9962 & 1.0000 \\
4 & 24 & 1.0000 & 0.9999 & 1.0000 & 1.0000 \\
5 & 35 & 1.0000 & 1.0000 & 1.0000 & 1.0000 \\
6 & 48 & 1.0000 & 1.0000 & 1.0000 & 1.0000 \\
\end{tabular}
\\[5pt]
\caption{Ratios between approximate and reference constants 
with respect to increasing $N$.}
\label{tab:const-convergence-from-basis}
\end{table}

%--------------------------------------------------------------------------------------%

%-------------------------------------------------------------------------------------------%
% \CPT with upper bounds
%--------------------------------------------------------------------------------------------%

\begin{figure}[!ht]
	\centering
	\subfloat[$\rho = \tfrac{\sqrt{2}}{2}$]{
  \includegraphics[scale=0.95]{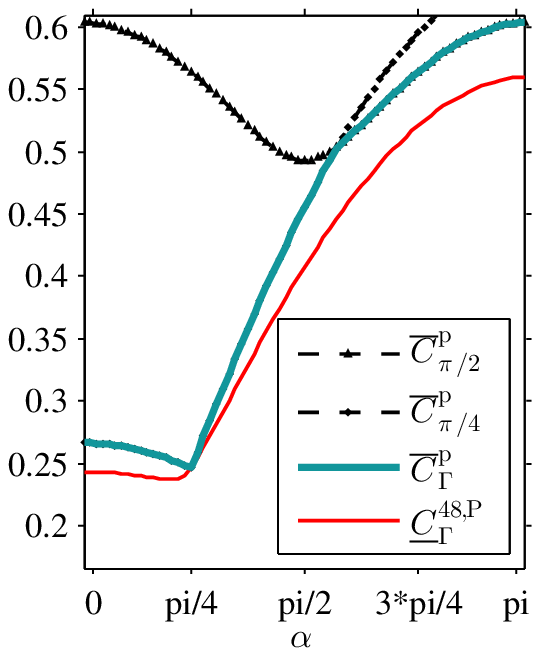}
	\label{fig:2d-cpt-rho-sqrt2-2}}
	\quad
	\subfloat[$\rho = \tfrac{\sqrt{2}}{2}$]{
	\begin{tikzpicture}[thick, scale=0.9]
	\draw[dashed] (2.1213,0) arc (0:180:2.1213cm);

	% coordinate axis Ox
	\draw[->] (0.0, 0.0) -- (3.5, 0.0) node[anchor=west] {$x_1$};
	\draw[->] (0.0, 0.0) -- (0.0, 3.5) node[anchor=west] {$x_2$};

	% another triangle
	\draw (0.0,0.0) node[anchor=north east] {$(0, 0)$} --
				(2.1213,0.0) node[anchor=north west] {$(1, 0)$} --
				(0.0,2.1213)  node[anchor=south west] {$\big(0, \tfrac{1}{\sqrt{2}}\big)$}-- (0,0);

	\filldraw [black] (0.0,0.0) circle (0.8pt)
										(2.1213,0.0) circle (0.8pt)
										(0.0,2.1213) circle (0.8pt);

	\draw (0,0.2) -- (0.2,0.2) -- (0.2, 0); 
	\draw (0.4, -0.6) node[anchor=south] {$\tfrac{\pi}{2}$};
	
	\draw [ultra thick] (0,0) -- (3,0); 
	\draw (1.5,-0.8) node[anchor=south] {${\Gamma}$};
	
	%---------------------------------------------------------------%
	
	\draw [ultra thick](0.0,0.0) -- (3.0, 0.0) -- (1.5, 1.5) -- (0,0);
	
	\filldraw [black] (0.0,0.0) circle (0.8pt)
										(3.0,0.0) circle (0.8pt)
										(1.5, 1.5) circle (0.8pt);
	\draw (1.0,0.0) .. controls (1.0,0.2) and (0. 95, 0.5) .. node[right] {$\tfrac{\pi}{4}$}  (0.8,0.72);

	%---------------------------------------------------------------%

	\draw (0.0,0.0) -- (3.0, 0.0) -- (1.0607, 1.8371) -- (0,0);
	
	\filldraw [black] (0.0,0.0) circle (0.8pt)
										(3.0,0.0) circle (0.8pt)
										(1.0607, 1.8371) circle (0.8pt);
	\draw (0.8, 0.0) .. controls (0.8, 0.4) and (0.7, 0.58) .. node[anchor=south] {$ $}  (0.49,0.82);
  \draw (0.3, 0.6) node[anchor=south] {$\tfrac{\pi}{3}$};
	
	%---------------------------------------------------------------%

	\draw (0.0,0.0) -- (3.0, 0.0) -- (-1.0607, 1.8371) -- (0,0);
	
	\filldraw [black] (0.0,0.0) circle (0.8pt)
										(3.0,0.0) circle (0.8pt)
										(-1.0607, 1.8371) circle (0.8pt);
	\draw (0.6, 0.0) .. controls (0.6,0.7) and (-0.1, 0.9) .. node[anchor=south] {$ $}  (-0.4,0.7);
  \draw (-0.2, 0.7) node[anchor=south] {$\tfrac{2\pi}{3}$};
	
	\end{tikzpicture}
	\label{eq:t-rho-sqrt2-2}
 	}\\
	\subfloat[$\rho = 1$]{
  \includegraphics[scale=0.95]{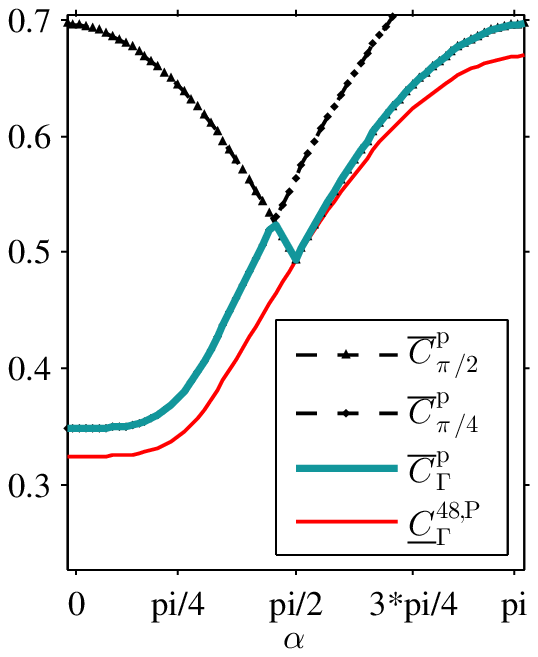}
	\label{fig:2d-cpt-rho-1}}
	\; \;
	\subfloat[$\rho = 1$]{
  \begin{tikzpicture}[thick, scale=0.9]

	\draw[dashed] (3,0) arc (0:180:3cm);

	% coordinate axis Ox
	\draw[->] (0.0, 0.0) -- (3.2, 0.0) node[anchor=west] {$x_1$};
	\draw[->] (0.0, 0.0) -- (0.0, 3.5) node[anchor=west] {$x_2$};

	% another triangle
	\draw[ultra thick] (0.0,0.0) node[anchor=north east] {$(0, 0)$} --
				(3.0,0.0) node[anchor=north west] {$(1, 0)$} --
				(0.0,3.0)  node[anchor=east] {$(0, 1)$}-- (0,0);

	\filldraw [black] (0.0,0.0) circle (0.8pt)
										(3.0,0.0) circle (0.8pt)
										(0.0,3.0) circle (0.8pt);

	\draw (0,0.2) -- (0.2,0.2) -- (0.2, 0); 
	\draw (0.4, -0.6) node[anchor=south] {$\tfrac{\pi}{2}$};
	
	\draw [ultra thick] (0,0) -- (3,0); 
	\draw (1.5,-0.8) node[anchor=south] {${\Gamma}$};
	
	%---------------------------------------------------------------%
	
	\draw (0.0,0.0) -- (3.0, 0.0) -- (2.5981,1.5) -- (0,0);
	
	\filldraw [black] (0.0,0.0) circle (0.8pt)
										(3.0,0.0) circle (0.8pt)
										(2.5981,1.5) circle (0.8pt);
	\draw (1.0,0.0) .. controls (1.0,0.2) and (0. 95, 0.4) .. node[right] {$\tfrac{\pi}{6}$}  (0.89,0.5);

	%---------------------------------------------------------------%

	\draw (0.0,0.0) -- (3.0, 0.0) -- (1.5,2.5981) -- (0,0);
	
	\filldraw [black] (0.0,0.0) circle (0.8pt)
										(3.0,0.0) circle (0.8pt)
										(1.5,2.5981) circle (0.8pt);
	\draw (0.8, 0.0) .. controls (0.8,0.4) and (0.7, 0.58) .. node[anchor=south] {$ $}  (0.49,0.82);
  \draw (0.78, 0.6) node[anchor=south] {$\tfrac{\pi}{3}$};
	
	%---------------------------------------------------------------%

	\draw (0.0,0.0) -- (3.0, 0.0) -- (-1.5,2.5981) -- (0,0);
	
	\filldraw [black] (0.0,0.0) circle (0.8pt)
										(3.0,0.0) circle (0.8pt)
										(-1.5,2.5981) circle (0.8pt);
	\draw (0.6, 0.0) .. controls (0.6,0.7) and (-0.1, 0.9) .. node[anchor=south] {$ $}  (-0.4,0.7);
  \draw (-0.2, 0.7) node[anchor=south] {$\tfrac{2\pi}{3}$};
	
	\end{tikzpicture}
	\label{eq:t-rho-1}}
	\\[5pt] 
	\caption{Two-sided bounds of $\CPT$
	for $\T$ with different $\rho$.}
	\label{fig:2d-cpt}
\end{figure}
%
%--------------------------------------------------------------------------------------------%
% \CGT with upper bounds
%--------------------------------------------------------------------------------------------%
\begin{figure}[!ht]
	\centering
	\subfloat[$\rho = \tfrac{\sqrt{2}}{2}$]{
  \includegraphics[scale=0.95]{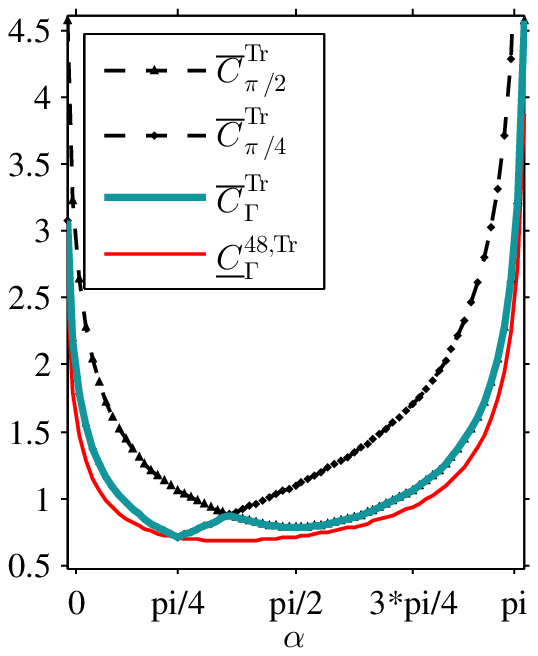}
	\label{fig:2d-cgt-rho-sqrt2-2}}
	\quad
	\subfloat[$\rho = 1$]{
  \includegraphics[scale=0.95]{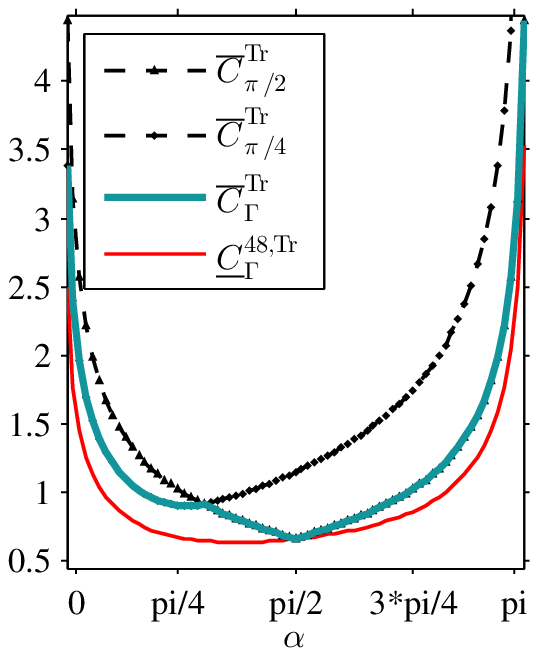}
	\label{fig:2d-cgt-rho-1}}
	\\[5pt] 
	\caption{Two-sided bounds of $\CGT$
	%$\approxCGT$ and upper bound 
	%$\CG = \min \,\{ {\CG}_{\rm I}, {\CG}_{\rm II} \}$, 
	for $\T$ with different $\rho$.}
	\label{fig:2d-cgt}
\end{figure}

In Figs. \ref{fig:2d-cpt-rho-sqrt2-2} and \ref{fig:2d-cpt-rho-1}, we depict $\approxCPT$ 
for $M(N) = 48$ (thin red line) for $T$ with $\rho = \tfrac{\sqrt{2}}{2}$, 
$1$, and $\alpha \in (0, \pi)$. Guaranteed upper bounds 
$\CPTleg = \cpleg \, \CPTrefleg$ and
$\CPThyp = \cphyp \, \CPTrefhyp$ are depicted by 
dashed black lines. Bold blue line illustrates $\CP = \min \Big \{\CPTleg, \CPThyp \Big \}$.
Analogously in Figs. \ref{fig:2d-cgt-rho-sqrt2-2} and \ref{fig:2d-cgt-rho-1}, a red marker denotes
the lower bound $\approxCGT$ (for $M(N) = 48$) of the constant $\CtrGamma$. It
is presented together with the upper bound $\CG$ (blue marker), which is defined as minimum of
$\CGTleg = \cgleg \, \CGTrefleg$ and $\CGThyp = \cghyp \, \CGTrefhyp$.
%Parameter $M$ is fixed to $48$ since the difference of order $1e\minus8$ 
%between $\approxCGT$ 
%(for bigger $M$) becomes unnoticeable.
Table \ref{tab:big-table} represents this information in the digital form.

%Figs \ref{eq:t-rho-sqrt2-2} and \ref{eq:t-rho-1} illustrate the change of the simplex
%$T$ with respect to the angle $\alpha$. For convenience of the reader, the data
%depicted in Figs \ref{fig:2d-cpt-rho-sqrt2-2}--\ref{fig:2d-cgt-rho-1} are represented
%in Table \ref{tab:big-table} .

%--------------------------------------------------------------------------------------%
% big table with c+ and c-
%--------------------------------------------------------------------------------------%
\begin{table}[!ht]
\centering
\footnotesize
\begin{tabular}{c|cc|cc|cc|cc}
\multicolumn{1}{c|}{$ $}
& \multicolumn{4}{c|}{ $\rho = \tfrac{\sqrt{2}}{2}$ }
& \multicolumn{4}{c}{  $\rho = 1$ } \\
\midrule
$\alpha$ &
$\underline{C}^{48, \mathrm{P}}_{\Gamma}$ & $\CP$ &
$\underline{C}^{48, \mathrm{Tr}}_{\Gamma}$ & $\CG$ &
$\underline{C}^{48, \mathrm{P}}_{\Gamma}$ & $\CP$ &
$\underline{C}^{48, \mathrm{Tr}}_{\Gamma}$ & $\CG$ \\
\midrule
$\pi/18$ & 0.2429 & 0.2657 & 1.2786 & 1.5386 & 0.3245 & 0.3486 & 1.2572 & 1.6971\\
$\pi/9$ & 0.2414 & 0.2627 & 0.9289 & 1.0838 & 0.3248 & 0.3493 & 0.9058 & 1.2116\\
$\pi/6$ & 0.2389 & 0.2577 & 0.7919 & 0.8792 & 0.3268 & 0.3527 & 0.7632 & 1.0118\\
$2\pi/9$ & 0.2379 & 0.2507 & 0.7259 & 0.7543 & 0.3339 & 0.3636 & 0.6906 & 0.9201\\
$5\pi/18$ & 0.2632 & 0.2722 & 0.6945 & 0.7503 & 0.3514 & 0.3884 & 0.6529 & 0.9003\\
$\pi/3$ & 0.3008 & 0.3220 & 0.6829 & 0.8348 & 0.3809 & 0.4269 & 0.6362 & 0.8634\\
$7\pi/18$ & 0.3382 & 0.3694 & 0.6840 & 0.8432 & 0.4173 & 0.4721 & 0.6332 & 0.7840\\
$4\pi/9$ & 0.3740 & 0.4140 & 0.6947 & 0.7973 & 0.4556 & 0.5187 & 0.6404 & 0.7162\\
$\pi/2$ & 0.4075 & 0.4554 & 0.7136 & 0.7801 & 0.4929 & 0.4929 & 0.6560 & 0.6560\\
$5\pi/9$ & 0.4382 & 0.4933 & 0.7409 & 0.7973 & 0.5280 & 0.5340 & 0.6797 & 0.7162\\
$11\pi/18$ & 0.4660 & 0.5165 & 0.7779 & 0.8432 & 0.5600 & 0.5710 & 0.7125 & 0.7840\\
$2\pi/3$ & 0.4905 & 0.5361 & 0.8274 & 0.9118 & 0.5884 & 0.6037 & 0.7569 & 0.8634\\
$13\pi/18$ & 0.5115 & 0.5552 & 0.8948 & 1.0040 & 0.6129 & 0.6318 & 0.8175 & 0.9607\\
$7\pi/9$ & 0.5289 & 0.5720 & 0.9898 & 1.1292 & 0.6332 & 0.6550 & 0.9033 & 1.0874\\
$5\pi/6$ & 0.5426 & 0.5856 & 1.1334 & 1.3107 & 0.6492 & 0.6733 & 1.0332 & 1.2673\\
$8\pi/9$ & 0.5524 & 0.5956 & 1.3796 & 1.6118 & 0.6607 & 0.6865 & 1.2565 & 1.5623\\
$17\pi/18$ & 0.5583 & 0.6017 & 1.9436 & 2.2851 & 0.6676 & 0.6944 & 1.7692 & 2.2179\\
\end{tabular}
\caption{Two-sided bounds of $\CPT$ and $\CGT$ for $\T$ for $\alpha \in (0, \pi)$ and 
different $\rho$.}
\label{tab:big-table}
\end{table}

%--------------------------------------------------------------------------------------%
Fig. \ref{fig:2d-cpt-rho-sqrt2-2} corresponds to the case $\rho =  \tfrac{\sqrt{2}}{2}$.
Notice that for $\alpha = \tfrac{\pi}{4}$ the constant $\CPGamma$ is known and the computed
lower bound $\approxCPT$ (red marker) practically coincides with it (see, e.g., Fig. \ref{eq:t-rho-sqrt2-2}). 
Since in this case, the mapping $\mathcal{F}_{\rfrac{\pi}{4}}$ is identical, the upper bound also
coincides with the exact value. An analogous coincidence can be observed for $\CtrGamma$
and $\approxCGT$ in Fig. \ref{fig:2d-cgt-rho-sqrt2-2}.
In Fig. \ref{fig:2d-cpt-rho-1}, the red curve, corresponding to $\approxCPT$, coincides with
the blue line of $\CPT$ at the point $\alpha = \tfrac{\pi}{2}$ (due to the fact that
for this angle $\mathcal{F}$ is the identical mapping and $T$ coincides with
$\Tref_{\rfrac{\pi}{2}}$ (see Fig. \ref{eq:t-rho-1})). Fig. \ref{fig:2d-cgt-rho-1} exposes 
similar results for ${\approxCGT}$ and $\CtrGamma$ ($\CGTleg$).

%-------------------------------------------------------------------------------------------%
% \CPT with upper bounds
%--------------------------------------------------------------------------------------------%

\begin{figure}[!ht]
	\centering
	\subfloat[$\rho = \tfrac{\sqrt{3}}{2}$]{
  \includegraphics[scale=0.9]{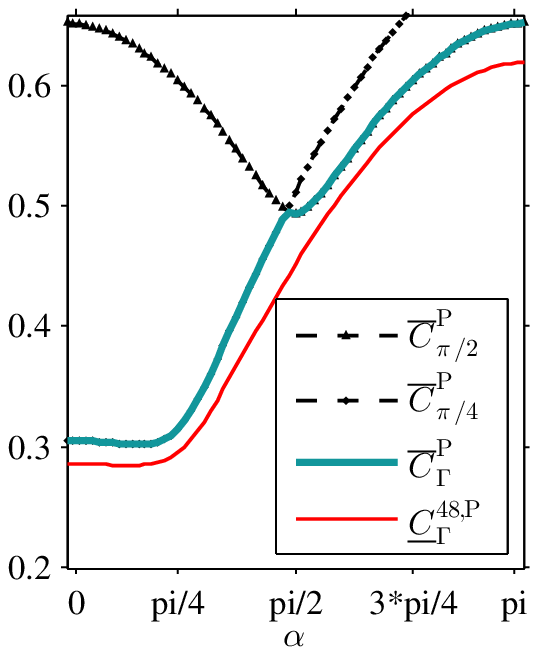}
	\label{fig:2d-cpt-rho-sqrt3-2}}
	\quad
	\subfloat[$\rho = \tfrac{3}{2}$]{
	\includegraphics[scale=0.9]{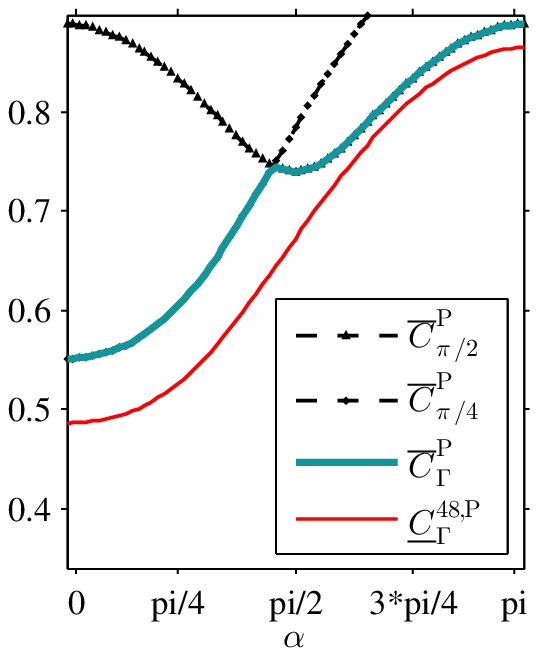}
	\label{fig:2d-cpt-rho-3-2}}
	\\[5pt] 
	\caption{Two-sided bounds of $\CPT$
	for $\T$ with different $\rho$.}
	\label{fig:2d-cpt-no-ref}
\end{figure}
%
%--------------------------------------------------------------------------------------------%
% \CGT with upper bounds
%--------------------------------------------------------------------------------------------%
\begin{figure}[!ht]
	\centering
	\subfloat[$\rho = \tfrac{\sqrt{3}}{2}$]{
  \includegraphics[scale=0.9]{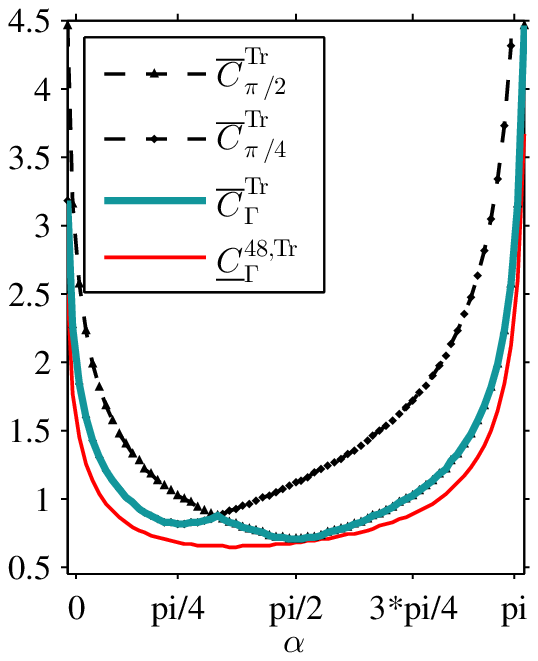}
	\label{fig:2d-cgt-rho-sqrt3-2}}
	\quad
	\subfloat[$\rho = \tfrac{3}{2}$]{
	\includegraphics[scale=0.9]{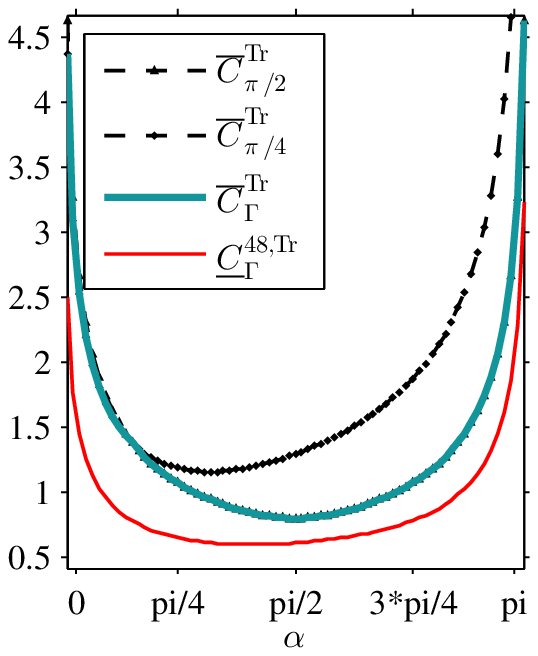}
	\label{fig:2d-cgt-rho-3-2}}
	\\[5pt] 
	\caption{Two-sided bounds of $\CGT$
	for $\T$ with different $\rho$.}
	\label{fig:2d-cgt-no-ref}
\end{figure}
%
%--------------------------------------------------------------------------------------%

Figs. \ref{fig:2d-cpt-no-ref} and \ref{fig:2d-cgt-no-ref} demonstrate the same bounds 
for $\rho = \tfrac{\sqrt{3}}{2}$ and $\tfrac{3}{2}$.
We see that estimates of $\CPT$ and $\CGT$ are very efficient. Namely, 
$\Ieff^{\rm P} := \tfrac{\CP}{\underline{C}^{48, \mathrm{P}}_{\Gamma}} \in [1.0463, 1.1300]$ 
for $\rho = \tfrac{\sqrt{3}}{2}$ and 
$\Ieff^{\rm P} \in [1.0249, 1.1634]$ for $\rho = \tfrac{3}{2}$. 
Analogously, 
$\Ieff^{\rm Tr} := \tfrac{\CG}{\underline{C}^{48, \mathrm{Tr}}_{\Gamma}} \in [1.0363, 1.3388]$ 
for $\rho = \tfrac{\sqrt{3}}{2}$ and $\Ieff^{\rm Tr} \in [1.2917, 1.7643]$ for $\rho = \tfrac{3}{2}$. 

%--------------------------------------------------------------------------------------%
\subsection{Two-sided bounds of ${C^{\mathrm P}_{T}}$}
%--------------------------------------------------------------------------------------%

The spaces $V^{N}_1$ and $ V^{N}_2$ can also be used for analysis of the quotient
%\linebreak
$\mathcal{R}_{T} [w] = \tfrac{\|\nabla w\|_{ T}}{\|w - \mean{w}_{ T} \|_{ T}}$,
which yields guaranteed lower bounds of the constant in 
(\ref{eq:classical-poincare-constant}). The respective values are denoted by $\approxC$. 
These bounds are compared with 
$\CPupper := \tfrac{\diam (\T)}{j_{1, 1}}$ and $\CPlower :=
\max \Big \{ \tfrac{\diam (\T)}{2 \,j_{0,1}},  \tfrac{P}{4 \,\pi} \Big \}$ 
(see \eqref{eq:cheng}--\eqref{eq:ls}, respectively) as well as the one derived in 
Lemma \ref{th:lemma-poincare-constants}.
%--------------------------------------------------------------------------------------%
%********************************************************************************************
%rho = 0.7, 0.86, 1, 1.5
%********************************************************************************************
\begin{figure}[!ht]
	\centering
	\subfloat[$\rho = \tfrac{\sqrt{2}}{2}$]{
	\includegraphics[scale=0.9]{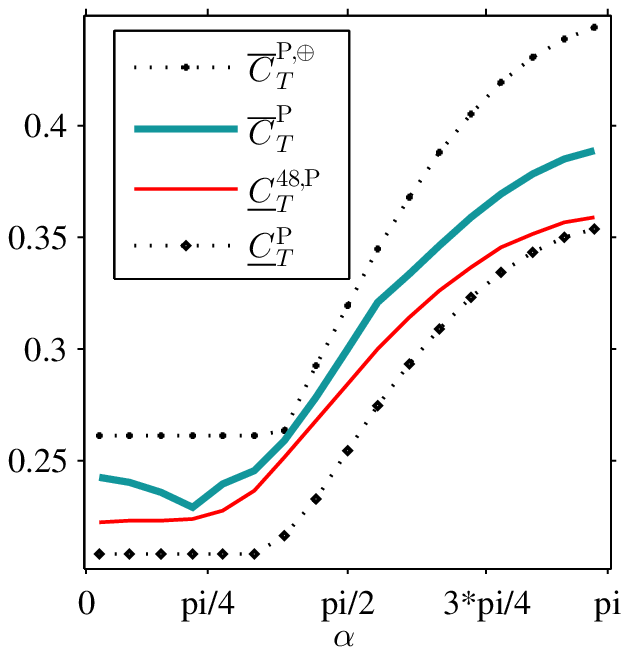}	\label{fig:cp-rho-sqrt2-2}
	}
	\subfloat[$\rho = \tfrac{\sqrt{3}}{2}$]{
	\includegraphics[scale=0.9]{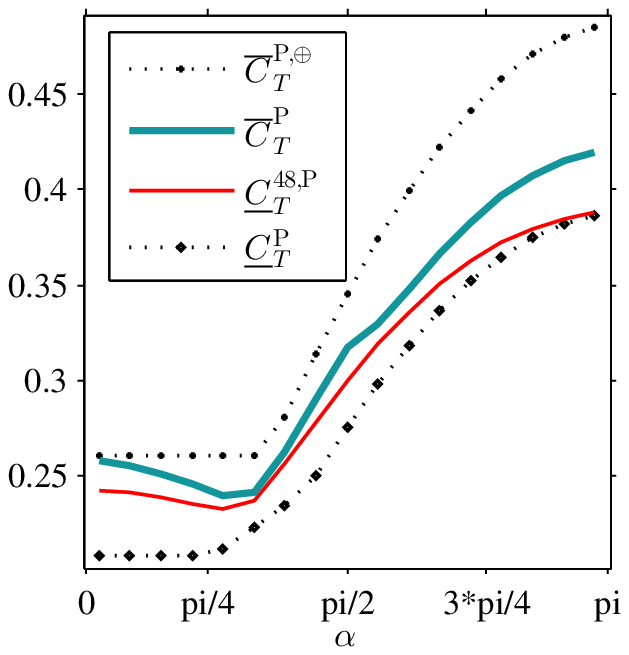}	\label{fig:cp-rho-sqrt3-2}
	}\\[-1pt]
	\subfloat[$\rho = 1$]{
	\includegraphics[scale=0.9]{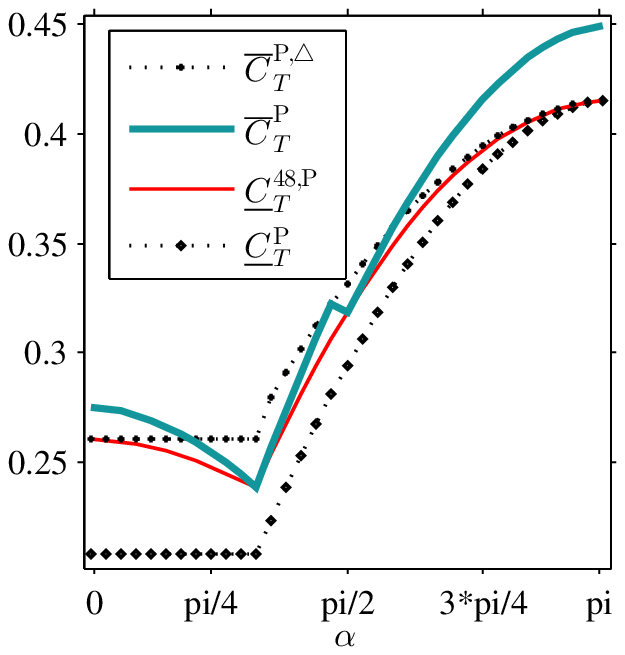}
	\label{fig:cp-rho-1}}
	\subfloat[$\rho = \tfrac{3}{2}$]{
	\includegraphics[scale=0.9]{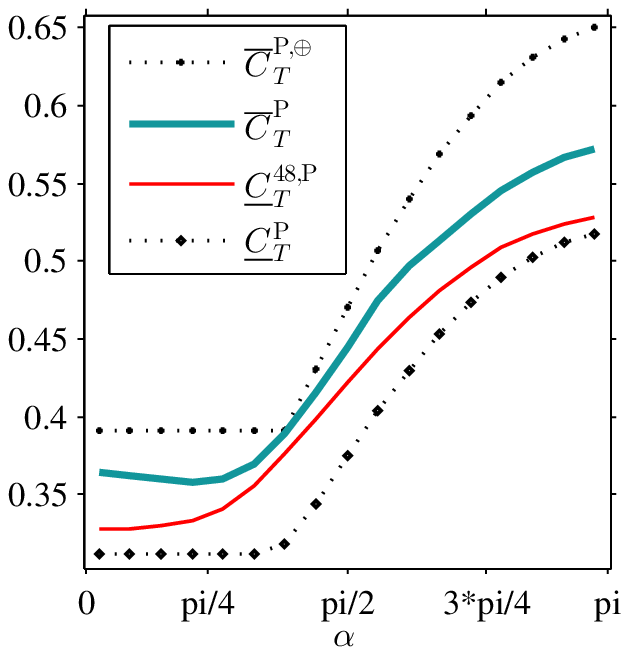}
	\label{fig:cp-rho-3-2}}
	\\[5pt] 
	\caption{$\underline{C}^{48}_{\T}$, $\CPLS$, $\CPMR$, and $\CPlower$ for $\T$ with 
	$\alpha \in (0, \pi)$ and different $\rho$.}	
	\label{fig:2d-cp-rho}
\end{figure}

In Figs. \ref{fig:cp-rho-sqrt2-2}, \ref{fig:cp-rho-sqrt3-2}, and \ref{fig:cp-rho-3-2}, 
we present $\approxC$ (in this case $M(N) = 48$) together with $\CPMR$ (blue think line), 
$\CPupper$ and $\CPlower$ for $\alpha \in (0, \pi)$, and  
$\rho = \tfrac{\sqrt{2}}{2}$, $\tfrac{\sqrt{3}}{2}$, and $\tfrac{3}{2}$. 
We see that $\underline{C}^{48, \rm P}_{\T}$ (red thin line) indeed lies within the admissible 
two-sided bounds. From these figures,
it is obvious that new upper bounds $\CPMR$ are sharper than $\CPupper$ for $\T$ with $\rho \neq 1$.
True values of the constant lie between the bold blue and thin red lines, but closer to the 
red one, which practically shows the constant (this follows from the fact 
that increasing $M(N)$ does not provide a noticeable change for the line, e.g., 
for $M(N) = 63$ 
maximal difference with respect to figure does not exceed $1e\minus8$). Also, we note that,
the lower bound $\CPlower$ (black dashed line) is quite efficient, and, moreover, asymptotically exact for
$\alpha \rightarrow \pi$. 

%\clearpage

Due to \cite{LaugesenSiudeja2010} and \cite{Bandle1980}, 
we know the improved upper bound $\CPLS$ 
(cf. \eqref{eq:improved-estimates}) for isosceles triangles. In Fig. \ref{fig:cp-rho-1}, 
we compare $\approxC$ ($M(N) = 48$) with both upper bounds  $\CPMR$ (from the Lemma 
\ref{th:lemma-poincare-constants}) and $\CPLS$ (black doted line). 
It is easy to see that $\CPLS$ (black dashed line) is 
rather accurate and for
$\alpha \rightarrow 0$ and $\alpha \rightarrow \pi$ provide 
almost exact estimates. $\CPMR$ (blues thick line) improves $\CPLS$ only for some $\alpha$. 
The lower bound $\underline{C}^{48}_{\T}$ (red thin line) indeed converges to
$\CPLS$ as $\T$ degenerates when $\alpha$ tends to  $0$.

%--------------------------------------------------------------------------------------%
\subsection{Shape of the minimizer}
%--------------------------------------------------------------------------------------%
%

Exact constants in (\ref{eq:Comega}) and (\ref{eq:Cgamma}) are generated by the minimal 
positive eigenvalues of \eqref{eq:eigenvalue-problem-cp}
and \eqref{eq:eigenvalue-problem-ctr}. This section presents results related to 
the respective eigenfunctions. In order to depict all of them in a unified form,
we use the barycentric coordinates $\lambda_i \in (0, 1)$, $i = 1, 2, 3$,
$\Sum_{i = 1}^{3}\lambda_i = 1$.
\begin{figure}[!ht]
	\centering
	\subfloat[$u^{48, {\rm p}}_{\Gamma}$, $\alpha = \tfrac{\pi}{6}$]{\includegraphics[scale=0.9]{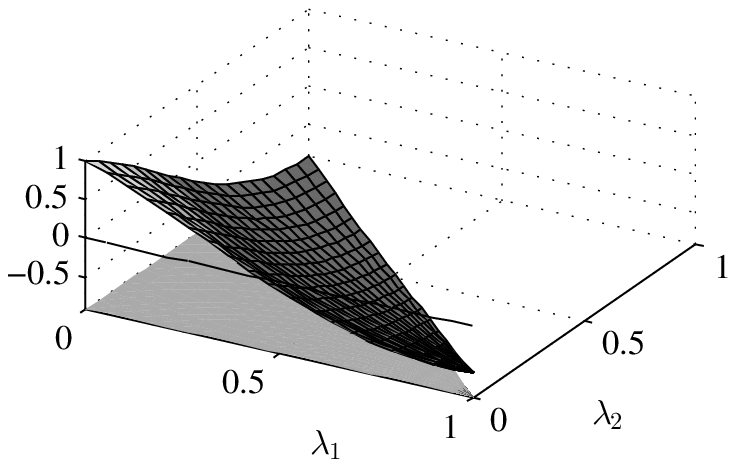}
	\label{fig:2d-up-eigenfunctions-pi-6}
	}
	\subfloat[$u^{48, {\rm p}}_{\Gamma}$, $\alpha = \tfrac{\pi}{3}$]{\includegraphics[scale=0.9]{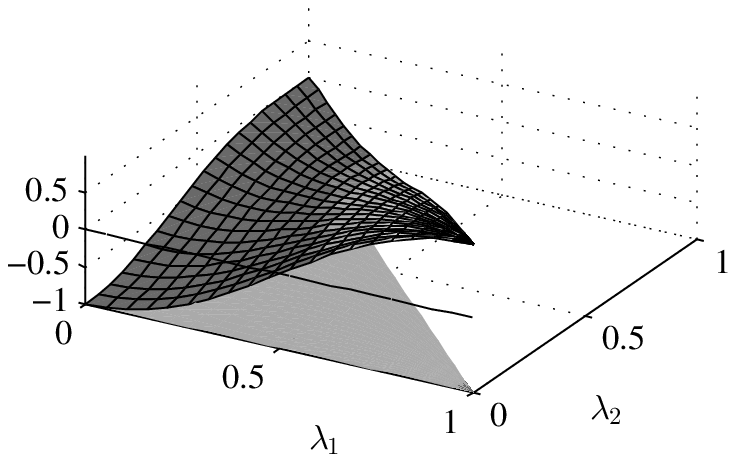}
	\label{fig:2d-up-eigenfunctions-pi-3}
	}\\
	\subfloat[$u^{48, {\rm p}}_{\Gamma}$, $\alpha = \tfrac{\pi}{2}$]{\includegraphics[scale=0.9]{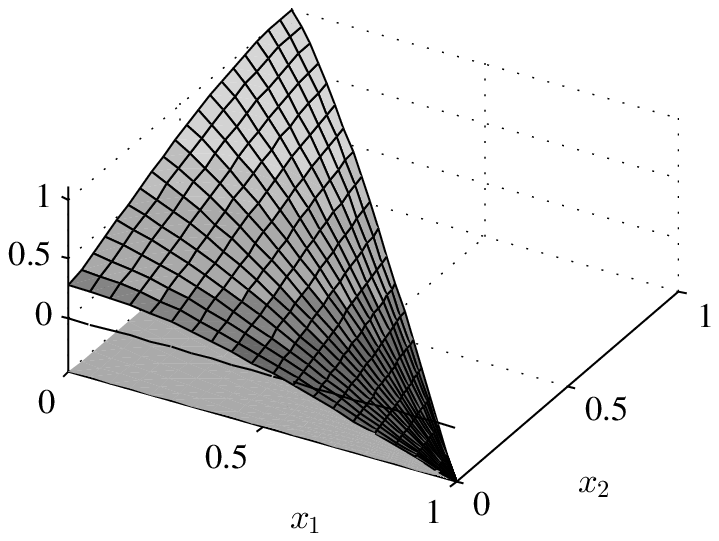}
	\label{fig:2d-up-eigenfunctions-pi-2}}	
	\subfloat[exact \; $u^{{\rm p}}_{\Gamma}$, $\alpha = \tfrac{\pi}{2}$]{\includegraphics[scale=0.9]{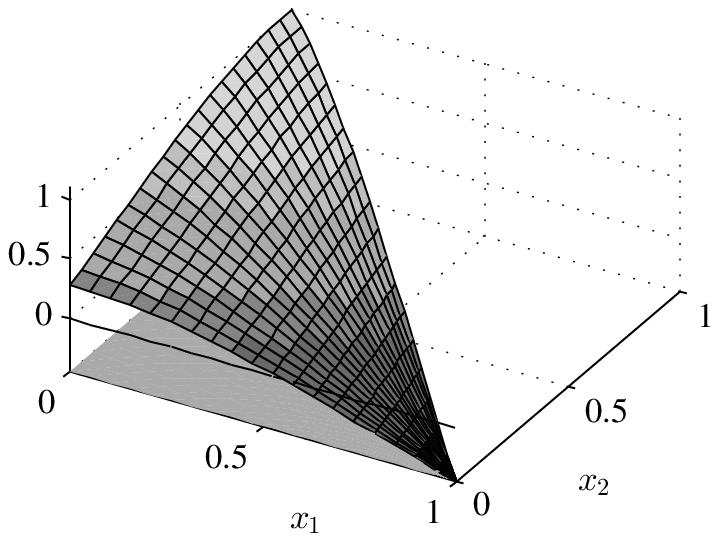}
	\label{fig:2d-exact-up-eigenfunctions-pi-2}}\\
	\subfloat[$u^{48, {\rm p}}_{\Gamma}$, $\alpha = \tfrac{2\pi}{3}$]{\includegraphics[scale=0.9]{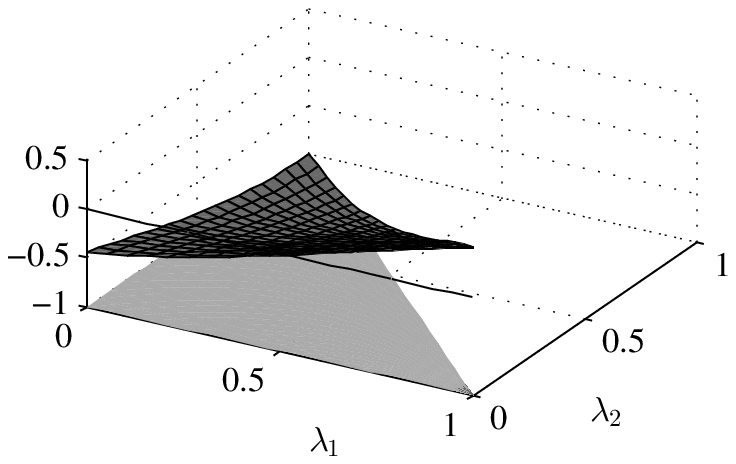}
	\label{fig:2d-up-eigenfunctions-3pi-4}}
	\subfloat[$u^{48, {\rm p}}_{\Gamma}$, $\alpha = \tfrac{3\pi}{4}$]{\includegraphics[scale=0.9]{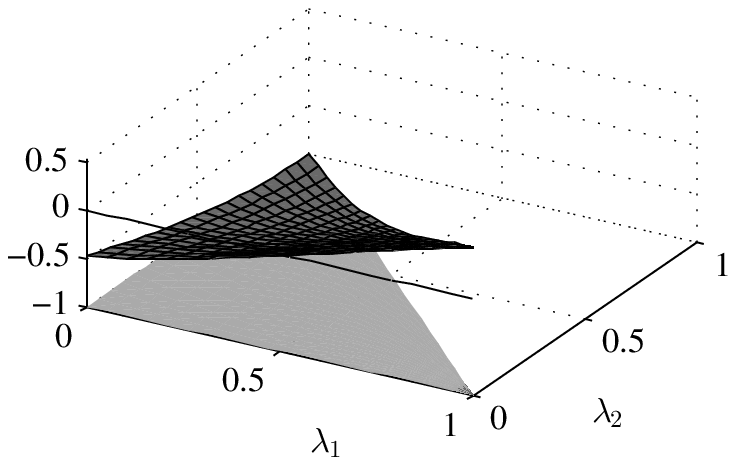}
	\label{fig:2d-up-eigenfunctions-5pi-6}}
	\caption{Eigenfunctions corresponding to $\approxCPT$ and 
	for $M = 48$ on simplex $\T$ with $\rho = 1$ and different $\alpha$.}	
	\label{fig:2d-eigenfunctions-cp}
\end{figure}
 
\begin{figure}[!ht]
	\centering
	\subfloat[$u^{48, {\rm Tr}}_{\Gamma}$, $\alpha = \tfrac{\pi}{6}$]{\includegraphics[scale=0.9]{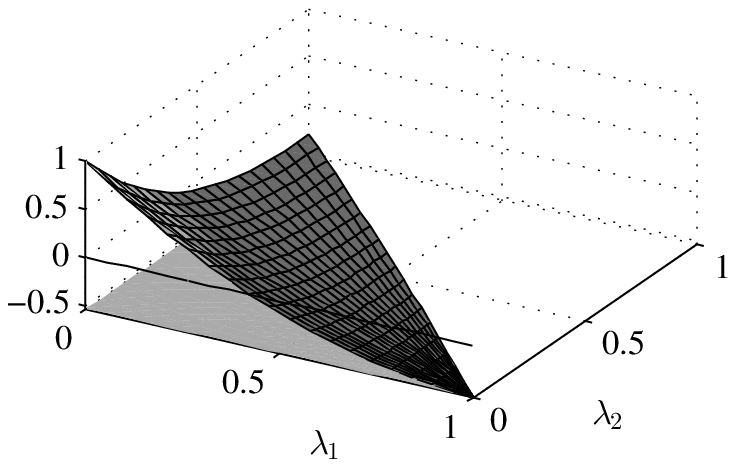}
	\label{fig:2d-ug-eigenfunctions-pi-6}
	} 
	\subfloat[$u^{48, {\rm Tr}}_{\Gamma}$, $\alpha = \tfrac{\pi}{3}$]{\includegraphics[scale=0.9]{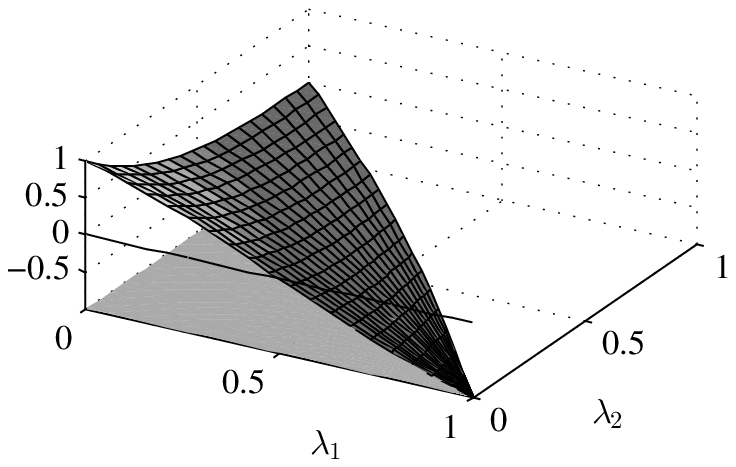}
	\label{fig:2d-ug-eigenfunctions-pi-3}}
	\\
	\subfloat[$u^{48, {\rm Tr}}_{\Gamma}$, $\alpha = \tfrac{\pi}{2}$]{\includegraphics[scale=0.9]{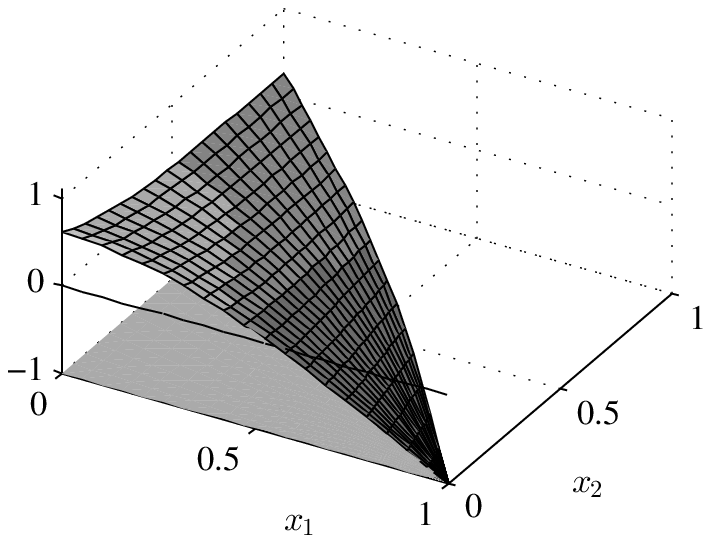}
	\label{fig:2d-ut-eigenfunctions-pi-2}
	}	
	\subfloat[exact \; $u^{{\rm Tr}}_{\Gamma}$, $\alpha = \tfrac{\pi}{2}$]{\includegraphics[scale=0.9]{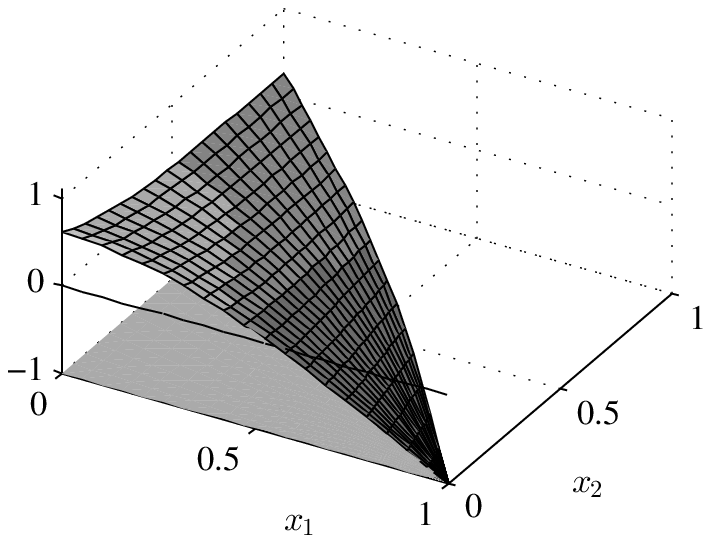}
	\label{fig:2d-exact-ut-eigenfunctions-pi-2}
	}\\	
	\subfloat[$u^{48, {\rm p}}_{\Gamma}$, $\alpha = \tfrac{2\pi}{3}$]{\includegraphics[scale=0.9]{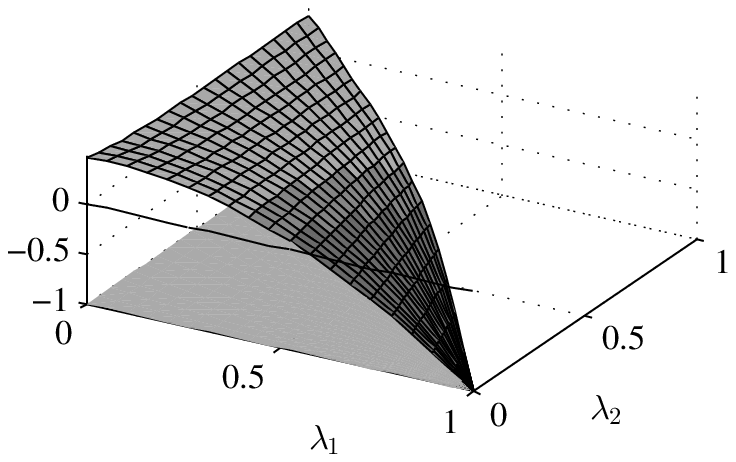}
	\label{fig:2d-ug-eigenfunctions-3pi-4}}
	\subfloat[$u^{48, {\rm p}}_{\Gamma}$, $\alpha = \tfrac{3\pi}{4}$]{\includegraphics[scale=0.9]{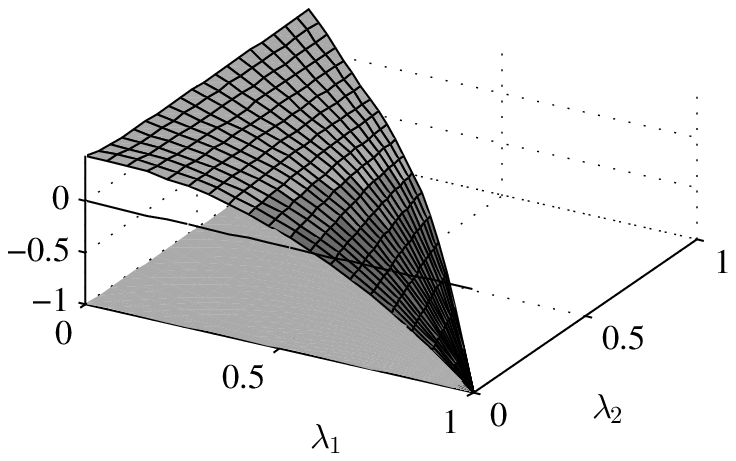}
	\label{fig:2d-ug-eigenfunctions-5pi-6}}
	\\
	\caption{Eigenfunctions corresponding to
	$\approxCGT$ for $M = 48$ on simplex $\T$ with $\rho = 1$ and different $\alpha$.}	
	\label{fig:2d-eigenfunctions-ctr}
\end{figure}

%%--------------------------------------------------------------------------------------%
%\begin{equation*}
    %\begin{matrix}
    %\begin{pmatrix}
    %x\\
    %y\\
    %\end{pmatrix}
    %=
    %\begin{pmatrix}
    %x_1 & x_2 & x_3 \\
    %y_1 & y_2 & y_3 \\
    %\end{pmatrix}
    %\,
    %(\lambda_1, \lambda_2, \lambda_3)^{\rm T},
    %\end{matrix}
    %\quad \mathrm{and} \quad
    %\begin{pmatrix}
    %\lambda_{1}\\
    %\lambda_{2}\\
    %\end{pmatrix}
    %= \mathbf{B}^{-1} (\mathbf{r} - \mathbf{r}_3),
%\end{equation*}
%
%where $\mathbf{r}_i = (x_i, y_i), \; i = 1, 2, 3,$ are the vertexes of $T$ and
%$\mathbf{r} = (x, y) \in T$. 
Figs. \ref{fig:2d-eigenfunctions-cp} and \ref{fig:2d-eigenfunctions-ctr} 
show the eigenfunctions computed for isosceles triangles with different angles $\alpha$
between two legs (zero mean condition is imposed on one of the legs).
The eigenfunctions have been computed in the process of finding $\approxCPT$ and
$\approxCGT$. The eigenfunctions are normalized so that the maximal
value is equal to $1$.
For $\alpha = \tfrac{\pi}{2}$, the exact eigenfunction associated 
with  the smallest positive eigenvalue
$\lambda^{\mathrm{P}}_{\Gamma} = \big( \tfrac{z_0}{h}\big)^{2}$
is known  (see \cite{NazarovRepin2014}):
\begin{equation*}
    u^{\rm P}_{\Gamma}
    = \cos(\tfrac{\zeta_0 x_1}{h}) + \cos \big(\tfrac{\zeta_0 (x_2 - h)}{h}\big).
    \label{eq:u-pt-exact}
\end{equation*}
Here, $\zeta_0$ is the root of the first equation in \eqref{eq:roots} 
(see Fig
\ref{fig:2d-exact-up-eigenfunctions-pi-2}).
We can compare $u^{\rm P}_{\Gamma}$ with the approximate eigenfunction $u^{M, \mathrm{P}}_{\Gamma}$ 
computed by minimization of
$\mathcal{R}^{\mathrm{P}}_{\Gamma} [w]$ (this function is depicted in Fig. \ref{fig:2d-up-eigenfunctions-pi-2}). 

Eigenfunctions related to the constant 
$\approxCGT$ are presented in Fig.  \ref{fig:2d-eigenfunctions-ctr}.
Again, for  $\alpha = \tfrac{\pi}{2}$ we know the exact eigenfunction
\begin{equation*}
    u^{\rm Tr}_{\Gamma} =
    \cos (\hat{\zeta}_0 x_1) \, \cosh \big(\hat{\zeta}_0 (x_2 - h)\big) +
    \cosh (\hat{\zeta}_0 \, x_1) \ \cos \big(\hat{\zeta}_0 (x_2 - h)\big),
    \label{eq:u-gt-exact}
\end{equation*}
%
%--------------------------------------------------------------------------------------%
where $\hat{\zeta}_0$ is the root of second equation in \eqref{eq:roots}
(see Fig. \ref{fig:2d-exact-ut-eigenfunctions-pi-2}). This function 
minimizes the quotient $\mathcal{R}^{\mathrm{Tr}}_{\Gamma} [w]$ and yields
the smallest positive eigenvalue
$\lambda^{\mathrm{Tr}}_{\Gamma} = \tfrac{\hat{\zeta}_0 \tanh (\hat{\zeta}_0)}{h}$.
It is easy to see that  numerical approximation
$\approxCGT$ (for $M(N) = 48$)
practically coincides with the exact function.
%

%----------------------------------------------------------------------------------------------% rho = 1
%----------------------------------------------------------------------------------------------%
\begin{figure}[!ht]
	\centering
	%----------------------------------------------------------------------------%
	% alpha = pi/3 - eps
	%----------------------------------------------------------------------------%
	\subfloat[ $u^{48}_{\T, 1}$, 
	$\alpha = \tfrac{\pi}{3} - \varepsilon$
	%, 
	%$\approxC = 0.2419$
	%$\lambda_1 = 17.0951$
	]{
	\includegraphics[scale=0.7]{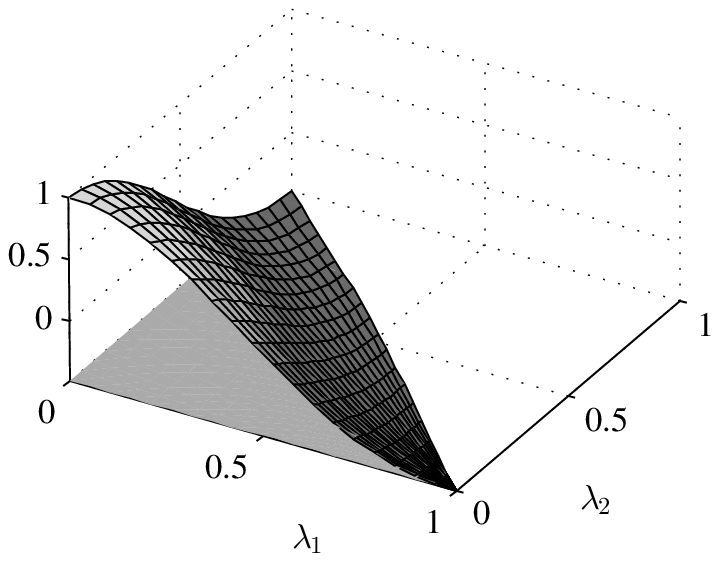}
	\label{fig:u1-rho-1-pi-3-minus-eps}}
	\subfloat[$u^{48}_{\T, 2}$, 
	$\alpha = \tfrac{\pi}{3} - \varepsilon$	
	%, 
	%$\approxC = 0.2229$
	%$\lambda_3 = 20.1216$
	]{
	\includegraphics[scale=0.7]{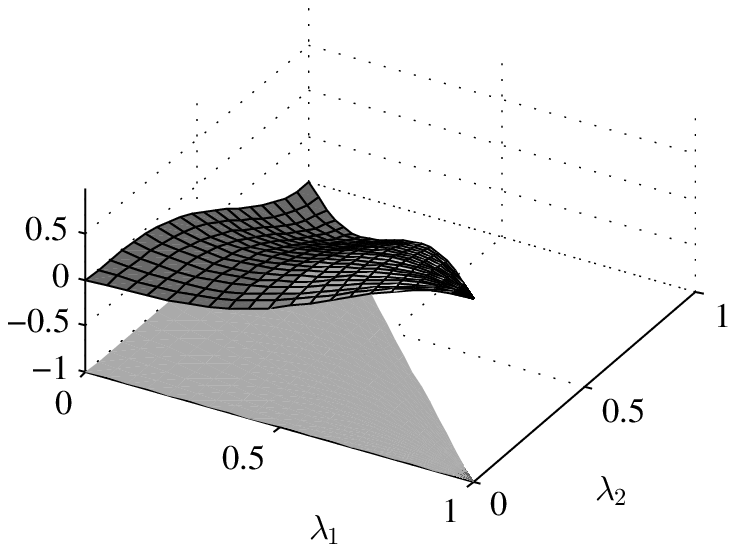}
	\label{fig:u2-rho-1-pi-3-minus-eps}}
	\subfloat[$u^{48}_{\T, 3}$, 
	$\alpha = \tfrac{\pi}{3} - \varepsilon$	
	%, 
	%$\approxC = 0.1353$
	%$\lambda_3 =  54.6024$
	]{
	\includegraphics[scale=0.7]{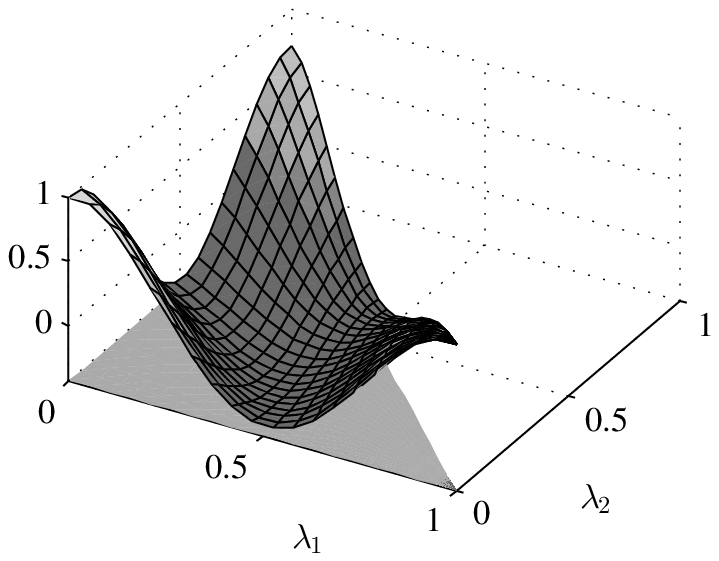}
	\label{fig:u3-rho-1-pi-3-minus-eps}}\\
	%----------------------------------------------------------------------------%
	% alpha = pi/3
	%----------------------------------------------------------------------------%
	\subfloat[	$u^{48}_{\T, 1}$, 
	$\alpha = \tfrac{\pi}{3}$
	%, 
	%$\approxC = 0.2387$
	%$\lambda_1 = 17.5463$
	]{
	\includegraphics[scale=0.7]{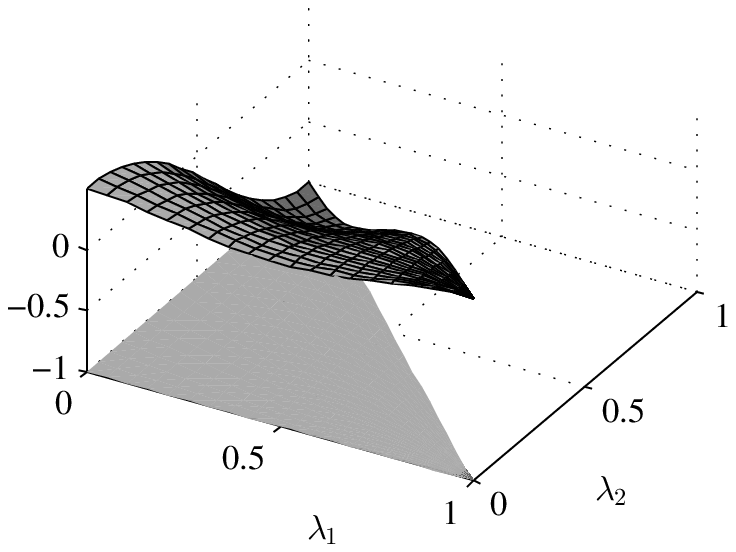}
	\label{fig:u1-rho-1-pi-3}}
	\subfloat[$u^{48}_{\T, 2}$, 
	$\alpha = \tfrac{\pi}{3}$
	%, 
	%$\approxC = 0.2387$
	%$\lambda_2 = 17.5463$
	]{
	\includegraphics[scale=0.7]{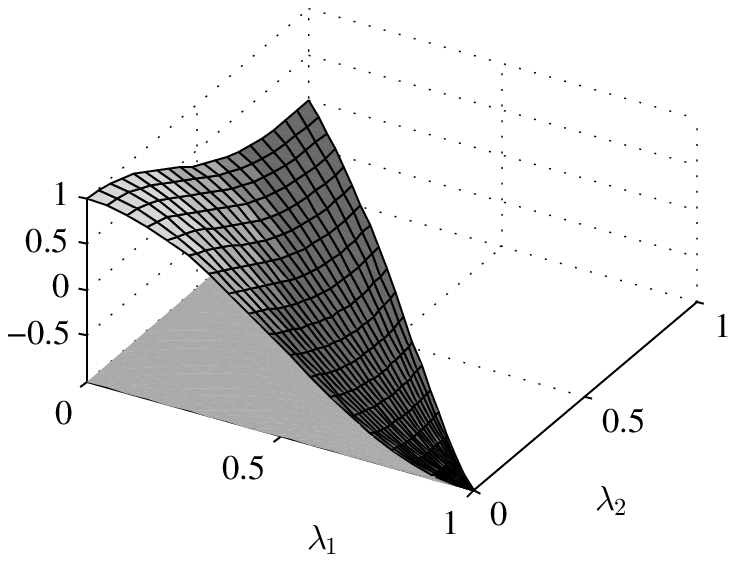}
	\label{fig:u2-rho-1-pi-3}}
	\subfloat[$u^{48}_{\T, 3}$, 
	$\alpha = \tfrac{\pi}{3}$	
	%, 
	%$\approxC = 0.1378$
	%$\lambda_3 =  52.6396$
	]{
	\includegraphics[scale=0.7]{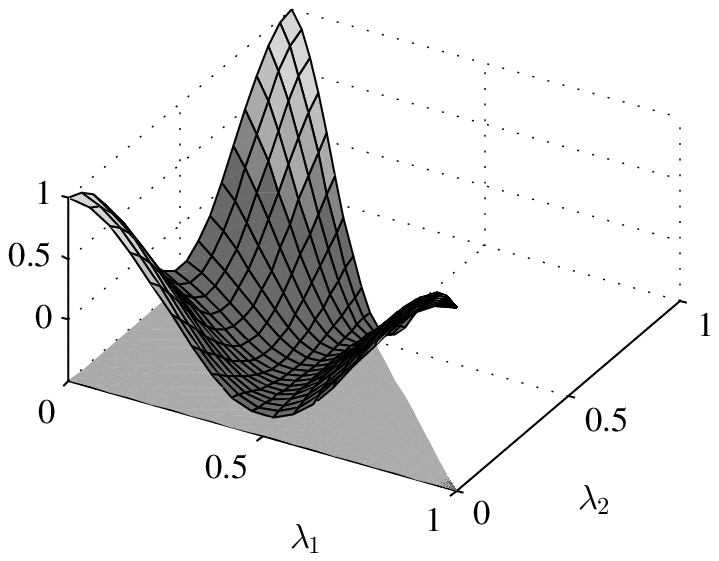}
	\label{fig:u3-rho-1-pi-3}}\\
	%----------------------------------------------------------------------------%
	% alpha = pi/3 + eps
	%----------------------------------------------------------------------------%
	\subfloat[$u^{48}_{\T, 1}$, 
	$\alpha = \tfrac{\pi}{3} + \varepsilon$	
	%, 
	%$\approxC = 0.2537$
	%$\lambda_1 =  15.5404$
	]{
	\includegraphics[scale=0.7]{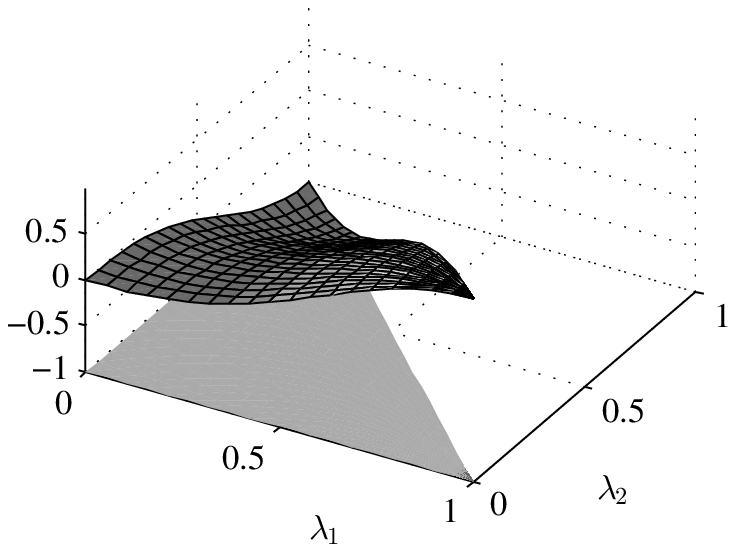}
	\label{fig:u1-rho-1-pi-3-plus-eps}}
	\subfloat[$u^{48}_{\T, 2}$, 
	$\alpha = \tfrac{\pi}{3} + \varepsilon$	
	%, 
	%$\approxC = 0.2355$
	%$\lambda_2 =  18.0309$
	]{
	\includegraphics[scale=0.7]{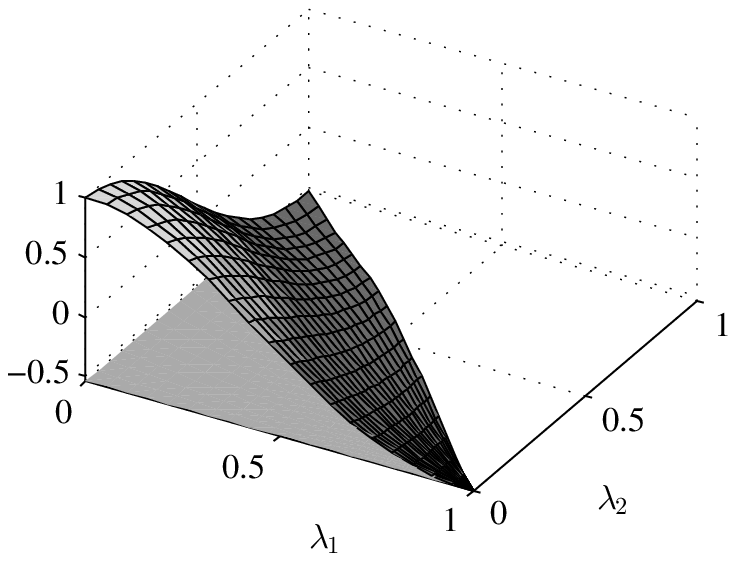}
	\label{fig:u2-rho-1-pi-3-plus-eps}}
	\subfloat[$u^{48}_{\T, 3}$, 
	$\alpha = \tfrac{\pi}{3} + \varepsilon$	
	%, 
	%$\approxC = 0.1422$
	%$\lambda_2 = 49.4818$
	]{
	\includegraphics[scale=0.7]{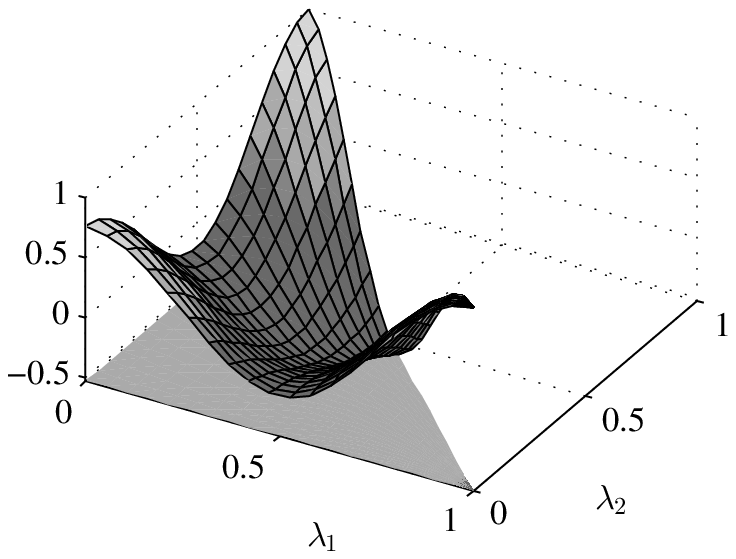}
	\label{fig:u3-rho-1-pi-3-plus-eps}}\\
	\caption{Eigenfunctions corresponding to $\approxC$ with $M = 48$
	on isosceles triangles $\T \in \Rtwo$ with $\alpha = \tfrac{\pi}{3}$, 
	$\tfrac{\pi}{3} - \varepsilon$, and 
	$\tfrac{\pi}{3} + \varepsilon$ in barycentric coordinates.}	
	\label{fig:u1-u2-u3-rho-1}	
\end{figure}

%--------------------------------------------------------------------------------------%
\begin{table}[!ht]
\centering
\footnotesize
\begin{tabular}{c|c|cc|cc|cc}
$$ & $$
& \multicolumn{2}{c|}{ $\tfrac{\pi}{3} - \varepsilon$ }
& \multicolumn{2}{c|}{ $\tfrac{\pi}{3}$ }
& \multicolumn{2}{c}{ $\tfrac{\pi}{3} + \varepsilon$ }\\
\midrule
$$ & $u^{M}_{\T, i}$ &
$\underline{C}^{48}_{\T, i}$ & $\lambda^{48}_{\T, i}$ &
$\underline{C}^{48}_{\T, i}$ & $\lambda^{48}_{\T, i}$ &
$\underline{C}^{48}_{\T, i}$ & $\lambda^{48}_{\T, i}$ \\
\midrule
\multirow{3}{*}{$\rho = 1$}
& $u^{48}_{\T, 1}$ & 0.2419 & 17.0951 & 0.2387 & 17.5463 & 0.2537 & 15.5404 \\
& $u^{48}_{\T, 2}$ & 0.2229 & 20.1216 & 0.2387 & 17.5463 & 0.2355 & 18.0309 \\
& $u^{48}_{\T, 3}$ & 0.1353 & 54.6024 & 0.1378 & 52.6396 & 0.1422 & 49.4818 \\
\midrule
\multirow{3}{*}{$\rho = \tfrac{\sqrt{2}}{2}$}
& $u^{48}_{\T, 1}$ & 0.23137 & 18.6804 & 0.23671 & 17.8471 & 0.24336 & 16.8850 \\
& $u^{48}_{\T, 2}$ & 0.17082 & 34.2707 & 0.17435 & 32.8970 & 0.17642 & 32.1295 \\
& $u^{48}_{\T, 3}$ & 0.1229  & 66.2058 & 0.12789 & 61.1402 & 0.13298 & 56.5493 \\
\midrule
\multirow{3}{*}{$\rho = \tfrac{3}{2}$}
& $u^{48}_{\T, 1}$ & 0.34714 & 8.2983 & 0.35523 & 7.9247 & 0.3648 & 7.5143 \\
& $u^{48}_{\T, 2}$ & 0.24485 & 16.6801 & 0.24885 & 16.1482 & 0.25125 & 15.8412 \\
& $u^{48}_{\T, 3}$  & 0.18258 & 29.9981 & 0.19084 & 27.4575 & 0.19845 & 25.3921 \\
\end{tabular}
\caption{$\approxC$ and $\lambda^{M}_{\T}$ corresponding to the first three 
eigenfunctions in Fig. \ref{fig:u1-u2-u3-rho-1}.}
\label{tab:approx-3-eigenvalues-and-constants}
\end{table}

Typically, the eigenfunctions associated with minimal positive
eigenvalues expose a continuous evolution with respect to $\alpha$.
However, this is not true for the quotient $\mathcal{R}_{\T}[w]$,
where the minimizer radically changes the profile. 
Fig. \ref{fig:cp-rho-1} indicates a possibility of such a rapid change
at $\alpha = \tfrac{\pi}{3}$, where the curve (related to $\underline{C}^{48}_{\T}$) 
obviously becomes
non-smooth. This happens because an equilateral triangle has double eigenvalue, therefore 
the minimizer of $\mathcal{R}_{\T}[w]$ over
$V^N_1$ changes its profile. 
Figs. \ref{fig:u1-rho-1-pi-3-minus-eps}--\ref{fig:u3-rho-1-pi-3-plus-eps} show
three eigenfunctions $u^{48}_{\T, 1}$, $u^{48}_{\T, 2}$,
and $u^{48}_{\T, 3}$ corresponding to three minimal eigenvalues $\lambda^{48}_{\T, 1}$, 
$\lambda^{48}_{\T, 2}$, and $\lambda^{48}_{\T, 3}$.
All functions are computed
for isosceles triangles and  are sorted in accordance with increasing 
values of the respective eigenvalues.
%Fig. \ref{fig:u1-u2-u3-rho-1} illustrates these three eigenfunctions for $\T$ with angles
%$\alpha = \tfrac{\pi}{3}$, $\tfrac{\pi}{3} + \eps$, and $\tfrac{\pi}{3} - \eps$, where 
%$\eps = \tfrac{\pi}{36}$.
It is easy to see that at $\alpha=\tfrac{\pi}{3}$ the first and the second eigenfunctions 
swap places.
Table \ref{tab:approx-3-eigenvalues-and-constants} presents the corresponding
results in the digital form.

It is worth noting that for equilateral triangles two minimal eigenfunctions
are known (see \cite{McCartin2002}):
\begin{alignat*}{2}
    u_1 & = \cos \Big(\tfrac{2\,\pi}{3} (2\, x_1 - 1)\Big)
    - \cos \Big(\tfrac{2\,\pi}{\sqrt{3}} x_2 \Big)
        \cos \Big(\tfrac{\pi}{3} (2\,x_1 - 1)\Big), \\
    u_2 & = \sin \Big(\tfrac{2\,\pi}{3} (2\, x_1 - 1)\Big)
    + \cos \Big(\tfrac{2\,\pi}{\sqrt{3}} x_2 \Big)
        \sin \Big(\tfrac{\pi}{3} (2\,x_1 - 1)\Big).
\end{alignat*}
These functions practically coincide with the functions $u^{48}_{\T, 1}$ and $u^{48}_{\T, 2}$ 
presented in Fig. \ref{fig:u1-rho-1-pi-3}.
%
%--------------------------------------------------------------------------------------%
Finally, we note that this phenomenon (change of the minimal eigenfunction) does not 
appear for $\rho = \tfrac{\sqrt{2}}{2}$ or $\rho=\tfrac{3}{2}$. 
%due to the fact 
%that non-quadrilateral triangles have simple lowest eigenvalue.
%(see Figs
%\ref{fig:u1-u2-u3-rho-sqrt22} and \ref{fig:u1-u2-u3-rho-32}).
The eigenvalues as well
as the constants corresponding to the eigenfunctions presented in Fig.
\ref{fig:u1-u2-u3-rho-1} %, \ref{fig:u1-u2-u3-rho-sqrt22} and \ref{fig:u1-u2-u3-rho-32}
are shown in the Table \ref{tab:approx-3-eigenvalues-and-constants}.

%======================================================================================%
%======================================================================================%
\section{Two-sided bounds of $\CPGamma$ and $\CtrGamma$ for tetrahedrons}
\label{eq:numerial-tests-3d}

\begin{figure}[ht!]
\centering
	\adjustbox{valign=b}{
	\begin{minipage}[b]{0.5\linewidth}
	\begin{tikzpicture}[scale=0.6]

\coordinate [label=below:$A$] (A) at (0,0,0);
\coordinate [label=above right:$$] (B) at (5,0,0);
\coordinate [label=below:$$] (C) at (0,0,7);
\coordinate [label=above left:$$] (D) at (1,5,3.5);
%\coordinate [label=below:$H$] (H) at (1,0,3.5);
\coordinate [label=above left:$$] (N) at (1,5,0);
%\coordinate [label=above:$M$] (M) at (0,5,3.5);
%\coordinate [label=below:$$] (L) at (0,0,3.5);
%\coordinate [label=below:$$] (O) at (1,0,0);
%\coordinate [label=below:$$] (K) at (0,5,0);

\draw[-, fill=black!10, opacity=.5, dashed] (A)--(B)--(C)--(A);

\draw (2.4, 1.2, 1) node[anchor=north west] {$\T$};
\draw (2, 0.0, 0.7) node[anchor=north east] {${\Gamma}$};

\draw[->] (B) -- (5.5,0,0) node[right] {$x_1$};
\draw[->] (0,0,0) -- (0,6,0) node[above] {$x_2$};
\draw[->] (C) -- (0,0,8) node[below] {$x_3$};
%
%\draw (4.85,-0.15) .. controls (4.85, 0.45) and (4.95, 0.45) .. node[right] {$$}  (5.11, 0.55);

% theta angle
\draw[->] (0.0, 0.7, 0.2) .. controls (0.0, 0.7, 0.7) and (0.0, 0.6, 1.0) .. node[left] {$\theta$}  (0.0, 0.0, 1.2);
%\draw (0.9, 0.5) node[anchor=north east] {$\theta$};

% alpha angle
\draw[->] (0.8,0.0) .. controls (0.7, 0.7) and (0.2, 0.6) .. node[right] {$\;\alpha$}  (0.1,0.6);

\draw (0,0.01,5) node[anchor=north west] {$C (0, 0, h_3)$};
\draw (4,0,0) node[anchor=north west] {$B (h_1, 0, 0)$};
%\draw (1,5,3.5) node[anchor=south west] {$D (h_2 \, \sin(\phi)\, \cos(\alpha), h_2 \, \sin(\phi)\, \sin(\alpha), h_2\, \cos(\alpha))$};
\draw (1,5,3.5) node[anchor=south east] {$D$};

%\draw[-, opacity=.5] (A)--(D);
\draw[-, opacity=.5] (B)--(C);
\draw[-, opacity=.5] (C)--(D);
\draw[-, opacity=.5] (B)--(D);

\draw (D)--(A);
\draw[dashed] (D)--(N);
%\draw[dashed] (D)--(M);
%\draw[dashed] (D)--(H);
%\draw[dashed] (D)--(L);
%
\draw[dashed] (N)--(A);
%\draw[dashed] (N)--(A);
%
%\draw[dashed] (H)--(L);
%\draw[dashed] (H)--(O);
%\draw[dashed] (M)--(L);
%\draw[dashed] (N)--(O);
%\draw[dashed] (K)--(N);
%\draw[dashed] (K)--(M);

\filldraw [black] (0,0,0) circle (0.8pt)
                  (5,0,0) circle (0.8pt)
                  (0,0,7) circle (0.8pt)
									(1,5,3.5) circle (0.8pt);

\end{tikzpicture}
	\caption{Simplex in $\Rthree$.}
	\label{eq:3d-simplex}
	\end{minipage}}
%-----------------------------------------------------------------------------------%
	\adjustbox{valign=b}{
	\begin{minipage}[b]{0.45\linewidth} 
	\begin{tikzpicture}[scale=0.8]

\coordinate [label=below:$A$] (A) at (0,0,0);
%\coordinate [label=above right:$$] (B) at (5,0,0);
%\coordinate [label=below:$$] (C) at (0,0,7);
\coordinate [label=above left:$$] (D) at (1,5,3.5);
\coordinate [label=below:$H$] (H) at (1,0,3.5);
\coordinate [label=above left:$N$] (N) at (1,5,0);
\coordinate [label=above:$M$] (M) at (0,5,3.5);
\coordinate [label=below:$$] (L) at (0,0,3.5);
\coordinate [label=below:$$] (O) at (1,0,0);
\coordinate [label=below:$$] (K) at (0,5,0);

%\draw[-, fill=black!10, opacity=.5, dashed] (A)--(B)--(C)--(A);

%
%\draw (2.4, 1.2, 1) node[anchor=north west] {$\T_{\theta, \alpha}$};
%\draw (2, 0.0, 0.7) node[anchor=north east] {${\Gamma}$};

\draw[->] (A) -- (2,0,0) node[right] {$x_1$};
\draw[->] (A) -- (0,6,0) node[above] {$x_2$};
\draw[->] (A) -- (0,0,5) node[below] {$x_3$};
%

% theta angle
\draw[->] (0.0, 0.7, 0.2) .. controls (0.0, 0.7, 0.7) and (0.0, 0.6, 1.0) .. node[left] {$\theta$}  (0.0, 0.0, 1.2);
%\draw (0.9, 0.5) node[anchor=north east] {$\theta$};

% alpha angle
\draw[->] (0.8,0.0) .. controls (0.7, 0.7) and (0.2, 0.6) .. node[right] {$\;\alpha$}  (0.1,0.6);

%\draw (0,0.01,5) node[anchor=north west] {$C (0, 0, h_3)$};
%\draw (4,0,0) node[anchor=north west] {$B (h_1, 0, 0)$};
%\draw (1,5,3.5) node[anchor=south west] {$D (h_2 \, \sin(\phi)\, \cos(\alpha), h_2 \, \sin(\phi)\, \sin(\alpha), h_2\, \cos(\alpha))$};
\draw (1,5,3.5) node[anchor=west] {$\quad D (D_x, D_y, D_z)$};

%\draw[-, opacity=.5] (A)--(D);
%\draw[-, opacity=.5] (B)--(C);
%\draw[-, opacity=.5] (C)--(D);
%\draw[-, opacity=.5] (B)--(D);

\draw (D)--(A);
\draw[dashed] (D)--(N);
\draw[dashed] (D)--(M);
\draw[dashed] (D)--(H);
\draw[dashed] (D)--(L);

\draw[dashed] (N)--(A);
\draw[dashed] (N)--(A);

\draw[dashed] (H)--(L);
\draw[dashed] (H)--(O);
\draw[dashed] (M)--(L);
\draw[dashed] (N)--(O);
\draw[dashed] (K)--(N);
\draw[dashed] (K)--(M);

%\draw[-] (5.8,0,1.5)--(5.8,0.2,1.5)--(6,0.2,1.5);
%\draw[dashed] (5.5, 0.0,1.5)--(6.5,0.0,1.5);
%\draw[dashed] (A)--(H);
%\draw[dashed] (H)--(M);
%\draw[dashed] (N)--(M);
%\draw[dashed] (N)--(B);
%\draw[dashed] (N)--(A);
%\draw[dashed] (L)--(D);
%\draw[dashed] (H)--(O);
%\draw[dashed] (D)--(O);

%\draw[-] (5.8,0,1.5)--(5.8,0.2,1.5)--(6,0.2,1.5);
%\draw[dashed] (5.5, 0.0,1.5)--(6.5,0.0,1.5);

% pi/2 angle 
%\draw[-] (5,0,2.6)--(5,0.3,2.1)--(5,0.3,1.3);

%
%\foreach \i in {B,C}
    %\draw[dashed] (A)--(\i);
%

\filldraw [black] (0,0,0) circle (0.8pt)
									(1,5,3.5) circle (0.8pt);

\end{tikzpicture}
	\caption{Coordinate of the vertex $D$.}
	\label{eq:3d-simplex-d-coordinate}
	\end{minipage}}
\end{figure}
We orient the coordinates as it is shown in Fig. \ref{eq:3d-simplex}
and define a non-degenerate simplex in $\Rthree$ with vertexes \linebreak 
$A = (0, 0, 0)$, $B = (h_1, 0, 0)$, $C = (0, 0, h_3)$, and 
$D = (h_2 \, \sin\theta\, \cos\alpha, h_2 \, \sin\theta\, \sin\alpha, h_2\, \cos\alpha)$,
%
%\begin{equation}
    %\T = {\rm conv} \big \{(0, 0, 0), (h_1, 0, 0), (0, 0, h_3), (D_{x_1}, D_{x_2}, D_{x_3}) \big\},
    %\label{eq:arbitrary-tetrahedron}
%\end{equation}
%%
%where
%$(D_{x_1}, D_{x_2}, D_{x_3}) =
%(h_2 \, \sin\theta\, \cos\alpha, h_2 \, \sin\theta\, \sin\alpha, h_2\, \cos\alpha)$, 
where
$h_1$ and $h_3$ are the scaling parameters along axis 
$O_{x_1}$ and $O_{x_3}$, respectively, $AD = h_2$,
$\alpha$ is a polar angle, and $\theta$ is an azimuthal angle
(see Fig. \ref{eq:3d-simplex-d-coordinate}).  
Let $\Gamma$ be defined by vertexes $A$, $B$, and $C$.
%
%--------------------------------------------------------------------------------------%

To the best of our knowledge, exact values of constants in Poincar\'{e}-type 
inequalities for simplexes in $\Rthree$ are unknown. Therefore,  we first consider 
four basic (reference) tetrahedrons with 
$h_2 = 1$, $\hat{\theta} = \tfrac{\pi}{2}$, and 
$\hat{\alpha}_1 = \tfrac{\pi}{4}$, 
$\hat{\alpha}_2 = \tfrac{\pi}{3}$, 
$\hat{\alpha}_3 = \tfrac{\pi}{2}$, and 
$\hat{\alpha}_4 = \tfrac{2\pi}{3}$. 
The respective constants are found numerically with high accuracy (see Table 
\ref{tab:const-convergence-from-basis-G-leg-alpha}, 
which shows convergence of the constants 
with respect to increasing $M(N)$).
%
%--------------------------------------------------------------------------------------%
Henceforth, $\Tref_{\hat{\theta}, \hat{\alpha}}$ denotes a reference tetrahedron, 
where $\hat{\theta}$ and $\hat{\alpha}$ are certain fixed angles.
By $\mathcal{F}_{\hat{\theta}, \hat{\alpha}}$ we denote the respective mapping 
$\mathcal{F}_{\hat{\theta}, \hat{\alpha}}: 
\Tref_{\hat{\theta}, \hat{\alpha}} \rightarrow \T$.
\begin{table}[!t]
\centering
\footnotesize
\begin{tabular}{c|cc|cc|cc|cc}
\multicolumn{1}{c|}{$ $}
& \multicolumn{2}{c|}{ $\hat{\alpha} = \tfrac{\pi}{4}$}
& \multicolumn{2}{c|}{ $\hat{\alpha} = \tfrac{\pi}{3}$}
& \multicolumn{2}{c|}{ $\hat{\alpha} = \tfrac{\pi}{2}$}
& \multicolumn{2}{c}{ $\hat{\alpha} = \tfrac{2\pi}{3}$}\\
\midrule
\multicolumn{1}{c|}{$ $} &
\multicolumn{2}{c|}{
\begin{tikzpicture}[line join = round, line cap = round]
\coordinate [label=left:] (A) at (0,0,0);
\coordinate [label=above right:] (B) at (1.2,0,0);
\coordinate [label=below:] (C) at (0,0,1.2);
\coordinate [label=above right:] (D) at ( 0.848, 0.848, 0);
\draw[->] (1.2, 0,   0) --   (1.6, 0, 0) node[below] {$\hat{x}_1$};
\draw[->] (0,   0,   0) --   (0, 1.4, 0) node[above] {$\hat{x}_2$};
\draw[->] (0,   0,   1.2) -- (0, 0, 1.6) node[right] {$\hat{x}_3$};
\draw (2.0, 1.05) node[anchor=north east] {$\Tref_{\rfrac{\pi}{4}}$};
\draw (1.0, -0.2) node[anchor=north east] {$\widehat{\Gamma}$};
\draw[-, fill=black!10, opacity=.5, dashed] (A)--(B)--(C)--(A);
\draw[-, opacity=.5] (B)--(C);
\draw[-, opacity=.5] (C)--(D);
\draw[-, opacity=.5] (B)--(D);
\draw (0.55,0.0) .. controls (0.5, 0.3) and (0.4, 0.3) .. node[right] {$\hat{\alpha}$}  (0.3,0.3);
\draw[dashed] (A)--(D);
\end{tikzpicture}
}
& \multicolumn{2}{c|}{
\begin{tikzpicture}[line join = round, line cap = round]
\coordinate [label=left:] (A) at (0,0,0);
\coordinate [label=above right:] (B) at (1.2,0,0);
\coordinate [label=below:] (C) at (0,0,1.2);
\coordinate [label=above right:] (D) at ( 0.6,  1.0392, 0);
\draw[->] (1.2, 0,   0) --   (1.6, 0, 0) node[below] {$\hat{x}_1$};
\draw[->] (0,   0,   0) --   (0, 1.4, 0) node[above] {$\hat{x}_2$};
\draw[->] (0,   0,   1.2) -- (0, 0, 1.6) node[right] {$\hat{x}_3$};
\draw (2.0, 1.05) node[anchor=north east] {$\Tref_{\rfrac{\pi}{3}}$};
\draw (1.0, -0.2) node[anchor=north east] {$\widehat{\Gamma}$};
\draw[-, fill=black!10, opacity=.5, dashed] (A)--(B)--(C)--(A);
\draw[-, opacity=.5] (B)--(C);
\draw[-, opacity=.5] (C)--(D);
\draw[-, opacity=.5] (B)--(D);
\draw (0.55,0.0) .. controls (0.5, 0.3) and (0.4, 0.3) .. node[right] {$\hat{\alpha}$}  (0.27,0.4);
\draw[dashed] (A)--(D);
\end{tikzpicture}
}
& \multicolumn{2}{c|}{
\begin{tikzpicture}[line join = round, line cap = round]
\coordinate [label=left:] (A) at (0,0,0);
\coordinate [label=above right:] (B) at (1.2,0,0);
\coordinate [label=below:] (C) at (0,0,1.2);
\coordinate [label=above right:] (D) at ( 0.0,  1.2, 0);
\draw[->] (1.2, 0,   0) --   (1.6, 0, 0) node[below] {$\hat{x}_1$};
\draw[->] (0,   1.2,   0) --   (0, 1.4, 0) node[above] {$\hat{x}_2$};
\draw[->] (0,   0,   1.2) -- (0, 0, 1.6) node[right] {$\hat{x}_3$};
\draw (2.0, 1.05) node[anchor=north east] {$\Tref_{\rfrac{\pi}{2}}$};
\draw (1.0, -0.2) node[anchor=north east] {$\widehat{\Gamma}$};
\draw[-, fill=black!10, opacity=.5, dashed] (A)--(B)--(C)--(A);
\draw[-, opacity=.5] (B)--(C);
\draw[-, opacity=.5] (C)--(D);
\draw[-, opacity=.5] (B)--(D);
\draw (0.3,0.0) .. controls (0.2, 0.3) and (0.1, 0.3) .. node[right] {$\hat{\alpha}$}  (0.0,0.3);
\draw[dashed] (A)--(D);
\draw[dashed] (0, 0, 0.5)--(0, 0.5, 0.5)--(0, 0.5, 0);
\draw (-0.1, 0.3) node[anchor=north east] {$\hat{\theta}$};

\end{tikzpicture}
}
& \multicolumn{2}{c}{
\begin{tikzpicture}[line join = round, line cap = round]
\coordinate [label=left:] (A) at (0,0,0);
\coordinate [label=above right:] (B) at (1.2,0,0);
\coordinate [label=below:] (C) at (0,0,1.2);
\coordinate [label=above right:] (D) at ( -0.6,  1.0392, 0);
\draw[->] (1.2, 0,   0) --   (1.6, 0, 0) node[below] {$\hat{x}_1$};
\draw[->] (0,   0,   0) --   (0, 1.4, 0) node[above] {$\hat{x}_2$};
\draw[->] (0,   0,   1.2) -- (0, 0, 1.6) node[right] {$\hat{x}_3$};
\draw (2.0, 1.05) node[anchor=north east] {$\Tref_{\rfrac{2\pi}{3}}$};
\draw (1.0, -0.2) node[anchor=north east] {$\widehat{\Gamma}$};
\draw[-, fill=black!10, opacity=.5, dashed] (A)--(B)--(C)--(A);
\draw[-, opacity=.5] (B)--(C);
\draw[-, opacity=.5] (C)--(D);
\draw[-, opacity=.5] (B)--(D);
\draw (0.3,0.0) .. controls (0.2, 0.3) and (0.1, 0.3) .. node[right] {$\hat{\alpha}$}  (-0.1,0.2);
\draw[dashed] (A)--(D);
\end{tikzpicture}
}\\
\midrule
$M(N)$ & $\CPTapproxhatalphahat$ & $\CGTapproxhatalphahat$
& $\CPTapproxhatalphahat$ & $\CGTapproxhatalphahat$
& $\CPTapproxhatalphahat$ & $\CGTapproxhatalphahat$
& $\CPTapproxhatalphahat$ & $\CGTapproxhatalphahat$\\
\midrule
%7   &  0.32430           & 0.76010       & 0.32599       & 0.65465       & 0.36053       & 0.6547 & 0.41521 &  0.68616 \\
%26  &  0.33854           & 0.82945       & 0.34027       & 0.76128       & 0.37367       & 0.7516 & 0.42748 &  0.86332 \\
%63  &  0.34112           & 0.83133       & 0.34256       & 0.76290       & 0.37559       & 0.7520 & 0.42864 &  0.86460 \\
%124 &  0.34115       & 0.83134       & 0.34259       & 0.76291       & 0.37560       & 0.7520 & 0.42867 &  0.86463 \\
%{\bf 215} &  {\bf 0.34115} & {\bf 0.83134} & {\bf 0.34259} & {\bf 0.76291} & {\bf 0.37560} & {\bf 0.75200} & {\bf 0.42867} &  {\bf 0.86463} \\
7   &  0.32431 & 0.760099 & 0.325985 & 0.654654 & 0.360532 & 0.654654 & 0.4152099 &  0.686161 \\
26  &  0.338539 & 0.829445 & 0.340267 & 0.761278 & 0.373669 & 0.751615 & 0.4274757 &  0.863324 \\
63  &  0.341122 & 0.831325 & 0.342556 & 0.762901 & 0.375590 & 0.751994 & 0.4286444 &  0.864595 \\
124 &  0.341147 & 0.831335 & 0.342589 & 0.762905 & 0.375603 & 0.751999 & 0.4286652 &  0.864630 \\
215 &  {\bf 0.341147} & {\bf 0.831335} & {\bf 0.342589} & {\bf 0.762905} & {\bf 0.375603} & {\bf 0.751999} & {\bf 0.4286652} &  {\bf 0.864630} \\
%---------------------------------------------------
%7   & 0.325985 & 0.654654 & 0.360532 & 0.654654 & 0.274407 & 0.830255 & 0.289421 & 0.717119 \\
%26  & 0.340267 & 0.761278 & 0.373669 & 0.751615 & 0.291311 & 0.930666 & 0.307103 & 0.840896 \\
%63  & 0.342556 & 0.762901 & 0.375590 & 0.751994 & 0.294844 & 0.931228 & 0.308907 & 0.840896 \\
%124 & 0.342589 & 0.762905 & 0.375603 & 0.751999 & 0.294872 & 0.931233 & 0.308924 & 0.840896 \\
%215 & 0.342589 & 0.762905 & 0.375603 & 0.751999 & 0.294872 & 0.931233 & 0.308924 & 0.840896 \\
\end{tabular}
\caption{$\CPTapproxhatalphahat$ and $\CGTapproxhatalphahat$ with
respect to $M(N)$ for $\Tref_{\hat{\theta}, \hat{\alpha}}$ with
$\rho = 1$, $\hat{\theta} = \tfrac{\pi}{2}$, and different 
$\hat{\alpha}$.}
\label{tab:const-convergence-from-basis-G-leg-alpha}
\end{table}
%
%--------------------------------------------------------------------------------------%
Then, for an arbitrary tetrahedron $\T$, we have
\begin{alignat}{2}
    \|v\|_{\T} \, & \leq\, \CPT \, h_2 \,\|\nabla v\|_{\T} \quad \mbox{and} \quad
    \|v\|_\Gamma \, \leq\, \CGT \, h_2^{\rfrac{1}{2}} \,\|\nabla v\|_{\T}
		\label{eq:general-poincare-type-inequalities-for-tetrahedron}
\end{alignat}
with approximate bounds
\begin{equation}
    \CPT \lessapprox \CPtetr = 
		\min_{\hat{\alpha} = \{\rfrac{\pi}{4}, \rfrac{\pi}{3}, \rfrac{\pi}{2}, \rfrac{2\pi}{3}\}} 
		\Big\{ \cpalphahat \, \CPThatalphahat \Big\}
\end{equation}		
and 
\begin{equation}
		\CGT  \lessapprox \CGtetr = 
		\min_{\hat{\alpha} = 
		\{\rfrac{\pi}{4}, \rfrac{\pi}{3}, \rfrac{\pi}{2}, \rfrac{2\pi}{3}\}} 
		\Big\{ \cgalphahat \, \CGThatalphahat \Big\},
\label{eq:}
\end{equation}
where $\CPThatalphahat$ and $\CGThatalphahat $ are the constants related to four reference tetrahedrons from Table \ref{tab:const-convergence-from-basis-G-leg-alpha}, and 
$\cpalphahat$ and $\cgalphahat$ (see \eqref{eq:ratio-3d}) are
generated by the mapping 
$\mathcal{F}_{\rfrac{\pi}{2}, \hat{\alpha}}$: 
$\Tref_{\rfrac{\pi}{2}, \hat{\alpha}} \rightarrow \T$. Here, the reference tetrahedrons are
defined based on $\hat{A} = (0, 0, 0)$, $\hat{B} = (1, 0, 0)$, 
$\hat{C} = (0, 0, 1)$, 
$\hat{D} = (\cos\hat{\alpha},  \sin\hat{\alpha}, 0)$ with $\hat{\alpha} = 
\{ \tfrac{\pi}{4},\tfrac{\pi}{3}, \tfrac{\pi}{2}, \tfrac{2\pi}{3}\}$,
and 
$\mathcal{F}_{\rfrac{\pi}{2}, \hat{\alpha}}(\hat{x})$ is presented by the relation
\begin{equation*}
    x = \mathcal{F}_{\rfrac{\pi}{2}, \hat{\alpha}}(\hat{x}) =
    B_{\rfrac{\pi}{2}, \hat{\alpha}} \hat{x}, \quad
    B_{\rfrac{\pi}{2}, \hat{\alpha}} = \{ b_{ij} \}_{i, j = 1, 2, 3} = h_2
    \begin{pmatrix}
        \: \tfrac{h_1}{h_2} & 
				\; \tfrac{\nu (\rho, \alpha)}{\sin \hat{\alpha}}  & 
				\; 0 \\[0.3em]
        \: 0 & 
				\; \, \tfrac{\sin\alpha \sin\theta}{\sin \hat{\alpha}} \: & 
				\: 0 \\[0.3em]
        \: 0 & 
				\; \, \tfrac{\cos\theta}{\sin \hat{\alpha}} \: & 
				\: \tfrac{h_3}{h_2} \\
    \end{pmatrix},
%\label{eq:3d-mapping}
\end{equation*}
%
%--------------------------------------------------------------------------------------%
where $\nu(\rho, \alpha) = \cos\alpha \sin\theta - \tfrac{h_1}{h_2} \cos \hat{\alpha}$, 
$\mathrm{det} \, B_{\rfrac{\pi}{2}, \hat{\alpha}} = h_1 \, h_2\, h_3 \, 
\tfrac{\sin\alpha \sin\theta}{\sin \hat{\alpha}}$.
By analogy with the two-dimensional case (see \eqref{eq:grad-hatv-lower-estimate-1}),
$\cpalphahat$ and $\cgalphahat$ depend on the maximum eigenvalue of the matrix 
\begin{alignat*}{2}
    A_{\rfrac{\pi}{2}, \hat{\alpha}} & := h_1^2
		\begin{pmatrix}
        \: b_{11}^2 + b^2_{12} \quad & b_{12} b_{22} \qquad & b_{12} b_{32} \\[0.5em]
      	b_{12} b_{22} & b^2_{22} & b_{22} b_{32} \\[0.5em]
				b_{12} b_{32} & b_{22} b_{32} & b_{33}^2 + b^2_{32} \\
    \end{pmatrix}
		.
		%\\[1em]
		%& = h_2^2 \cdot
		%\begin{pmatrix}
        %\: 1 + \tfrac{1}{\sin^2 \hat{\alpha}} \, 
				         %\nu^2 \quad &
				%
        %\tfrac{\rho \sin\alpha \sin\theta}{\sin^2 \hat{\alpha}} 
				%\nu \quad &
				%
        %\tfrac{\rho \cos\theta }{\sin^2 \hat{\alpha}} 
				%\nu \\[0.5em]
      	%\tfrac{ \rho \sin\alpha \sin\theta }{\sin^2 \hat{\alpha}} 
				%\nu \; &
        %\tfrac{ \rho^2 \, \sin^2 \alpha \, \cos^2 \theta}{\sin^2 \hat{\alpha}}  \; &
        %\tfrac{\rho^2 \sin\alpha \sin2\theta }{2 \sin^2 \hat{\alpha}}\\[0.5em]
        %\tfrac{\rho \cos\theta}{\sin^2 \hat{\alpha}} 
				%\nu \; &
        %\tfrac{\rho^2 \sin\alpha \sin2\theta }{2 \sin^2 \hat{\alpha}} \; &
        %\tfrac{h_3^2}{h_1^2} + \rho^2 \tfrac{\cos^2 \theta}{\sin^2 \hat{\alpha}}  \\
    %\end{pmatrix}.
    %\label{eq:a-3d-alpha-beta}
\end{alignat*}
%
%where $b_{12}$, $b_{22}$, $b_{32}$ are elements of $B_{\rfrac{\pi}{2}, \hat{\alpha}}$ %
%in \eqref{eq:3d-mapping}. 
%--------------------------------------------------------------------------------------%
The maximal eigenvalue of the matrix  $A_{\rfrac{\pi}{2}, \hat{\alpha}}$ is defined by 
the relation $\lambda_{\rm max} (A_{\rfrac{\pi}{2}, \hat{\alpha}}) 
=	h_2^2 \,  \mu_{\alpha, \theta, \hat{\alpha}}$ with
\begin{equation*}
	 \mu_{\alpha, \theta, \hat{\alpha}}  
	= \Big(\mathcal{E}_5^{\rfrac{1}{3}} 
	- \mathcal{E}_3 \mathcal{E}_5^{-\rfrac{1}{3}} 
	+ \tfrac{1}{3} \mathcal{E}_1 \Big),  
\end{equation*}
where
\begin{alignat*}{2}
	\mathcal{E}_1 & = b_{11}^2 + b_{12}^2 + b_{22}^2 + b_{32}^2 + b_{33}^2, \\
	\mathcal{E}_2 & = b_{11}^2 \, b_{22}^2 + b_{11}^2 \, b_{32}^2 + b_{11}^2 \, b_{33}^2 + b_{12}^2 \, b_{33}^2 + b_{22}^2\, b_{33}^2, \\
	\mathcal{E}_3 & = \tfrac{\mathcal{E}_2}{3} - \big(\tfrac{\mathcal{E}_1}{3})^2, \\
	\mathcal{E}_4 & = \big(\tfrac{\mathcal{E}_1}{3})^3 - \tfrac{\mathcal{E}_1 \, \mathcal{E}_2}{3} + \tfrac{1}{2} \, b_{11}^2 \, b_{22}^2 \, b_{33}^2, \\
	\mathcal{E}_5 & = \mathcal{E}_4 + (\mathcal{E}_3^3 + \mathcal{E}_4^2)^{\rfrac{1}{2}}.
\end{alignat*}
Therefore, $\cpalphahat$ and $\cgalphahat$ in 
\eqref{eq:general-poincare-type-inequalities-for-tetrahedron} are as follows:
\begin{equation}
  \cpalphahat = \mu^{\rfrac{1}{2}}_{\pi/2, \hat{\alpha}}, \quad
  \cgalphahat = 
	\Big(  \tfrac{\sin \hat{\alpha}}{\rho \sin\alpha \sin\theta}\Big)^{\rfrac{1}{2}} \,\cpalphahat.
	\label{eq:ratio-3d}
\end{equation}

Lower bounds of the constants $\CPGamma$ and
$\CtrGamma$ are computed by minimization of
$\mathcal{R}^{\mathrm{P}}_{\Gamma} [w]$ and $\mathcal{R}^{\mathrm{Tr}}_{\Gamma} [w]$
over the set    $V^N_3 \subset \H{1}(\T)$, where
\begin{equation*}
	V^{N}_3 := \Big\{\:\varphi_{ijk} = x^{i} y^{j} z^k, \quad i, j, k = 0, \ldots, N,
	\;\; (i, j, k) \neq (0, 0, 0) \: \Big\}
	%\label{eq:power-series-3d}
\end{equation*}
and $\mathrm{dim} V^{N}_3 = M(N) := (N + 1)^3 - 1$.
%--------------------------------------------------------------------------------------%
%Here, $V^{N}_3$ is finite dimensional subspace of $\H{1}(\T)$ formed with the help of%
%suitable anzats of trial functions in $\Rthree$, and where
%$\mathrm{dim} V^{N}_3 = M := (N + 1)^3 - 1$.

\begin{figure}[!ht]
	\centering
	\stackunder[1pt]{
	\includegraphics[scale=0.8]{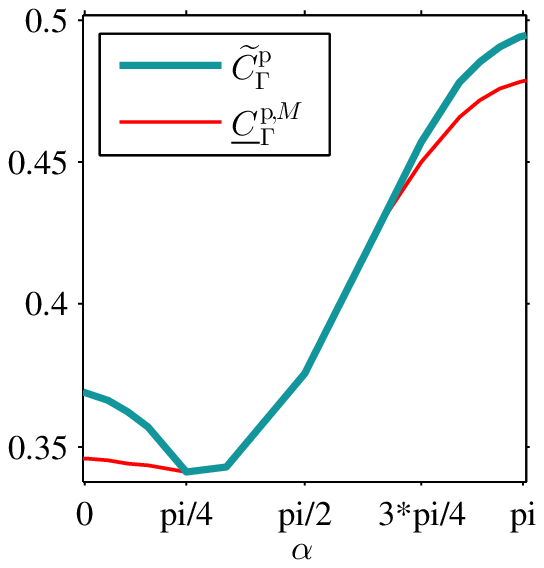}	}{}\quad
	\stackunder[1pt]{
  \includegraphics[scale=0.8]{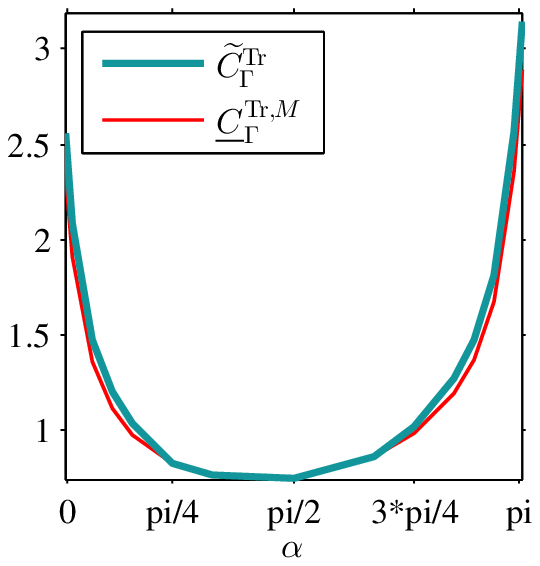}}{}\\  
	\caption{$\CPT$ and $\CGT$ for $\T \in \Rthree$ with $H = 1$, $\rho = 1$ with estimate based on 
	four reference tetrahedrons.}
	\label{fig:3d-cpt-cgt-constants-with-estimates-4-refs}
\end{figure}

%The numerical results with 
%$\underline{C}^{M, \mathrm{p}}_{\Gamma}$ and $\underline{C}^{M, \mathrm{Tr}}_{\Gamma}$ as well as $\%CPGamma$ and $\CtrGamma$ are presented with respect to the angle 
%$\alpha \in (0, \pi)$ and $\theta \in (0, \pi)$
%in Tables \ref{tab:cp-3d-approx-exact-constants} and \ref{tab:cg-3d-approx-exact-constants},
%respecxvely. One can generalize the upper bounds based on even more reference tetrahedrons, which %%would improve the estimates.
%--------------------------------------------------------------------------------------%
%
The respective results are presented
in Tables \ref{tab:cp-3d-approx-exact-constants} and \ref{tab:cg-3d-approx-exact-constants}
for $\T$ with $h_1 = 1$, $h_3 = 1$, and $\rho = 1$.
We note that exact values of constants are probably closer to the numbers presented in 
left-hand side columns. For $\theta = \rfrac{\pi}{2}$, we also present 
estimates of $\approxCPT$ and $\approxCGT$ (red lines) graphically in
Fig. \ref{fig:3d-cpt-cgt-constants-with-estimates-4-refs}.

\begin{table}[!ht]
\centering
\footnotesize
\begin{tabular}{c|cc|cc|cc|cc}
$$
& \multicolumn{2}{c|}{ $\alpha = \tfrac{\pi}{6}$}
& \multicolumn{2}{c|}{ $\alpha = \tfrac{\pi}{4}$}
& \multicolumn{2}{c|}{ $\alpha = \tfrac{\pi}{3}$}
& \multicolumn{2}{c}{ $\alpha = \tfrac{\pi}{2}$} \\
\midrule
$\theta$
    & $\underline{C}^{M, \mathrm{P}}_{\Gamma}$ & $\CPtetr$ &
    $\underline{C}^{M, \mathrm{P}}_{\Gamma}$ & $\CPtetr$ &
    $\underline{C}^{M, \mathrm{P}}_{\Gamma}$ & $\CPtetr$ &
        $\underline{C}^{M, \mathrm{P}}_{\Gamma}$ & $\CPtetr$ \\
\midrule
$\pi/6$ & 0.23883 & 0.49035& 0.24621 & 0.49841& 0.25870 & 0.51054& 0.29484 & 0.51308\\
$\pi/4$ & 0.23883 & 0.45388& 0.24621 & 0.46173& 0.25870 & 0.47683& 0.29484 & 0.49075\\
$\pi/3$ & 0.29666 & 0.41958& 0.31194 & 0.42259& 0.33489 & 0.43724& 0.38976 & 0.46002\\
$\pi/2$ & 0.34302 & 0.35667& 0.34112 & 0.34115& 0.34256 & 0.34259& 0.37559 & 0.37560\\
$2\pi/3$ & 0.40428 & 0.41958& 0.40562 & 0.42259& 0.40927 & 0.43724& 0.42867 & 0.46002\\
$3\pi/4$ & 0.42890 & 0.45388& 0.43110 & 0.46173& 0.43505 & 0.47683& 0.45017 & 0.49075\\
$5\pi/6$ & 0.44964 & 0.49035& 0.45193 & 0.49841& 0.45539 & 0.51054& 0.46607 & 0.51308\\
\midrule
$$
& \multicolumn{2}{c|}{ $\alpha = \tfrac{\pi}{2}$}
& \multicolumn{2}{c|}{ $\alpha = \tfrac{2\pi}{3}$}
& \multicolumn{2}{c|}{ $\alpha = \tfrac{3\pi}{4}$}
& \multicolumn{2}{c}{ $\alpha = \tfrac{5\pi}{6}$} \\
\midrule
$\theta$
    & $\underline{C}^{M, \mathrm{P}}_{\Gamma}$ & $\CPtetr$ &
    $\underline{C}^{M, \mathrm{P}}_{\Gamma}$ & $\CPtetr$ &
    $\underline{C}^{M, \mathrm{P}}_{\Gamma}$ & $\CPtetr$ &
        $\underline{C}^{M, \mathrm{P}}_{\Gamma}$ & $\CPtetr$ \\
\midrule
$\pi/6$ & 0.29484 & 0.51308& 0.33069 & 0.51792& 0.34468 & 0.52253& 0.35499 & 0.52694\\
$\pi/4$ & 0.29484 & 0.49075& 0.33069 & 0.50261& 0.34468 & 0.51308& 0.35499 & 0.52253\\
$\pi/3$ & 0.38976 & 0.46002& 0.43880 & 0.48413& 0.45742 & 0.50261& 0.47106 & 0.51792\\
$\pi/2$ & 0.37559 & 0.37560& 0.42865 & 0.42867& 0.45017 & 0.45731& 0.46607 & 0.47811\\
$2\pi/3$ & 0.42867 & 0.46002& 0.45997 & 0.48413& 0.47457 & 0.50261& 0.48598 & 0.51792\\
$3\pi/4$ & 0.45017 & 0.49075& 0.47204 & 0.50261& 0.48239 & 0.51308& 0.49064 & 0.52253\\
$5\pi/6$ & 0.46607 & 0.51308& 0.47972 & 0.51792& 0.48607 & 0.52253& 0.49115 & 0.52694\\
\end{tabular}
\\[5pt]
\caption{$\underline{C}^{M, \mathrm{P}}_{\Gamma}$ ($M(N)$ = 124) and $\CPtetr$.}
\label{tab:cp-3d-approx-exact-constants}
\end{table}
%--------------------------------------------------------------------------------------%
%--------------------------------------------------------------------------------------%

%\begin{figure}[!ht]
	%\centering
	%\stackunder[5pt]{\includegraphics[scale=0.83]{pics/3d-exact-and-approx-C-PT-H-0-5-rho-1-h-1-1-power-series}}{}\\ 
	%\stackunder[5pt]{\includegraphics[scale=0.83]{pics/3d-exact-and-approx-C-GT-H-0-5-rho-1-h-1-1-power-series}}{}\\ 
	%\caption{$\approxCPT$ and $\approxCGT$ for $\T \in \Rthree$ with $H = 0.5$, $\rho = 1$ and estimate 
	%based on two reference tetrahedrons.}
	%\label{fig:3d-cpt-cgt-H-1-2-constants-with-estimates-2-refs}
%\end{figure}
%%
%\begin{figure}[!ht]
	%\centering
	%\stackunder[1pt]{
  %\includegraphics[scale=0.83]{pics/3d-exact-and-approx-C-PT-H-1-rho-1-h-1-1-power-series}
	%}{}\\ 
	%\stackunder[1pt]{
  %\includegraphics[scale=0.83]{pics/3d-exact-and-approx-C-GT-H-1-rho-1-h-1-1-power-series}}{}\\ 
	%\caption{$\approxCPT$ and $\approxCGT$ for $\T \in \Rthree$ with $H = 1$, $\rho = 1$ and estimate 
	%based on two reference tetrahedrons..}
	%\label{fig:3d-cpt-cgt-H-1-constants-with-estimates-2-refs}
%\end{figure}
%

\begin{table}[!ht]
\centering
\footnotesize
\begin{tabular}{c|cc|cc|cc|cc}
$$
& \multicolumn{2}{c|}{ $\alpha = \tfrac{\pi}{6}$}
& \multicolumn{2}{c|}{ $\alpha = \tfrac{\pi}{4}$}
& \multicolumn{2}{c|}{ $\alpha = \tfrac{\pi}{3}$}
& \multicolumn{2}{c}{ $\alpha = \tfrac{\pi}{2}$}\\
\midrule
$\theta$    &
    $\underline{C}^{M, \mathrm{Tr}}_{\Gamma}$ & $\CGtetr$ &
    $\underline{C}^{M, \mathrm{Tr}}_{\Gamma}$ & $\CGtetr$ &
    $\underline{C}^{M, \mathrm{Tr}}_{\Gamma}$ & $\CGtetr$ &
        $\underline{C}^{M, \mathrm{Tr}}_{\Gamma}$ & $\CGtetr$ \\
\midrule
$\pi/6$ & 1.09760 & 3.78259& 0.96245 & 2.71866& 0.91255 & 2.27382& 0.93123 & 2.05449\\
$\pi/4$ & 1.09760 & 2.43897& 0.96245 & 1.78094& 0.91255 & 1.50166& 0.93123 & 1.38951\\
$\pi/3$ & 0.89122 & 1.74467& 0.79146 & 1.31130& 0.75950 & 1.12431& 0.78904 & 1.06349\\
$\pi/2$ & 0.98017 & 1.22920& 0.83132 & 0.83133& 0.76290 & 0.76291& 0.75199 & 0.75200\\
$2\pi/3$ & 1.17698 & 1.74467& 0.99473 & 1.31130& 0.90578 & 1.12431& 0.86463 & 1.06349\\
$3\pi/4$ & 1.35195 & 2.43897& 1.14144 & 1.78094& 1.03737 & 1.50166& 0.98220 & 1.38951\\
$5\pi/6$ & 1.65317 & 3.78259& 1.39424 & 2.71866& 1.26490 & 2.27382& 1.19017 & 2.05449\\
\midrule
$$
& \multicolumn{2}{c|}{ $\alpha = \tfrac{\pi}{2}$}
& \multicolumn{2}{c|}{ $\alpha = \tfrac{2\pi}{3}$}
& \multicolumn{2}{c|}{ $\alpha = \tfrac{3\pi}{4}$}
& \multicolumn{2}{c}{ $\alpha = \tfrac{5\pi}{6}$} \\
\midrule
$\theta$    &
    $\underline{C}^{M, \mathrm{Tr}}_{\Gamma}$ & $\CGtetr$ &
    $\underline{C}^{M, \mathrm{Tr}}_{\Gamma}$ & $\CGtetr$ &
    $\underline{C}^{M, \mathrm{Tr}}_{\Gamma}$ & $\CGtetr$ &
        $\underline{C}^{M, \mathrm{Tr}}_{\Gamma}$ & $\CGtetr$ \\
\midrule
$\pi/6$ & 0.93123 & 2.05449& 1.07244 & 2.39471& 1.21573 & 2.95902& 1.47044 & 4.21999\\
$\pi/4$ & 0.93123 & 1.38951& 1.07244 & 1.64324& 1.21573 & 2.01841& 1.47044 & 2.80588\\
$\pi/3$ & 0.78904 & 1.06349& 0.91773 & 1.27423& 1.04309 & 1.50833& 1.26357 & 2.11790\\
$\pi/2$ & 0.75199 & 0.75200& 0.86459 & 0.86463& 0.98220 & 1.12971& 1.19017 & 1.67033\\
$2\pi/3$ & 0.86463 & 1.06349& 0.96174 & 1.27423& 1.08134 & 1.50833& 1.30191 & 2.11790\\
$3\pi/4$ & 0.98220 & 1.38951& 1.07921 & 1.64324& 1.20686 & 2.01841& 1.44721 & 2.80588\\
$5\pi/6$ & 1.19017 & 2.05449& 1.29582 & 2.39471& 1.44268 & 2.95902& 1.72383 & 4.21999\\
\end{tabular}
\\[5pt]
\caption{$\underline{C}^{M, \mathrm{Tr}}_{\Gamma}$ ($M(N)$ = 124)
and $\CGtetr$ for different $\theta, \alpha \in (0, \pi)$.}
\label{tab:cg-3d-approx-exact-constants}
\end{table}
\section{Example}
\label{sec:example}

%--------------------------------------------------------------------------------------%
Constants in the Friedrichs', Poincar\'{e}, and other functional inequalities  arise in various problems of numerical analysis, where we need to 
know values of the respective constants associated with particular domains. 
Constants in projection type estimates arise in a priori analysis 
(see, e.g., \cite{Braess2001, Ciarlet1978, Mikhlin1986}). Constants in Clement's 
interpolation inequalities are important for residual type a posteriori estimates 
(see, e.g., \cite{AinsworthOden2000, Verfurth1996}, and \cite{CarstensenFunken2000}, 
where these constants have been evaluated).
%For example, 
%results related to extension and projection type estimates for FEM can be found in 
%\cite{Mikhlin1986,Ciarlet1978} (and 
%many other publications). 
Concerning constants in the trace inequalities associated 
with polygonal domain, we mention the paper \cite{CarstensenSauter2004}. 
Constants in functional (embedding) inequalities arise in a posteriori error estimates
of the functional type (error majorants). The details concerning last application can be 
found \cite{RepinDeGruyter2008, LangerRepinWolfmayr2014, 
MatculevichNeitaanmakiRepin2015,ReVarInq2010,RepinBoundaryMeanTrace2015, Repin2000, 
RepinXanthis1996} and other references cited therein.
%various error estimates. 
Below, we deduce an advanced version of an error majorant,
%in \cite{RepinDeGruyter2008}
which uses constants in Poincar\'{e}-type inequalities 
for functions with zero mean traces on inter-element boundaries.
This is done in order to maximally extend the space of admissible fluxes. 
However, first, we shall discuss the reasons that invoke 
Poincar\'{e}-type constants in a posteriori estimates.  
%Below we address 
%the latter case and first explain reasons that invoke the constants in general terms.

Let $u$ denote the exact solution of an elliptic boundary value problem generated by the pair
of conjugate operators $\rm grad$ and $-\dvrg$ 
(e.g.,  the problem (\ref{eq:reactdiff1})--(\ref{eq:reactdiff4}) considered below)
and $v$ be a function in the energy space satisfying the prescribed (Dirichlet)
boundary conditions. Typically, the error $e := u - v$ is measured in terms of the energy
norm $\|\nabla\,e\|$ (or some other equivalent norm), whose square is bounded
from above by the quantities
\begin{equation*}
\IntO R (v, \dvrg q) \,e\, \dx,\quad
\IntO D (\nabla v, q) \cdot \nabla\,e\, \dx,\;
{\rm and}\;
\Int_{\Gamma_N} R_{\Gamma_N}(v, q\cdot n)\,e\,\ds,
\end{equation*}
%
%--------------------------------------------------------------------------------------%
where $\Omega$ Lipschitz bounded domain, $\Gamma_N$ is the Neumann part of the boundary $\partial \Omega$ with the
outward unit normal vector $n$, and $q$ is an approximation of the dual variable (flux). 
The terms
$R$, $D$, and $R_{\Gamma_N}$ represent residuals of the differential (balance)
equation, constitutive (duality) relation, and Neumann boundary condition, respectively.
Since $v$ and $q$ are known from a numerical solution, fully computable estimates
can be obtained if these integrals are estimated by the H\"{o}lder, Friedrichs, and trace
inequalities (which involve the corresponding constants). However, for $\Omega$ with 
piecewise smooth (e.g., polynomial) boundaries these
constants may be unknown. A way to avoid these difficulties is suggested by
modifications of the estimates using ideas of domain decomposition. Assume that $\Omega$
is a polygonal (polyhedral) domain decomposed into a collection of non-overlapping convex
polygonal sub-domains $\Omega_i$, i.e.,
\begin{equation*}
    \overline{\Omega} :=
    \bigcup\limits_{ \Omega_i \in \, \mathcal{O}_\Omega}
    {\overline{\Omega}}_i,
    \quad
    \mathcal{O}_\Omega :=
    \Big\{ \; \Omega_i \in \Omega \; \big| \;
                        {\Omega}_{i'} \, \cap \, {\Omega}_{i''} = \emptyset, \;
                                                                                                    i' \neq i'', \;
                                                                                                    i = 1, \ldots, N \; \Big \}.
    %\label{eq:omega-representation}
\end{equation*}
%
%--------------------------------------------------------------------------------------%
We denote the set of all edges (faces) by ${\mathcal G}$ and the set of all interior
faces by ${\mathcal G}_{\rm int}$ (i.e., $\Gamma_{ij} \in {\mathcal G}_{\rm int}$,
if $\Gamma_{ij} = \overline{\Omega}_i \, \cap \, \overline{\Omega}_j$). Analogously,
${\mathcal G}_{N}$ denotes the set of edges on $\Gamma_N$. The latter set is decomposed
into $\Gamma_{N_k}:=\Gamma_N\cap\partial\Omega_k$ (the number of faces that belongs to 
${\mathcal G}_{N}$ is $K_N$).
Now, the integrals associated with $R $ and $R_{\Gamma_N}$ can be replaced by sums of
local quantities
\begin{equation*}
\Sum_{\Omega_i \in {\mathcal G}} \Int_{\Omega_i}R_{\Omega} (v,\dvrg q)\,e \, \dx,\quad{\rm and}\quad
\Sum_{\Gamma_{N_k} \in {\mathcal G}_{N}} \Int_{\Gamma_{N_k}} R_{\Gamma_{N}} (v,q\cdot n)\,e\,\ds.
\end{equation*}
If the residuals satisfy the conditions
\begin{equation*}
		\Int_{\Omega_i} R_{\Omega_i}(v,\, \dvrg q) \dx = 0
    \quad \forall i = 1, \ldots, N, \quad 
    %\label{eq:balanceOmega}
\end{equation*}
and
\begin{equation*}
		\Int_{\Gamma_{N_k}} R_{\Gamma_{N}}(v,\, q\cdot n) \ds = 0,
    \quad \forall k = 1, \ldots, K_N,
    %\label{eq:balanceGammaN}
\end{equation*}
then
\begin{equation}
    \Int_{\Omega_i} R_{\Omega}(v, \dvrg q)\,e\, \dx
    \leq C^{{\mathrm P}}_{\Omega_i} \|R_{\Omega_i}(v,\dvrg q)\|_{\Omega_i}\,
             \|\nabla\,e\|_{\Omega_i} \quad
    \label{eq:estOmega}
\end{equation}
and 
\begin{equation}
    \Int_{\Gamma_{N_k}} R_{\Gamma_{N}}(v, q\cdot n)\,e\,\ds
    \leq C^{{\mathrm Tr}}_{\Gamma_{N_k}} \|R_{\Gamma_{N}}(v,q\cdot n)\|_{\Gamma_{N_k}}
             \|\nabla\,e\|_{\Omega_k}.
    \label{eq:estGammaN}
\end{equation}
%
%--------------------------------------------------------------------------------------%
Hence, we can deduce a computable upper bound of the error that contains local constants
$C^{{\mathrm P}}_{\Omega_i}$ and $C^{{\mathrm Tr}}_{\Gamma_{N_k}}$ for simple subdomains
(e.g., triangles or tetrahedrons) instead of the global constants associated with
$\Omega$. 

The constant $\CPoincare$ may arise if, e.g., nonconforming approximations are used.
For example, if $v$ does not exactly satisfy the Dirichlet boundary condition on
$\Gamma_{D_k}$, then in the process of estimation it may be necessary to evaluate terms of
the type
\begin{equation*}
    \Int_{\Gamma_{D_k}} G_{D} (v)\,e\,\ds,\quad k = 1, \ldots, K_D,
\end{equation*}
%
%--------------------------------------------------------------------------------------%
where $\Gamma_{D_k}$ is a part of $\Gamma_D$ associated with a certain $\Omega_k$, and
$G_D(v)$ is a residual generated by inexact satisfaction of the boundary condition. If
we impose the requirement that the Dirichlet boundary condition is satisfied in a
weak sense, i.e.,
$
\mean{G_{D}(v)}_{\Gamma_{D_k}} = 0,
$
then each boundary integral can be estimated as follows:
\begin{equation}
    \Int_{\Gamma_{D_k}} G_{D}(v)\,e\,\ds
    \leq C^{{\mathrm P}}_{\Gamma_{D_k}} \|G_D(v)\|_{\Gamma_{D_k}}\|\nabla\,e\|_{\Omega_k}.
    \label{eq:estGammaD}
\end{equation}
After summing up (\ref{eq:estOmega}), (\ref{eq:estGammaN}), and (\ref{eq:estGammaD}), we
obtain a product of weighted norms of localized residuals (which are known) and
$\|\nabla e\|_\Omega$. Since the sum is bounded from below by the squared energy norm,
we arrive at computable error majorant.
%--------------------------------------------------------------------------------------%

Now, we discuss elaborately these questions within the paradigm of the 
following
boundary value problem: find $u$ such that
\begin{eqnarray}
& - \dvrg {p} + \varrho^2 u  = f, \quad & {\rm in} \; \Omega, \quad
\label{eq:reactdiff1}\\
& p = A \nabla  u, \quad & {\rm in} \; \Omega, \quad
\label{eq:reactdiff2}\\
& u = u_D,         \quad & {\rm on} \; \Gamma_D,\
%\label{eq:reactdiff3}
\\
& A \nabla u \cdot { n} = F \quad & {\rm on} \; \Gamma_N.
\label{eq:reactdiff4}
\end{eqnarray}
%
%--------------------------------------------------------------------------------------%
Here $f \in \L{2}(\Omega)$, $F \in \L{2}(\Gamma_N)$, $u_D \in \H{1}(\Omega)$, and $A$
is a symmetric positive definite matrix with bounded coefficients satisfying the
condition $\lambda_1 |\xi|^2 \, \leq \, A \xi \cdot \xi$, where $\lambda_1$ is a
positive constant independent of $\xi$. The generalized solution of
\eqref{eq:reactdiff1}--\eqref{eq:reactdiff4} exists and is unique in the set $V_0 + u_D$,
where 
$V_0 := \Big \{ w \in \H{1}(\Omega) \; \mid \; w = 0 \; {\rm on} \; \Gamma_D \Big\}$.

%--------------------------------------------------------------------------------------%
Assume that $v \in V_0 + u_D$ is a conforming approximation of $u$.
We wish to find a computable majorant of the error norm
\begin{equation}
\mid\!\mid\!\mid\! e \!\mid\!\mid\!\mid^2 \,:= \| \nabla e \|^2_{A} + \| \varrho\, e \|^2,
\end{equation}
where
$
\| \nabla e \|^2_{A} := \Int_\Omega A \nabla e \cdot \nabla e \dx.
$
First, we note that the integral identity that defines $u$ can be rewritten in the form
\begin{equation}
    \Int_\Omega A \nabla e \cdot \nabla w \, \dx \,
    + \Int_\Omega \varrho^2 e \, w \, \dx
    = \Int_\Omega (f w - \varrho^2 v \, w - A \nabla v \cdot \nabla w) \dx \,
    + \Int_{\Gamma_N} F w \ds, \quad \forall w \in V_0.
    \label{eq:weak-statement}
\end{equation}
%
%--------------------------------------------------------------------------------------%
It is well known (see Section 4.2 in \cite{RepinDeGruyter2008}) that this relation 
yields a computable majorant of $\mid\!\mid\!\mid\! e \!\mid\!\mid\!\mid^2 $, 
if we introduce a vector-valued function
${q} \in H(\Omega, \dvrg)$, such that ${q} \cdot n \in L^2(\Omega)$, and transform 
\eqref{eq:weak-statement} by means of integration by parts relations. The majorant has 
the form
\begin{equation}
    \mid\!\mid\!\mid\! e \!\mid\!\mid\!\mid\,
    \leq\,\|D_\Omega (\nabla v,q)\|_{A^{-1}} +
    C_1 \| R (v,\dvrg {q}) \|_{\Omega} + C_2 \| R_{\Gamma_N} (v,q\cdot n) \|_{\Gamma_N},
    \label{eq:global-estimate}
\end{equation}
where $C_1$ and $C_2$ are positive constants explicitly defined by $\lambda_1$, 
the Friedrichs' constant $C^{\mathrm F}_{\Omega}$ in inequality \linebreak
$\|v\|_\Omega \leq \, C^{\mathrm F}_{\Omega} \|\nabla v\|_\Omega$ for functions vanishing
on $\Gamma_D$, and constant $C^{\mathrm Tr}_{\Gamma_N}$ in the
trace inequality associated with $\Gamma_N$. The integrands are defined by the
relations
%
%--------------------------------------------------------------------------------------%
\begin{equation*}
    D (\nabla v,q) := A\nabla v - q, \quad
    R (v,\dvrg q) := \dvrg q + f - \varrho^2 v, \quad {\rm and} \quad
    R_{\Gamma_N} (v, q \cdot n) := q\cdot n - F.
\end{equation*}
In general, finding $C^{\mathrm F}_{\Omega}$ and $C^{\mathrm Tr}_{\Gamma_N}$ may not be 
an easy
task. We can exclude $C_2$ if $ q$ additionally satisfies the condition ${ q}\cdot n = F$.
Then, the last term in \eqref{eq:global-estimate} vanishes. However, this condition is
difficult to satisfy, if $F$ is
a complicated nonlinear function. In order to exclude $C_1$ together with $C_2$, 
we can apply domain
decomposition technique and use (\ref{eq:estOmega}) instead of the global estimate. Then, the
estimate will operate with the constants $C^{{\mathrm P}}_{\Omega_i}$ (whose upper bounds
are known for convex domains). Moreover, it is shown below that by using the inequalities
(\ref{eq:Comega}) and (\ref{eq:Cgamma}), we can essentially weaken the assumptions
required for the variable $q$.

%--------------------------------------------------------------------------------------%
Define the space of vector-valued functions
\begin{alignat*}{2}
 \hat{H} (\Omega, {\mathcal O}_{\Omega}, \dvrg) := \, \Big\{{q} \in \L{2} (\Omega, \Rd) \; \mid & \quad {q} = {q}_i \in H(\Omega_i, \dvrg),
                    \;\; \mean{\dvrg {q}_i + f - \varrho^2 \, v}_{\Omega_i} = 0, \quad \forall\,\Omega_i \in {\mathcal O}_{\Omega}, \\
										& \;\; \mean{({q}_i - {q}_j) \cdot { n}_{ij}}_{\Gamma_{ij}} = 0 , \quad \forall\, \Gamma_{ij} \in {\mathcal G}_{\rm int},\\
                    & \;\; \mean{{ q}_i \cdot { n}_k - F}_{ \Gamma_{N_k}} = 0, \quad \forall \, k = 1, \ldots, K_N \Big \}.
\end{alignat*}
We note that the space $\hat{H} (\Omega, {\mathcal O}_{\Omega},\dvrg)$ is wider than
$H(\Omega, \dvrg)$ (so that we have more flexibility in determination of 
optimal reconstruction
of numerical fluxes). Indeed, the vector-valued functions in $H(\Omega,\dvrg)$ must
have continuous normal components on all $\Gamma_{ij} \in {\mathcal G}_{\rm int}$ 
and satisfy the Neumann boundary
condition in the pointwise sense. The functions in
$\hat{H} (\Omega, {\mathcal O}_{\Omega},\dvrg)$ satisfy much weaker conditions: namely,
the normal components are continuous  only in terms of mean values (integrals) and
the Neumann condition must hold in the integral sense only.

%--------------------------------------------------------------------------------------%
We reform (\ref{eq:weak-statement}) by means of the integral identity
\begin{equation*}
    \Sum_{\Omega_i \in {\mathcal O}_{\Omega}} \int_{\Omega_i}
    \left( {q} \cdot \nabla w + \dvrg {q} \, w \right) \dx
    = \Sum_{\Gamma_{ij} \, \in \, {\mathcal G}_{\rm int}} \;
        \Int_{\Gamma_{ij}}\; ({ q}_i - { q}_j) \cdot { n}_{ij} \, w \ds
        + \, \Sum_{\Gamma_{N_k} \, \in \, \Gamma_N} \;
        \Int_{\Gamma_{N_k}} { q}_i \cdot { n}_{i} \, w \ds,
\end{equation*}
which holds for any $ w \in V_0$ and
${q} \in \hat{H} (\Omega, {\mathcal O}_{\Omega},\dvrg)$. By setting $w = e$ in
\eqref{eq:weak-statement} and applying the H\"{o}lder inequality, we find that
%
%--------------------------------------------------------------------------------------%
\begin{multline*}
    \mid\!\mid\!\mid\! e \!\mid\!\mid\!\mid^2\! \;\;
    \leq \|D (\nabla v,q) \|_{A^{-1}} \|\nabla e\|_{A}
    + \Sum_{\Omega_i \in {\mathcal O}_{\Omega}} \|R(v,\dvrg q)\|_{\Omega_i}
                             \big\|e - \mean{e}_{\Omega_i}\big\|_{\Omega_i}   \\
    + \Sum_{\Gamma_{ij} \in \, {\mathcal G}_{\rm int}}
                      r_{ij}(q) \big\| e - \mean{e}_{\Gamma_{ij}}\big\|_{\Gamma_{ij} }
    + \Sum_{\Gamma_{N_k} \in \Gamma_N} \rho_k(q)
                         \Big\| e - \mean{e}_{\Gamma_{N_k}}\Big\|_{\Gamma_{N_k}},
   % \label{eq:estimate-1}
\end{multline*}
where
\begin{equation*}
    r_{ij}(q) := \| ({ q}_i - { q}_j) \cdot { n}_{ij} \|_{\Gamma_{ij}} 
		\quad \mbox{and} \quad
    \rho_k(q) := \| { q}_k \cdot { n}_k - F \|_{\Gamma_{N_k}}.
\end{equation*}
%
%--------------------------------------------------------------------------------------%
In view of (\ref{eq:classical-poincare-constant}) and (\ref{eq:Cgamma}), we obtain
\begin{multline}
    \mid\!\mid\!\mid\! e \!\mid\!\mid\!\mid^2\! \;\;
    \leq \|D (\nabla v,q)  \|_{A^{-1}} \|\nabla e\|_{A}
    + \Sum_{\Omega_i \in {\mathcal O}_{\Omega}}
      \|R(v,\dvrg q)\|_{\Omega_i} C^{\mathrm P}_{\Omega_i}
                         \|\nabla e \|_{\Omega_i}   \\
    + \Sum_{\Gamma_{ij} \in \, {\mathcal G}_{\rm int}}
                      r_{ij}(q) C^{\mathrm{Tr}}_{\Gamma_{ij}}
                                            \|\nabla e \|_{\Omega_i}
    + \Sum_{\Gamma_{N_k} \in \Gamma_N}
                         \rho_k(q) C^{\mathrm{Tr}}_{\Gamma_{N_k}}
                                             \|\nabla e \|_{\Omega_i}.
    \label{eq:estimate-2}
\end{multline}
%
%--------------------------------------------------------------------------------------%
The second term in the right hand side is estimated by the quantity
$\Re_1(v,q)\,\|\nabla e\|_\Omega$, where
\begin{equation*}
\Re^2_1(v,q) := \,\Sum_{\Omega_i \in {\mathcal O}_{\Omega}}
\tfrac{( \diam \, \Omega_i)^2}{\pi^2}\|R(v,\dvrg q)\|^2_{\Omega_i}.
\end{equation*}
We can represent any $\Omega_i \in {\mathcal O}_{\Omega}$ as a sum of simplexes such
that each simplex has one edge on $\partial\Omega_i$. Let $C^{\mathrm{Tr}}_{i, {\rm max}}$
denote the largest constant in the respective Poincar\'{e}-type inequalities
\eqref{eq:Cgamma} associated with all edges of $\partial\Omega_i$. Then, the last two
terms of (\ref{eq:estimate-2}) can be estimated by the quantity
$\Re_2 (v,q)\,\|\nabla e\|_\Omega$, where
\begin{equation*}
    \Re^2_2(q):=\,\Sum_{\Omega_i \in {\mathcal O}_{\Omega}}
    (C^{\mathrm{Tr}}_{i, {\rm max}})^2 \, \eta^2_i, \quad { \rm with} \quad
%\end{equation*}
%
%and
%
%\begin{equation*}
    \eta^2_{i} =
    \Sum_{\myatop{\Gamma_{ij} \in {\mathcal G}_{\rm int}}{\Gamma_{ij} \cap \partial\Omega_i \neq \varnothing}}\tfrac{1}{4} r^2_{ij}(q) +
    \Sum_{\myatop{\Gamma_{k} \in {\mathcal G}_{N}}{\Gamma_{k} \,\cap \,\partial\Omega_i \neq \varnothing}}
    \rho^2_k(q).
\end{equation*}
Then, (\ref{eq:estimate-2}) yields the estimate
\begin{equation*}
    \mid\!\mid\!\mid\! e \!\mid\!\mid\!\mid^2\! \;\;
    \leq \|D(\nabla v,q)\|_{A^{-1}} \|\nabla e\|_{A}
    + (\Re_1(v,q) + \Re_2(q) )\,\|\nabla e\|_\Omega,
    %\label{eq:estimate-3}
\end{equation*}
which shows that
%--------------------------------------------------------------------------------------%
\begin{equation}
    \mid\!\mid\!\mid\! e \!\mid\!\mid\!\mid\! \;\;
    \leq \|D(\nabla v,q) \|_{A^{-1}}
    + \tfrac{1}{\lambda_1} \Big( \Re_1(v,q) + \Re_2(q) \Big).
    \label{eq:estimate-final}
\end{equation}
%
%--------------------------------------------------------------------------------------%
Here, the term $\Re_2(q)$ controls violations of conformity of $q$ (on interior edges)
and inexact satisfaction of boundary conditions (on edges related to $\Gamma_N$). It is
easy to see that $\Re_2(q) = 0$, if and only if the quantity $q \cdot n$ is continuous on
${\mathcal G}_{\rm int}$ and exactly satisfies the boundary condition. Hence, $\Re_2(q)$ 
can be viewed as a measure of  the `flux nonconformity'. Other terms have the same
meaning as in well-known a posteriori estimates of the functional type, namely, the 
first term measures the
violations of the relation ${q} = A \nabla v$ (cf. (\ref{eq:reactdiff2})), and
$\Re_1(v,q)$  measures inaccuracy in the equilibrium (balance) equation
(\ref{eq:reactdiff1}). The right-hand side of \eqref{eq:estimate-final} contains known
functions (approximations $v$ and $q$ of the exact solution and exact flux). The 
constants $C^{\mathrm{Tr}}_{i, {\rm max}}$ can be easily computed using results of Section \ref{sc:arbitrary-triangle}-\ref{eq:numerial-tests-3d}.
Finally, we note that estimates similar to \eqref{eq:estimate-final} were derived in 
\cite{ReVarInq2010} for elliptic variational inequalities and in 
\cite{MatculevichNeitaanmakiRepin2015} for a class of parabolic problems.

%\begin{acknowledgement}
%This work was supported by the Suomalainen Tiedeakatemia foundation and 
%RFBR Grant $14\minus01\minus00162\backslash 15$.
%\end{acknowledgement}

\bibliographystyle{plain} 
\bibliography{lib,my_lib}

\def\cprime{$'$}
\begin{thebibliography}{10}

\bibitem{AinsworthOden2000}
M.~Ainsworth and J.~T. Oden.
\newblock {\em A posteriori error estimation in finite element analysis}.
\newblock Wiley and Sons, New York, 2000.

\bibitem{BabuskaAziz1976}
I.~Babu{\v{s}}ka and A.~K. Aziz.
\newblock On the angle condition in the finite element method.
\newblock {\em SIAM J. Numer. Anal.}, 13(2):214--226, 1976.

\bibitem{Bandle1980}
C.~Bandle.
\newblock {\em Isoperimetric inequalities and applications}, volume~7 of {\em
  Monographs and Studies in Mathematics}.
\newblock Pitman (Advanced Publishing Program), Boston, Mass.-London, 1980.

\bibitem{BanuelosKulczyckiPolterovichSiudeja2010}
R.~Ba{\~n}uelos, T.~Kulczycki, I.~Polterovich, and B.~Siudeja.
\newblock Eigenvalue inequalities for mixed {S}teklov problems.
\newblock In {\em Operator theory and its applications}, volume 231 of {\em
  Amer. Math. Soc. Transl. Ser. 2}, pages 19--34. Amer. Math. Soc., Providence,
  RI, 2010.

\bibitem{Berard1980}
P.~H. B{\'e}rard.
\newblock Spectres et groupes cristallographiques. {I}. {D}omaines euclidiens.
\newblock {\em Invent. Math.}, 58(2):179--199, 1980.

\bibitem{Braess2001}
D.~Braess.
\newblock {\em Finite elements}.
\newblock Cambridge University Press, Cambridge, second edition, 2001.
\newblock Theory, fast solvers, and applications in solid mechanics, Translated
  from the 1992 German edition by Larry L. Schumaker.

\bibitem{CarstensenFunken2000}
C.~Carstensen and S.~A. Funken.
\newblock Fully reliable localized error control in the {FEM}.
\newblock {\em SIAM J. Sci. Comput.}, 21(4):1465--1484, 2000.

\bibitem{CarstensenGedicke2014}
C.~Carstensen and J.~Gedicke.
\newblock Guaranteed lower bounds for eigenvalues.
\newblock {\em Math. Comp.}, 83(290):2605--2629, 2014.

\bibitem{CarstensenSauter2004}
C.~Carstensen and S.~A. Sauter.
\newblock A posteriori error analysis for elliptic {PDE}s on domains with
  complicated structures.
\newblock {\em Numer. Math.}, 96(4):691--721, 2004.

\bibitem{Cheng1975}
S.~Y. Cheng.
\newblock Eigenvalue comparison theorems and its geometric applications.
\newblock {\em Math. Z.}, 143(3):289--297, 1975.

\bibitem{Ciarlet1978}
P.~G. Ciarlet.
\newblock {\em The finite element method for elliptic problems}.
\newblock North-Holland Publishing Co., Amsterdam-New York-Oxford, 1978.
\newblock Studies in Mathematics and its Applications, Vol. 4.

\bibitem{Dohrmann2008}
C.~R. Dohrmann, A.~Klawonn, and O.~B. Widlund.
\newblock Domain decomposition for less regular subdomains: overlapping
  {S}chwarz in two dimensions.
\newblock {\em SIAM J. Numer. Anal.}, 46(4):2153--2168, 2008.

\bibitem{FoxKuttler1983}
D.~W. Fox and J.~R. Kuttler.
\newblock Sloshing frequencies.
\newblock {\em Z. Angew. Math. Phys.}, 34(5):668--696, 1983.

\bibitem{ArxivGirouardPolterovich}
A.~Girouard and I.~Polterovich.
\newblock Spectral geometry of the steklov problem.
\newblock {\em arXiv.org}, math/1411.6567, 2014.

\bibitem{HoshikawaUrakawa2010}
Y.~Hoshikawa and H.~Urakawa.
\newblock Affine {W}eyl groups and the boundary value eigenvalue problems of
  the {L}aplacian.
\newblock {\em Interdiscip. Inform. Sci.}, 16(1):93--109, 2010.

\bibitem{Klawonnatall2008}
A.~Klawonn, O.~Rheinbach, and O.~B. Widlund.
\newblock An analysis of a {FETI}-{DP} algorithm on irregular subdomains in the
  plane.
\newblock {\em SIAM J. Numer. Anal.}, 46(5):2484--2504, 2008.

\bibitem{KozlovKuznetsov2004}
V.~Kozlov and N.~Kuznetsov.
\newblock The ice-fishing problem: the fundamental sloshing frequency versus
  geometry of holes.
\newblock {\em Math. Methods Appl. Sci.}, 27(3):289--312, 2004.

\bibitem{KozlovKuznetsovMotygin2004}
V.~Kozlov, N.~Kuznetsov, and O.~Motygin.
\newblock On the two-dimensional sloshing problem.
\newblock {\em Proc. R. Soc. Lond. Ser. A Math. Phys. Eng. Sci.},
  460(2049):2587--2603, 2004.

\bibitem{KuznetsovKulczyckiKwasnickiNazarovPoborchiPolterovichSiudeja2014}
N.~Kuznetsov, T.~Kulczycki, M.~Kwa{\'s}nicki, A.~Nazarov, S.~Poborchi,
  I.~Polterovich, and B.~Siudeja.
\newblock The legacy of {V}ladimir {A}ndreevich {S}teklov.
\newblock {\em Notices Amer. Math. Soc.}, 61(1):9--22, 2014.

\bibitem{LangerRepinWolfmayr2014}
U.~Langer, S.~Repin, and M.~Wolfmayr.
\newblock Functional a posteriori error estimates for parabolic time-periodic
  boundary value problems.
\newblock {\em CMAM}, 15(3):353--372, 2015.

\bibitem{LaugesenSiudeja2009}
R.~S. Laugesen and B.~A. Siudeja.
\newblock Maximizing {N}eumann fundamental tones of triangles.
\newblock {\em J. Math. Phys.}, 50(11):112903, 18, 2009.

\bibitem{LaugesenSiudeja2010}
R.~S. Laugesen and B.~A. Siudeja.
\newblock Minimizing {N}eumann fundamental tones of triangles: an optimal
  {P}oincar\'e inequality.
\newblock {\em J. Differential Equations}, 249(1):118--135, 2010.

\bibitem{XuefengOishi2013}
X.~Liu and S.~Oishi.
\newblock Guaranteed high-precision estimation for {$P_0$} interpolation
  constants on triangular finite elements.
\newblock {\em Jpn. J. Ind. Appl. Math.}, 30(3):635--652, 2013.

\bibitem{LoggMardalWells2012}
A.~Logg, K.-A. Mardal, and G.~N. Wells, editors.
\newblock {\em Automated solution of differential equations by the finite
  element method}, volume~84 of {\em Lecture Notes in Computational Science and
  Engineering}.
\newblock Springer, Heidelberg, 2012.
\newblock The FEniCS book.

\bibitem{MatculevichNeitaanmakiRepin2015}
S.~Matculevich, P.~Neittaanm{\"a}ki, and S.~Repin.
\newblock A posteriori error estimates for time-dependent reaction-diffusion
  problems based on the {P}ayne--{W}einberger inequality.
\newblock {\em AIMS}, 35(6), 2015.

\bibitem{McCartin2002}
B.~J. McCartin.
\newblock Eigenstructure of the equilateral triangle. {II}. {T}he {N}eumann
  problem.
\newblock {\em Math. Probl. Eng.}, 8(6):517--539, 2002.

\bibitem{Mikhlin1986}
S.~G. Mikhlin.
\newblock {\em Constants in some inequalities of analysis}.
\newblock A Wiley-Interscience Publication. John Wiley and Sons, Ltd.,
  Chichester, 1986.
\newblock Translated from the Russian by Reinhard Lehmann.

\bibitem{NakaoYamamoto2001I}
M.~T. Nakao and N.~Yamamoto.
\newblock A guaranteed bound of the optimal constant in the error estimates for
  linear triangular element.
\newblock In {\em Topics in numerical analysis}, volume~15 of {\em Comput.
  Suppl.}, pages 165--173. Springer, Vienna, 2001.

\bibitem{NazarovRepin2014}
A.~I. Nazarov and S.~I. Repin.
\newblock Exact constants in {P}oincare type inequalities for functions with
  zero mean boundary traces.
\newblock {\em Mathematical Methods in the Applied Sciences}, 2014.
\newblock Puplished in arXiv.org in 2012, math/1211.2224.

\bibitem{PayneWeinberger1960}
L.~E. Payne and H.~F. Weinberger.
\newblock An optimal {P}oincar\'e inequality for convex domains.
\newblock {\em Arch. Rational Mech. Anal.}, 5:286--292 (1960), 1960.

\bibitem{Pinsky1980}
M.~A. Pinsky.
\newblock The eigenvalues of an equilateral triangle.
\newblock {\em SIAM J. Math. Anal.}, 11(5):819--827, 1980.

\bibitem{Poincare1890}
H.~Poincare.
\newblock Sur les {E}quations aux {D}erivees {P}artielles de la {P}hysique
  {M}athematique.
\newblock {\em Amer. J. Math.}, 12(3):211--294, 1890.

\bibitem{Poincare1894}
H.~Poincare.
\newblock Sur les {E}quations de la {P}hysique {M}athematique.
\newblock {\em Rend. Circ. Mat. Palermo}, 8:57--156, 1894.

\bibitem{Repin2000}
S.~Repin.
\newblock A posteriori error estimation for variational problems with uniformly
  convex functionals.
\newblock {\em Math. Comput.}, 69(230):481--500, 2000.

\bibitem{RepinDeGruyter2008}
S.~Repin.
\newblock {\em A posteriori estimates for partial differential equations},
  volume~4 of {\em Radon Series on Computational and Applied Mathematics}.
\newblock Walter de Gruyter GmbH \& Co. KG, Berlin, 2008.

\bibitem{ReVarInq2010}
S.~Repin.
\newblock Estimates of deviations from exact solutions of variational
  inequalities based upon payne-weinberger inequality.
\newblock {\em J. Math. Sci. (N. Y.)}, 157(6):874--884, 2009.

\bibitem{RepinBoundaryMeanTrace2015}
S.~Repin.
\newblock Estimates of constants in boundary-mean trace inequalities and
  applications to error analysis.
\newblock In {\em Numerical Mathematics and Advanced Applications - ENUMATH
  2013}, volume 103 of {\em Lecture Notes in Computational Science and
  Engineering}, pages 215--223. Springer, Switzerland, 2015.

\bibitem{RepinXanthis1996}
S.~I. Repin and L.~S. Xanthis.
\newblock A posteriori error estimation for elastoplastic problems based on
  duality theory.
\newblock {\em Comput. Methods Appl. Mech. Engrg.}, 138(1-4):317--339, 1996.

\bibitem{Stekloff1902}
V.~A. Steklov.
\newblock Sur les probl\'{e}mes fondamentaux de la physique mathematique.
\newblock {\em Ann. Sci. Ec. Norm. Sup\'{e}r}, 3(19):191--259, 455--490, 1902.

\bibitem{Matlab}
Inc. \textcopyright 1994-2015 The~MathWorks.
\newblock Mathworks. products and services, 2015.

\bibitem{ToselliWidlund2005}
A.~Toselli and O.~Widlund.
\newblock {\em Domain decomposition methods---algorithms and theory}, volume~34
  of {\em Springer Series in Computational Mathematics}.
\newblock Springer-Verlag, Berlin, 2005.

\bibitem{Verfurth1996}
R.~Verf{\"u}rth.
\newblock {\em A review of a posteriori error estimation and adaptive
  mesh-refinement techniques}.
\newblock Wiley and Sons, Teubner, New-York, 1996.

\end{thebibliography}

\end{document}